\documentclass[%
 notitlepage,superscriptaddress,
 amsmath,amssymb,
 aps
 amsmath,amssymb,
aps,pra,onecolumn
]{revtex4-2}

\usepackage{graphicx}
\usepackage{dcolumn}
\usepackage{bm}
\usepackage{mathtools}
\usepackage{tikz}
\usepackage{hyperref}
\usepackage{xcolor}

\begin{document}
\newcommand\rnumber{\operatorname{L_2-norm}}
\newcommand{\thickhyphen}{\rule[0.5ex]{0.5em}{1pt}}
\newcommand{\thickhyphent}{\rule[0.5ex]{0.35em}{1pt}}
\preprint{APS/123-QED}

\title{Compact LABFM: a framework for meshless methods with spectral-like resolving power
}

\author{H. M. Broadley}
\email{henry.broadley@manchester.ac.uk}

\author{J. R. C. King}%

\affiliation{%
 Department of Mechanical and Aerospace Engineering, The University of Manchester}%

\author{S. J. Lind}

\affiliation{School of Engineering, Cardiff University}%

\date{\today}

\begin{abstract}
Meshless methods are often used in numerical simulations of systems of partial differential equations (PDEs), particularly those which involve complex geometries or free surfaces. Here we present a novel compact scheme based on the local anisotropic basis function method (LABFM), a meshless method which provides approximations to spatial operators to arbitrary polynomial consistency. Our approach mimics compact finite-differences by using implicit stencils to optimise the resolving power of each operator, whilst retaining diagonal dominance of the resulting global sparse linear system. The new method is demonstrated to provide improved approximations by a series of convergence tests and resolving power analysis, before solutions to canonical PDEs are computed. Significant gains in accuracy are observed, in particular for solutions containing high wavenumber components. Our compact meshfree method provides a pathway to high-order simulations of PDEs in complex geometries with spectral-like resolving power, and has the potential to lead to a step-change in the accuracy of numerical solutions to such problems.

\end{abstract}

\maketitle

\section{Introduction}\addvspace{10pt}
\label{sec::Introduction}
Approaches to numerical simulations of partial differential equations (PDEs) typically take one of two forms when discretising spatial derivatives: explicit or implicit. Explicit discretisations require information from only current values of a given function to calculate the desired derivative, whereas implicit approximations require information about derivatives from other nodal points, thus a global linear system must be solved in order to obtain the approximations at all points simultaneously. Despite this extra computational cost, there are many scenarios where implicit approximations are preferable to explicit computations. For example, when solving Poisson-type problems which require a global linear system to be solved for both explicit and implicit discretisations, an implicit scheme is often preferable. In addition, implicit schemes often have better resolving power characteristics, a key metric which measures how well a numerical method resolves all wavelengths admitted by a given discretisation, than their explicit counterparts for a given order of polynomial consistency. This property has lead to common use of compact schemes in codes like Xcompact3D\cite{Xcompact3D} which simulate (for example) turbulent flows, where accurate resolution of the smallest scales is a necessity. 

In general, compact schemes can be separated into two different types: the Lele-type schemes of \cite{Lele}, and the Weinan-type schemes of \cite{Weinan}. Weinan-type schemes require information to be gleaned from the underlying PDE of interest to improve the convergence properties of approximation, enabling more accurate computations of (for example) Poisson equations. In contrast, Lele-type schemes develop implicit formulas for derivatives which often eschew the highest order polynomial consistency possible for a given stencil, instead preferring to maximise resolving power of a given operator. This approach allows for applicability to arbitrary operators, and each operator in a system of PDEs can be treated individually in this manner. For both types of schemes, in finite-difference settings analytic formulae for weights can be derived which are the same for each node, and the structured mesh allows for the construction of banded matrices that enable efficient solution of the global linear system. Thus, both approaches can achieve spectral-like accuracy at limited computational cost.

Finite-difference type schemes are by their nature limited to simple geometries. If we wish to solve PDEs in more complex geometries, meshless methods, which can discretise the domain trivially, are much more suitable. For an overview of meshless methods in general, we refer the reader to \cite{Li_Review,Garg_Review}. The most widely used meshless method is Smoothed Particle Hydrodynamics (SPH), an interpolation scheme originally developed for astrophysical applications \cite{gingold,Lucy}. Since the advent of modern high-performance computing, the highly scalable nature of SPH has led to its widespread use across the fluid dynamics community, in particular for free surface flows \cite{Monaghan_Review}. Traditional implementations of SPH for incompressible fluid dynamics problems used the weakly compressible formulation of the Navier--Stokes equations \cite{Monaghan_94}, a purely explicit approach to solving the governing equations. More recently, a truly incompressible formulation of SPH (ISPH) has been developed \cite{ISPH_method}, utilising a pressure-projection method and thus requiring the solution of a Poisson equation at each timestep. This requires the solution of a global, sparse linear system to obtain the pressure field. Given the large number of particles often associated with SPH simulation (often of the order $10^8$), this has led to development of massively parallel codes such as \cite{Xiaohu_MPSPH} which efficiently solve such systems using sparse iterative solvers.

The low-order convergence properties of SPH has led to the development of other, high-order meshless methods. The radial basis function (RBF) method, first used for cartography \cite{RBF_Cartography} is a global meshless method with spectral-like accuracy. However, the global nature of the RBF method limits its applicability, as weights are obtained by solving a dense global linear system which must remain well-conditioned. This restriction is difficult to satisfy on large grids and led to the development of the radial basis function-finite-differences method (RBF-FDs) \cite{shu,Cecil,Wright}, where a local stencil is used to approximate spatial operators. Consequently, only a localised linear system is required to be solved to generate weights for approximations, which are required to allow for exact reconstruction of derivatives of polynomials up to a given degree (typically second order polynomial consistency). The greater flexibility of RBF-FDs in allowing for numerical simulations of larger problems comes at the cost of the spectral resolution characteristics provided by global RBFs. 

More recently, compact RBF-FDs (often termed RBF-HFDs), an implicit method first developed by \cite{WrightFornberg_compactFDs}, has been used in solutions of various PDE systems (see \cite{Fornberg&Flyer2015} and citations therein for examples). This method uses Hermite RBF interpolation to obtain high-order approximation using the same size local stencils as RBF-FDs. Whilst these methods can provide highly accurate approximations, to our knowledge these have so far been limited to Weinan-type schemes and therefore must be tailored to the PDE system of interest. In addition, this approach typically precludes compact treatment of (for example) advection terms in the Navier--Stokes equations which must be treated explicitly.

Another high-order meshless method is the generalised moving least squares (GMLS) of \cite{TraskGMLS}. Here polynomial reconstruction is used to approximate operators, and when applied to Stokes flow in complex geometries up to $4^{th}$ order convergence was observed. More recently, a compact formulation of moving least squares (CMLS) has been developed \cite{TraskCMLS}. Therein both Weinan-type and Lele-type compact schemes were used, the first time a meshless version of Lele-type approximations had been generated. Similar Hermite interpolation techniques to RBF-HFDs were used to generate weights, and although accuracy gains of an order of magnitude were observed in convergence studies no comparable resolving power analysis to that found in \cite{Lele} was performed.

In this paper we take a different approach to generating compact approximations. We aim to reformulate the local anisotropic basis function method (LABFM), an explicit meshless method developed by some of the authors \cite{king_2020}, in a compact fashion. Explicit LABFM uses weighted sums of anisotropic basis functions (typically orthogonal polynomials multiplied by a smoothing kernel \cite{king_2022}) to generate approximations to differential operators, with weights obtained by solving a local linear system. Our goal is to provide Lele-type compact approximations of operators by maximising the resolving power of our approximations. Given the disordered nodal distributions inherent to meshless discretisations, no analytic formula for coefficients can be prescribed and these will differ from node to node. In addition, we have previously shown \cite{MKLABFM_Arxiv} that, unlike finite-differences, resolving power analysis of a meshless method in a two-dimensional domain has a dependency on both components of the wavenumber vector. Any attempt to maximise resolving power must therefore consider the whole of wavenumber space, not just the line corresponding to the direction of the derivative as in \cite{Lele}. Unlike other compact meshless methods, our criteria for the choice of implicit coefficients is based on optimisation of resolution characteristics of an approximation across all wavelengths permitted by a given discretisation. 

The rest of this paper is organised as follows. In \S\ref{sec:Explicit_LABFM} we briefly present explicit LABFM with the minimum description necessary for reproduction of the method. In \S\ref{sec:Compact} we formulate compact LABFM, describing how we choose the nodes in the implicit stencil for a given operator, and how the implicit coefficients are chosen to maximise resolving power whilst maintaining diagonal dominance of the global linear system. In \S\ref{sec:Resolving Power} the improvements in resolving power obtained using the implicit approximations are demonstrated. In \S\ref{sec:Convergence_Studies} we compare the convergence of compact operators on a test function and analyse their stability. In \S\ref{sec:PDE_solutions} we solve test PDE systems with the new formulation. For the Lele-type derivatives considered in this manuscript, all compact approximations require the solution of a global linear system. Such approximations are of course implicit, and we use the terms implicit and compact interchangeably throughout the manuscript.

\section{Explicit LABFM}
\label{sec:Explicit_LABFM}
Before introducing implicit LABFM approximations (the main focus of this work), we briefly outline the standard explicit implementation of LABFM as first presented in \cite{king_2020}. We refer to collocation points in each domain as nodes and define our local stencil at each node $i$ to be the set of nodes $\mathcal{N}_i$. Consider a set of $N$ nodes in a domain $\Omega\subset\mathbb{R}^2$. We define the position of each node $i$ through $\boldsymbol{r}_i=(x_i,y_i)$, and to approximate derivatives at $\boldsymbol{r}_i$ we use computational stencils of radius $2h_i$. All points $j$ within the computational stencil of the point $i$ are referred to as the neighbours of $i$, with the set of neighbours of $i$ denoted $\mathcal{N}_i$ and the average nodal spacing within the stencil of $i$ denoted $s_i$. For ease of exposition we denote $(\cdot)_{ji}=(\cdot)_j-(\cdot)_i$ throughout. We first define a vector of monomials
\begin{equation}
    \boldsymbol{X}=\left[x,\ y,\ \frac{x^2}{2},\ xy,\ \frac{y^2}{2},\ \frac{x^3}{6}, \cdots\right]^T
\end{equation}
and a vector operator of partial derivatives
\begin{equation}
    \boldsymbol{\partial}(\cdot) = \left[\frac{\partial (\cdot)}{\partial x},\ \frac{\partial (\cdot)}{\partial y},\ \frac{\partial^2 (\cdot)}{\partial x^2},\ \frac{\partial^2 (\cdot)}{\partial x\partial y},\ \frac{\partial^2 (\cdot)}{\partial y^2},\  \frac{\partial^3 (\cdot)}{\partial x^3},\cdots  \right]^T.
\end{equation}
An approximation to the differential operator $\mathcal{L}(\cdot)$ acting at node $i$ may be expressed as
\begin{equation}
    \mathcal{L}(\cdot)\Big|_i=\sum_{j\in\mathcal{N}_i}(\cdot)_{ji}w_{ji}^d
\end{equation}
where $w_{ji}^d$ are a set of weights to be determined and the superscript $d$ identifies the derivative operator being approximated; $d=x$ is the first derivative with respect to $x$, $d=y$ is the first derivative with respect to $y$, $d=L$ is the Laplacian operator. A given operator $\mathcal{L}(\cdot)$ can be expressed as $\mathcal{L}(\cdot)\Big|_i=\boldsymbol{C}^d\cdot \boldsymbol{\partial}(\cdot)\Big|_i$ where $\boldsymbol{C}^d$ is a vector with its only nonzero elements pointing to derivatives appearing in $\mathcal{L}(\cdot)$. In LABFM, we express the weights as a linear sum of a set of Anisotropic Basis Functions (ABFs) $W(\boldsymbol{r}_{ji}/h_i)$
\begin{equation}
    w_{ji}^d=\boldsymbol{W}_{ji}\cdot \mathbf{\Psi}_i^d
\end{equation}
where $\boldsymbol{W}_{ji}=[W^1_{ji}, W^2_{ji}, W^3_{ji},\cdots]^T$ is the vector of (thus far unspecified) basis functions which must depend on $\boldsymbol{r}_{ji}$ and $\mathbf{\Psi}_i^d$ is a vector of linear coefficients. Defining the moments matrix $\boldsymbol{M}_i$ as the tensor product 
\begin{equation}
\boldsymbol{M}_i=\sum_{j\in\mathcal{N}_i} \boldsymbol{X}_{ji}\otimes \boldsymbol{W}_{ji},
\end{equation}
then for a given choice of ABFs (which consequently gives $\boldsymbol{W}_{ji}$) we can obtain the weights for approximating the operator $\mathcal{L}(\cdot)$ by solving the linear system
\begin{equation}
\label{eq:Explicit_LinearSystem}
    \boldsymbol{M}_i\mathbf{\Psi}_i^d=\boldsymbol{C}^d.
\end{equation}
To obtain an approximation of the operator $\mathcal{L}(\cdot)$ which are of polynomial consistency $m$, we require the size of $\boldsymbol{X}, \boldsymbol{W}, \boldsymbol{C}$ to be $p=(m^2+3m)/2$ in two-dimensions. This places a practical requirement on the size of the computational stencil, as the size of the set $\mathcal{N}_i$ must be greater than $p$ in general. Further issues which need to be accounted for are unacceptable accuracy and poor conditioning of the moments matrix. To overcome this, \cite{king_2022} used orthogonal polynomials to construct series of ABFs and preconditioning which improved residuals by the one order in $h$, both of which we shall also use in this study. We choose to use the Hermite polynomials with a Wendland C2 smoothing kernel to generate weights for the same reasons as outlined in \cite{king_2022}.

\section{Compact LABFM}
\label{sec:Compact}
\subsection{Formulation}
\label{subsec:Formulation}
We now turn to implicit or compact approximations to differential operators using LABFM. Again, let us consider a domain $\Omega\subset \mathbb{R}^2$ which we discretise using $N$ nodal points. We wish to find approximations for spatial differential operators $\mathcal{L}(\cdot)$ acting on a function $\phi$ at each point $i$ within the domain. We may write this as an implicit system
\begin{equation}
\label{eq:Implicit_general_operator}
    \sum_{q\in \mathcal{M}_i}\alpha_{q,i}^d\mathcal{L}(\phi)|_q=\sum_{j\in \mathcal{N}_i}\phi_{ji}w_{ji}^d
\end{equation}
where $\mathcal{M}_i$ is the implicit (left-hand side) stencil, $\mathcal{N}_i$ is the right-hand side stencil, $\alpha_{q,i}^d$ are constant coefficients which can be chosen freely, $\phi_{ji}=\phi_j-\phi_i$, and $w_{ji}^d$ are weights, chosen to ensure the approximation has polynomial consistency of (prescribed) order $m$. Note that we attain explicit approximations when $\alpha_{q,i}^d$ is the Kronecker delta.

In what follows, we focus on extending LABFM to implicit approximations and use the same notation as above, though the ideas presented in this manuscript are applicable to all meshless methods which can be expressed in the same form as \eqref{eq:Implicit_general_operator}. We again define the vector of monomials $\boldsymbol{X}_{ji}$, the vector of partial derivatives $\boldsymbol{\partial} \phi$, a vector of ABFs $\boldsymbol{W}(x_{ji})=\boldsymbol{W}_{ji}$, and a vector of linear coefficients $\boldsymbol{\Psi}$ as in $\S$\ref{sec:Explicit_LABFM}. For given right-hand side and left-hand side stencils, the weights $w_{ji}^d$ are dependent on the choice of $\alpha_{q,i}^d$ and must be determined through the solution of a linear system for each point $i$. Performing a Taylor expansion of the left-hand side of \eqref{eq:Implicit_general_operator} at $x=x_i$ we have
\begin{equation}
\label{Compact_expansion}
    \mathcal{L}(\phi)\Big|_i\sum_{q\in\mathcal{M}_i} \alpha_{q,i}^d+\boldsymbol{\partial}\left(\mathcal{L}\left(\phi\right)\right) \Big|_i \cdot\sum_{q\in\mathcal{M}_i} \boldsymbol{X}_{qi}\alpha_q^i=\sum_{j\in \mathcal{N}_i}\phi_{ji}w_{ji}^d.
\end{equation}
In order to determine the weights $w_{ji}^d$ we choose our series of ABFs and solve the linear system
\begin{equation}
    \label{eq:Implicit_LinearSystem}
    \boldsymbol{M}_i\mathbf{\Psi}_i^d=\boldsymbol{\Tilde{C}}^d
\end{equation}
for $\mathbf{\Psi}_i$, where the moments matrix $\boldsymbol{M}_i$ is the same as in \eqref{eq:Explicit_LinearSystem}, but here $\boldsymbol{C}^d$ has been replaced by the vector $\boldsymbol{\Tilde{C}}^d$ which depends on the choice of implicit coefficients through
\begin{equation}
    \boldsymbol{\Tilde{C}}^d=\sum_{q\in\mathcal{M}_i}\alpha_{q,i}^d \mathcal{L}(\boldsymbol{X}_{qi}).
\end{equation}
It is clear that we can recover the explicit formulation by setting $\alpha_{q,i}^d=\delta_{iq}$ where $\delta_{iq}$ is the Kronecker delta.

For a given choice of ABFs, $\mathcal{M}_i$, $\mathcal{N}_i$ and polynomial consistency $m$, we are free to choose the $\alpha_{q,i}^d$, from which $\boldsymbol{\Psi}_i$ can then be uniquely determined and used to generate the weights $w_{ji}^d$. Once weights at each point have been obtained, approximations to the operator $\mathcal{L}$ acting on a specific function $\phi$ are determined at all points in domain simultaneously by solving a global linear system. We now move to a discussion of how to choose the implicit coefficients.

\subsection{Choice of Implicit Stencil \& Constraints on Coefficients}
\label{subsec:Stencil_choice}
With the above formulation of an implicit version of LABFM, it is necessary to choose the optimal implicit stencil $\mathcal{M}_i$ and coefficients $\alpha_q^i$. To do so, we must define a measurement of what is meant by `optimal'. For any given choice of coefficients we impose that the weights $w_{ji}^d$ are chosen to yield the same desired polynomial consistency, therefore choosing a higher order approximation (as may done in compact finite-differences) is not possible. Instead, we choose to borrow another measurement of the accuracy of derivative approximations from finite-differences, that of resolving power \cite{Lele,MKLABFM_Arxiv}. However, before maximising the resolving power of our approximation we must first choose the implicit stencil $\mathcal{M}_i$ at each point $i$. 

The choice of $\mathcal{M}_i$ requires a choice of the number of nodes in this stencil and their relative location to the point $i$. As implicit approximations require a global linear system to be solved, maintaining as sparse a system as possible is desirable to minimise computational cost. We therefore choose to limit the number of neighbours in the left-hand side stencil to a maximum of 9 for gradient operators, and a maximum of 17 for Laplacians (including the central node $i$ in both cases). Dependency on this choice will be examined in convergence studies in $\S$\ref{sec:Resolving Power} and \ref{sec:Convergence_Studies}.

\begin{figure}
\begin{center}
\setlength{\unitlength}{1cm}
  \begin{tikzpicture}(16,9.25)(0,0.25)
    \filldraw[red] (0,0) circle (0.1);
    \node[] at (0.2,-0.2) {$i$};

    \draw(0,0) circle (2.5);
    \draw[<->, thick] (-2.5,0)--(-0.1,0); 
    \node[] at (-1.25, -0.3) {$2h_i$};

    \draw(1.8,1.6) circle (0.1);
    \draw(0.3,1.7) circle (0.1);
    \draw(0.8,1) circle (0.1);
    \draw(1.5,0.7) circle (0.1);
    \filldraw[red] (2.2,0.2) circle (0.1);
    \filldraw[red] (0.6,-0.1) circle (0.1);
    \draw(1.2,-0.4) circle (0.1);
    \draw(1.4,-1.2) circle (0.1);
    \draw(0.2,-2.2) circle (0.1);
    \draw(-0.3,-1.4) circle (0.1);
    \draw(-1.2,-1.7) circle (0.1);
    \draw(-1.8,-0.8) circle (0.1);
    \filldraw[red] (-0.5,-0.35) circle (0.1);
    \filldraw[red] (-1.7,0.25) circle (0.1);
    \draw(-1.4,0.8) circle (0.1);
    \draw(-0.5,0.6) circle (0.1);
    \draw(-0.8,2) circle (0.1);
    \draw(-0.05,1.3) circle (0.1);

    \draw(1.25,1.95) circle (0.1);
    \draw(-1.85,1.4) circle (0.1);
    \draw(-1,1.25) circle (0.1);
    \draw(0.35,0.6) circle (0.1);
    \draw(-0.8,-0.9) circle (0.1);
    \draw(-0.6,-2.05) circle (0.1);
    \draw(0.55,-0.95) circle (0.1);
    \draw(0.9,-1.65) circle (0.1);
    \draw(2,-0.65) circle (0.1);

    \filldraw[blue] (6,0) circle (0.1);
    \node[] at (6.2,-0.2) {$i$};

    \draw(6,0) circle (2.5);
    \draw[<->, thick] (3.5,0)--(5.9,0); 
    \node[] at (4.75, -0.3) {$2h_i$};

    \draw(7.8,1.6) circle (0.1);
    \filldraw[blue] (6.3,1.7) circle (0.1);
    \draw(6.8,1) circle (0.1);
    \draw(7.5,0.7) circle (0.1);
    \draw(8.2,0.2) circle (0.1);
    \draw(6.6,-0.1) circle (0.1);
    \draw(7.2,-0.4) circle (0.1);
    \draw(7.4,-1.2) circle (0.1);
    \filldraw[blue] (6.2,-2.2) circle (0.1);
    \filldraw[blue] (5.7,-1.4) circle (0.1);
    \draw(4.8,-1.7) circle (0.1);
    \draw(4.2,-0.8) circle (0.1);
    \draw(5.5,-0.35) circle (0.1);
    \draw(4.3,0.25) circle (0.1);
    \draw(4.6,0.8) circle (0.1);
    \draw(5.5,0.6) circle (0.1);
    \draw(5.2,2) circle (0.1);
    \filldraw[blue] (5.95,1.3) circle (0.1);

    \draw(7.25,1.95) circle (0.1);
    \draw(4.15,1.4) circle (0.1);
    \draw(5,1.25) circle (0.1);
    \draw(6.35,0.6) circle (0.1);
    \draw(5.2,-0.9) circle (0.1);
    \draw(5.4,-2.05) circle (0.1);
    \draw(6.55,-0.95) circle (0.1);
    \draw(6.9,-1.65) circle (0.1);
    \draw(8,-0.65) circle (0.1);

 \filldraw[black] (12,0) circle (0.1);
    \node[] at (12.2,-0.2) {$i$};

    \draw(12,0) circle (2.5);
    \draw[<->, thick] (9.5,0)--(11.9,0); 
    \node[] at (10.75, -0.3) {$2h_i$};

 \draw(13.8,1.6) circle (0.1);
    \filldraw[black] (12.3,1.7) circle (0.1);
    \draw(12.8,1) circle (0.1);
    \draw(13.5,0.7) circle (0.1);
    \filldraw[black] (14.2,0.2) circle (0.1);
    \filldraw[black] (12.6,-0.1) circle (0.1);
    \draw(13.2,-0.4) circle (0.1);
    \draw(13.4,-1.2) circle (0.1);
    \filldraw[black] (12.2,-2.2) circle (0.1);
    \filldraw[black] (11.7,-1.4) circle (0.1);
    \draw(10.8,-1.7) circle (0.1);
    \draw(10.2,-0.8) circle (0.1);
    \filldraw[black] (11.5,-0.35) circle (0.1);
    \filldraw[black] (10.3,0.25) circle (0.1);
    \draw(10.6,0.8) circle (0.1);
    \draw(11.5,0.6) circle (0.1);
    \draw(11.2,2) circle (0.1);
    \filldraw[black] (11.95,1.3) circle (0.1);

    \draw(13.25,1.95) circle (0.1);
    \draw(10.15,1.4) circle (0.1);
    \draw(11,1.25) circle (0.1);
    \draw(12.35,0.6) circle (0.1);
    \draw(11.2,-0.9) circle (0.1);
    \draw(11.4,-2.05) circle (0.1);
    \draw(12.55,-0.95) circle (0.1);
    \draw(12.9,-1.65) circle (0.1);
    \draw(14,-0.65) circle (0.1);

    \node[] at (-2.5,-2.5) {(a)};
    \node[] at (3.5,-2.5) {(b)};
    \node[] at (9.5,-2.5) {(c)};
   \end{tikzpicture}
\end{center}
\caption{An illustrative example of how an implicit stencil is chosen from within the right-hand stencil for compact LABFM, with $Q_i=5$. (a) Gradient Operator with respect to $x$; (b) Gradient Operator with respect to $y$; (c) Laplacian Operator. Solid colour nodes correspond to those in the implicit stencil.}
\label{fig:Implicit_stencils}
\end{figure}
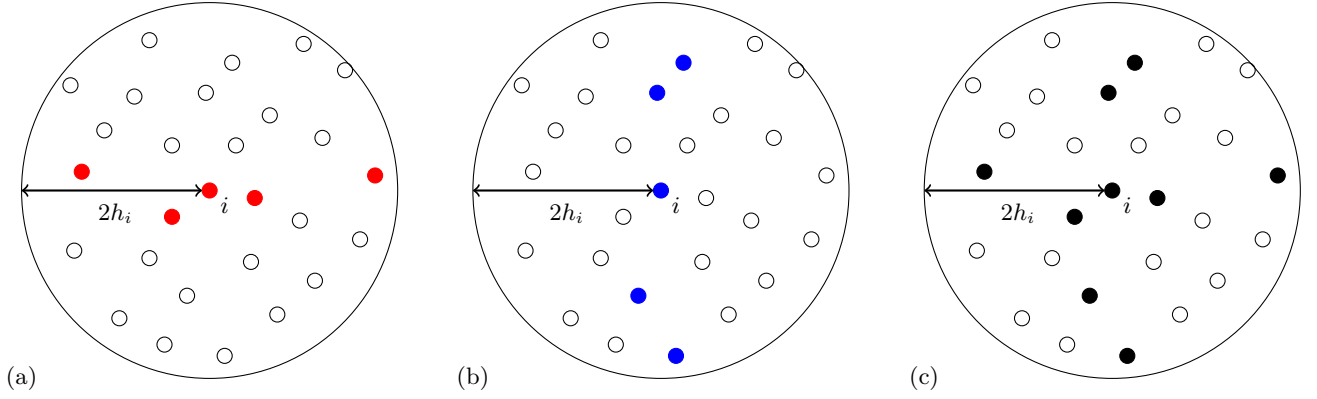

We now turn to a discussion of how we choose $\mathcal{M}_i$ for various operators. Let us first consider operators which only require derivatives with respect to one variable. We denote the number of nodes in the left-hand side stencil as $Q_i$, and we restrict potential choices in all cases to nodes within the right-hand side stencil $\mathcal{N}_i$. We introduce two sets, $A_i^x, A_i^y$ defined through
\begin{equation}
    \label{eq:sets}
    A_i^x=\{x_j-x_i | j\in \mathcal{N}_i  \}, \quad A_i^y=\{y_j-y_i | j\in \mathcal{N}_i  \}.
\end{equation}
To approximate derivatives with respect to $x$ we use the nodes corresponding to the $Q_i$ smallest values in $A_i^y$, and correspondingly for derivatives with respect to $y$ we take nodes with the $Q_i$ smallest values in $A_i^x$. For operators which require derivatives in two variables such as Laplacians, we take a union of the choices for single variable coefficients, leading to $2Q_i-1$ points in the stencil. Graphical examples of the choices of nodes for the implicit stencil can be found in Fig.\ref{fig:Implicit_stencils}.

Before optimising our choice of implicit coefficients, we first note that there are practical constraints on the  $\alpha_{q,i}^d$ values to guarantee acceptable accuracy of the approximation. An important constraint is to require that the implicit matrix (i.e. the global matrix containing each $\boldsymbol{\alpha}_i$ on row $i$) is well-conditioned, as well as the matrix being diagonally dominant, with negative coefficients also typically undesirable. The simplest way to ensure these constraints are met is for the implicit coefficients to be of a prescribed form which prevents such issues. To this end, we set 
\begin{equation}
\label{eq:implicit_coeffs}
    \alpha_{q,i}^d=e^{-(a_x^2 x_{qi}^2+a_y^2 y_{qi}^2)/s_i}
\end{equation}
and optimise our choice of $a_x$ and $a_y$. To prevent the distribution from becoming too flat (which results in issues due to poor condition of the implicit matrix) we limit the minimum value of $a_x, a_y$ as outlined below.

\subsection{Resolving Power Optimisation}
\label{subsec:Optimise_coefficients}
Taking into account the restraints described above, we now choose our implicit coefficients in order to maximise resolving power. We showed in \cite{MKLABFM_Arxiv} that in meshfree methods the notion of resolving power requires a multi-dimensional Fourier analysis due to the scattered nodal distributions. We begin by considering a gradient operator with respect to $x$. We may perform a similar analysis to that of \cite{MKLABFM_Arxiv}, though we must now retain the extra implicit terms on the left-hand side of \eqref{eq:Implicit_general_operator}. We consider a generic function $\phi$ in a doubly periodic two-dimensional domain $(x,y)\in[0,2\pi]\times[0,2\pi]$. This function can be expressed in terms of a two-dimensional Fourier series 
\begin{equation}
\label{eq::fourier_series}
    \phi(x,y)=\sum_{k_x\in \mathbb{Z},}\sum_{k_y\in \mathbb{Z}} c_{k_x,k_y}e^{ i(k_xx+k_yy)},
\end{equation}
where $k_x, k_y$ are real wavenumbers $\Big(\sqrt{k_x^2+k_y^2}\leq k_{Ny}$, where $k_{Ny}\coloneqq\pi/s_i$ is the Nyquist wavenumber$\Big)$. Substituting this expression into \eqref{eq:Implicit_general_operator} for a gradient operator with respect to $x$ allows us to compare the analytic wavenumber $k_x$ to the numerically approximated effective wavenumber $k_{eff}$ for each wavenumber pair $(k_x,k_y)$. Rearranging \eqref{eq:Implicit_general_operator} with $\phi$ given by \eqref{eq::fourier_series} we find that the effective wavenumber for a gradient operator with respect to $x$ is given by
\begin{gather}
\label{eq:k_eff}
    k_{eff}(k_x,k_y)=\frac{\gamma_1 \lambda_2+\gamma_2\lambda_1 +i(\gamma_2\lambda_2-\gamma_1\lambda_1)}{\lambda_1^2+\lambda_2^2},
\end{gather}
where $\gamma_{1},\gamma_2,\lambda_1,\lambda_2$ are defined through
\begin{gather}
    \gamma_1=\sum_{j\in\mathcal{N}_i}\sin(k_xx_{ji}+k_y y_{ji})w_{ji}^x,\quad \gamma_2=\sum_{j\in\mathcal{N}_i}[1-\cos(k_xx_{ji}+k_y y_{ji})]w_{ji}^x,\\
    \lambda_1=\sum_{q\in\mathcal{M}_i}\sin(k_xx_{qi}+k_yy_{qi})\alpha_{q,i}^x,\quad  \lambda_2=\sum_{q\in\mathcal{M}_i}\cos(k_xx_{qi}+k_yy_{qi})\alpha_{q,i}^x.
\end{gather}
The error in effective wavenumber obviously has a real and imaginary part. As will be shown in \S\ref{sec:Resolving Power} the imaginary part of $k_{eff}$ is typically smaller than $|k_x-k_{eff}|$, therefore we choose $a_x$, $a_y$ in order to minimise the deviation of the $\Re\{k_{eff}\}$ from $k_x$. Initial numerical experiments revealed that smaller values of $a_x, a_y$ tended to minimise $\Re\{k_x-k_{eff}\}$ up to a point whereupon the global linear system became poorly conditioned. We therefore maximise resolving power at each point $i$ by systematically reducing $a_x$ through the following procedure:
\begin{enumerate}
 \item Set $a_y$ such that $\sum_{q\in \mathcal{M}_i,q\neq i}\alpha_{q,i}^x\leq 2$ (note $a_y$ varies on a stencil-by-stencil basis).
    \item Set $a_x$, to a sufficiently large enough value to ensure a well-conditioned local linear system.
    \item Reduce $a_x=a_x-0.01$ and calculate $w_{ji}^x$.
    \item Calculate resolving power using \eqref{eq:k_eff} at a range of wavenumbers in $(k_x,k_y)$ space.
    \item If $\frac{\Re\{k_{eff}\}}{k_x}\leq 1.005$ for all $k_x, k_y$ then repeat steps 3 and 4.
    \item If $\frac{\Re\{k_{eff}\}}{k_x}> 1.005$, terminate algorithm and use previous value of $a_x$ to calculate $\alpha_{q,i}^x$ and $w_{ji}^x$.
\end{enumerate}
This algorithm ensures that diagonal dominance of our implicit system is maintained and resolving power is maximised whilst ensuring an excitation of the real part of the wavenumber of less than $0.5\%$ for all $(k_x,k_y)$. For first derivatives with respect to $y$, the procedure is the same but with subscripts $x$ and $y$ swapped.

A similar procedure can be performed for the Laplacian operator. In that case we find that the effective wavenumber $q_{eff}^2$ at each point $i$ can be expressed in the form
\begin{equation}
\label{eq:qeff}
    q_{eff}^2=\frac{\hat{\gamma}_1\hat{\lambda}_1 +\hat{\gamma}_2\hat{\lambda}_2-i(\hat{\gamma}_2\hat{\lambda}_1+\hat{\gamma}_1\hat{\lambda}_2)}{\hat{\lambda}_1^2+\hat{\lambda}_2^2},
\end{equation}
where we have defined
\begin{gather}
    \hat{\gamma}_1=\sum_{j\in\mathcal{N}_i}[1-\cos(k_xx_{ji}+k_y y_{ji})]w_{ji}^L,\quad \hat{\gamma}_2=\sum_{j\in\mathcal{N}_i}\sin(k_xx_{ji}+k_y y_{ji})w_{ji}^L,\\
    \hat{\lambda}_1=\sum_{q\in\mathcal{M}_i}\cos(k_xx_{qi}+k_yy_{qi})\alpha_{q,i}^L,\quad  \hat{\lambda}_2=\sum_{q\in\mathcal{M}_i}\sin(k_xx_{qi}+k_yy_{qi})\alpha_{q,i}^L.
\end{gather}
For spectral resolving power we would obtain $q^2_{eff}=q^2\coloneqq k_x^2+k_y^2$. Note that the choice of implicit stencil $\mathcal{M}_i$ and coefficients $\alpha_{q,i}^d$ here differ from the choices made for the gradient operator above. Again considering only the real part in the error of the effective wavenumber we maximise the resolving power of the Laplacian operator using a similar procedure to that detailed for the gradient operator. For the Laplacian operator however, we omit Step 1 and instead require
\begin{equation}
    a_x=a_y
\end{equation}
as well as appropriately changing steps 4,5, and 6 to compare $\Re\{q_{eff}^2\}$ to $q^2$.

The resolving power of a scheme is resolution independent. However, it is clear from \eqref{eq:k_eff} \& \eqref{eq:qeff} that the resolving power of a meshless method is dependent on the specific discretisation as it is different at each point. To obtain measures of the resolving power of an approximation to a differential operator which are valid across the whole nodal distribution we take a root mean square average of both the real and imaginary part of the effective wavenumber at all points.

For all the compact formulations outlined above, a separate global linear system must be solved for each operator in order to obtain the implicit approximation. For time-dependent problems evolved using explicit time integration, global linear systems must be solved at each time-step. In this study, we use a GMRES algorithm with AMG preconditioner via the HYPRE library \cite{hypre} to perform a sparse iterative solve on these systems. 

\section{Resolving Power}
\label{sec:Resolving Power}
As is clear from the formulation above, using compact LABFM instead of explicit LABFM does not improve the formal convergence properties of approximations. Therefore, to determine the improvement in approximations of derivative operators provided by compact LABFM across all scales it is necessary to consider the resolving power of each operator. 

To assess resolving power, we discretise a periodic domain $\Omega$ using a disordered nodal distribution with a constant average nodal spacing $s=s_i$ generated in the same manner as in \cite{king_2022}, using a propagating front algorithm before performing ten iterations of a shifting algorithm as in \cite{Fornberg&Flyer_disc} to avoid clustering of particles. Throughout this section and the rest of this manuscript we present results in terms of orders of convergence for our operators, not polynomial consistency, noting that the same order of convergence requires a Laplacian operator to be polynomially consistent to one degree higher than gradient operators. This results in consistent convergence properties across spatial operators in our studies of PDE systems, though we require different size stencils for different operators. Table \ref{tab:stencil_sizes} shows stencil sizes used for different operators at each order of approximation.
\begin{table}[]
    \centering
    \begin{tabular}{ |*{7}{c|} }
    \hline
\quad Order of Convergence \quad  & \multicolumn{2}{c|}{2}
            & \multicolumn{2}{c|}{4}               \\
    \hline
Operator   & \quad   Gradient \quad   & \quad Laplacian  \quad & \quad  Gradient \quad & \quad  Laplacian \quad \\
    \hline
 Stencil size $h/s_i$   &   1.2  &   1.35  &   1.4  &   1.7    \\
    \hline
    \end{tabular}
    \caption{Stencil sizes for gradient and Laplacian operators for different orders of convergence. Laplacian operators require larger stencils as a higher polynomial consistency is required for convergence to a given order.}
    \label{tab:stencil_sizes}
\end{table}
In this manuscript we will consider compact approximations with various sizes of implicit stencil. The choices of $Q_i$ we make are given in Table \ref{tab:stencils}. The labelling of the schemes remains constant throughout the manuscript, and when approximations to both types of operators are required (for example in \S\ref{subsec:Burgers_Equation}) we are consistent in choice of scheme (e.g. the scheme (h) will choose (h) in Table \ref{tab:stencils} for both gradients and Laplacians).
\begin{table}[]
    \centering
    \begin{tabular}{ |*{3}{c|} }
     \multicolumn{3}{c}{Gradients}
\vspace{0.5em}
\\

   \quad Scheme \quad & \quad   Order \quad   & \quad Implicit Stencil \quad\\ 
    \hline
 (a)  &$2^{nd}$& Explicit\\
 (b)  &$2^{nd}$& $Q_i=3$\\
 (c)  &$2^{nd}$& $Q_i=5$\\
 (d) &$2^{nd}$& $Q_i=7$ \\
 (e) & $4^{th}$ & Explicit\\
 (f) & $4^{th}$ & $Q_i=5$\\
 (g) &$4^{th}$ & $Q_i=7$\\
 (h) & $4^{th}$ & $Q_i=9$\\ 
    \end{tabular}\hspace{4em}
   \begin{tabular}{ |*{3}{c|} }
          \multicolumn{3}{c}{Laplacians}
  \vspace{0.5em}
            \\
    
   \quad Scheme \quad & \quad   Order\quad   & \quad Implicit Stencil \quad  \\
    \hline
 (a)  &$2^{nd}$& Explicit\\
 (b)  &$2^{nd}$& $2Q_i-1=5$ \\
 (c)  &$2^{nd}$& $2Q_i-1=9$\\
 (d)  &$2^{nd}$& $2Q_i-1=13$\\
 (e) & $4^{th}$ & Explicit \\
 (f) & $4^{th}$ & $2Q_i-1=9$ \\
 (g) & $4^{th}$ & $2Q_i-1=13$ \\
 (h) & $4^{th}$ & $2Q_i-1=17$ \\
    \end{tabular}
    \caption{Table showing discretisation schemes used throughout the manuscript (labelling constant throughout). Left table: gradient operators. Right table: Laplacians. Note that order here refers to order of convergence of the approximation.}
    \label{tab:stencils}
\end{table}

We begin by considering a gradient operator with respect to $x$. For an accurate assessment of the resolving power of the approximation it is necessary to consider all of wavenumber space $(k_x,k_y)$, as opposed to the one-dimensional analysis common to finite-difference studies. Consequently, improved resolving power across all values of $(k_x,k_y)$ typically requires a larger implicit stencil than compact finite-difference schemes, and here we present results up to $Q_i=9$. 
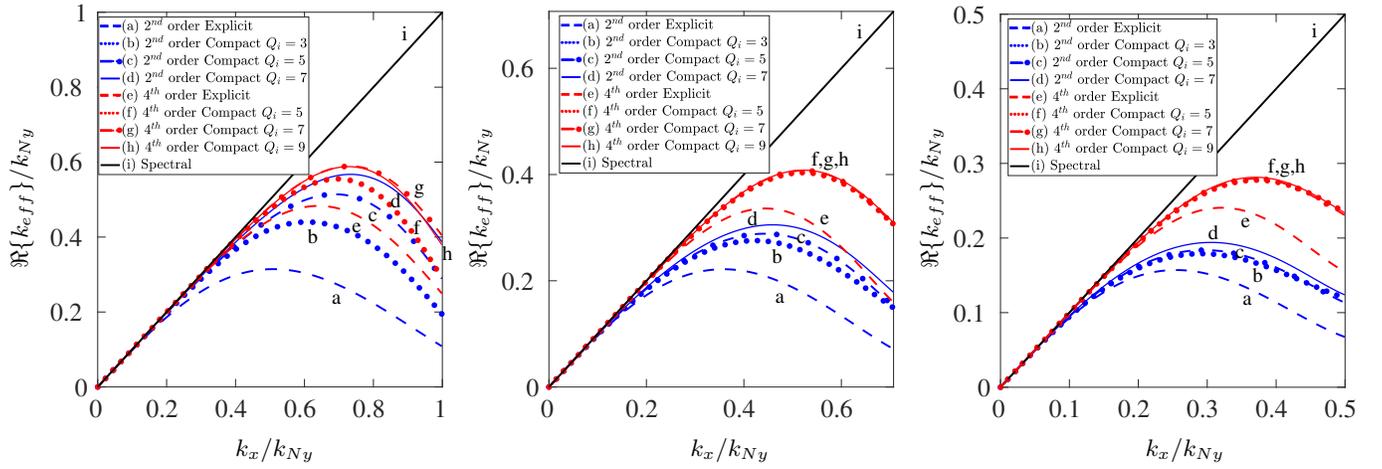
\begin{figure}[t]
    \begin{center}
    \setlength{\unitlength}{1cm}
    \begin{picture}(18,5)(0,0)
    \put(0.3,0){\includegraphics[width=0.33\linewidth]{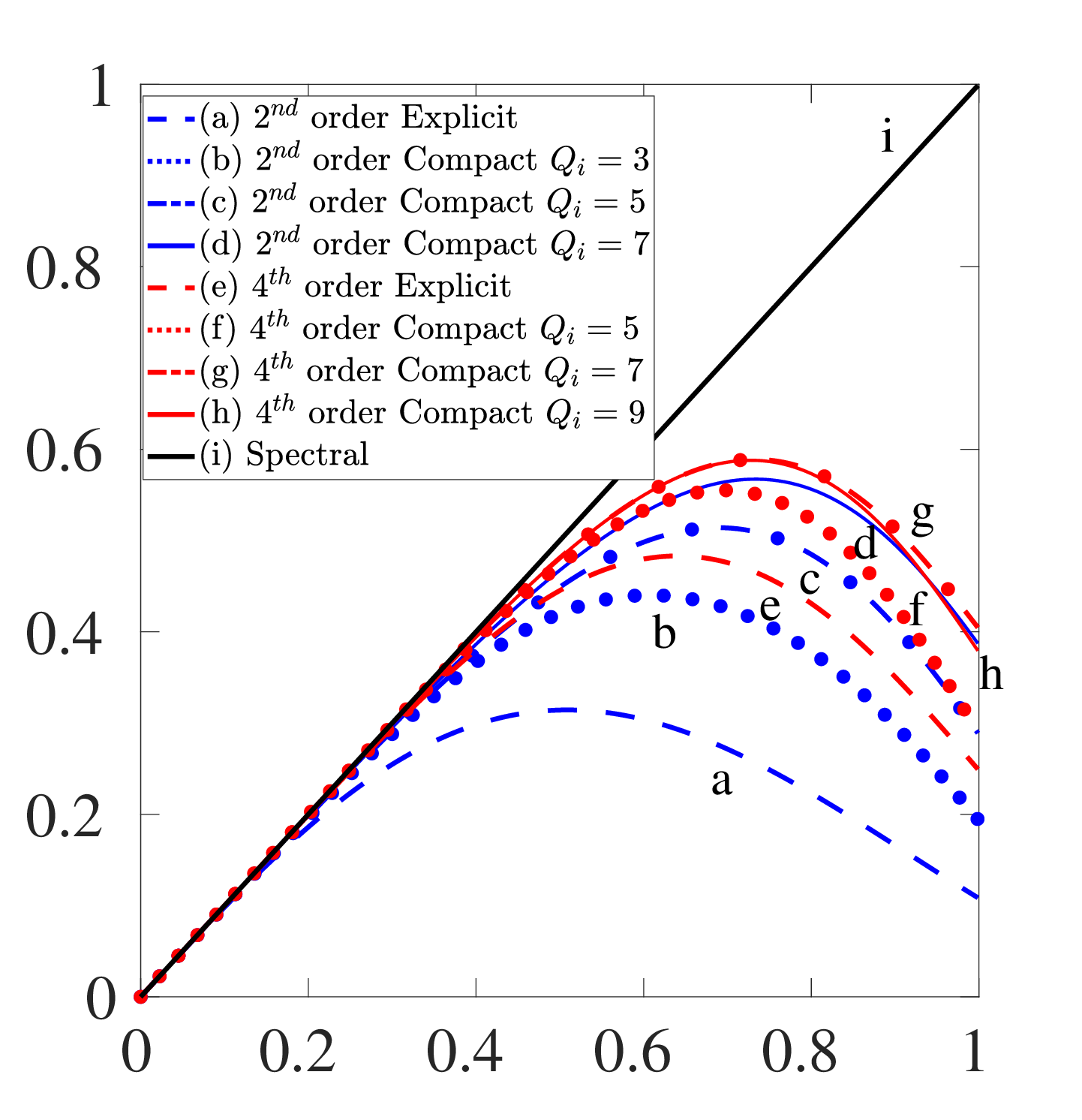}}
    \put(2.9,-0.25){$k_x/k_{Ny}$}
    \put(6.3,00){\includegraphics[width=0.33\linewidth]{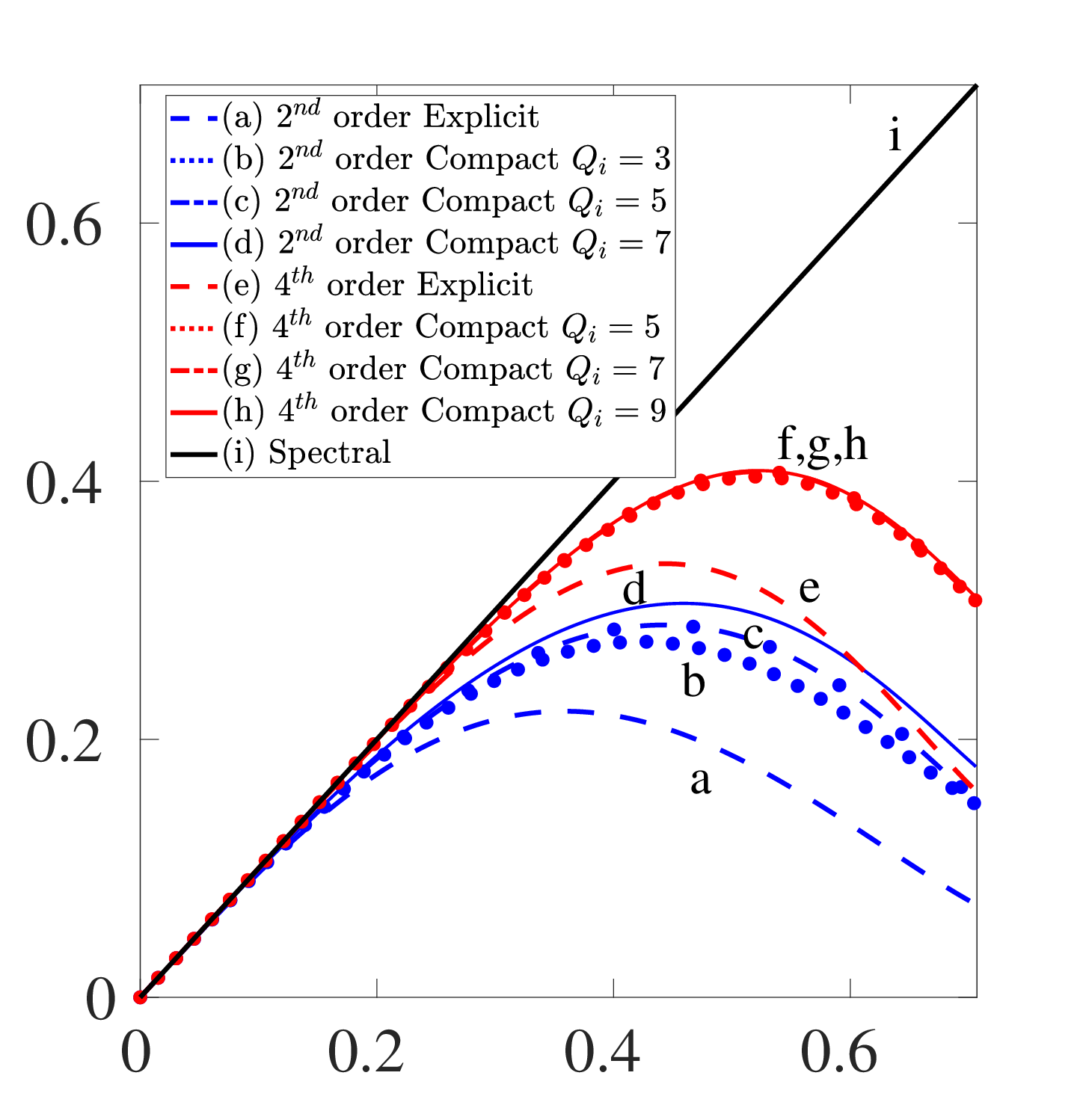}}
    \put(8.95,-0.25){$k_x/k_{Ny}$}
    \put(12.3,00){\includegraphics[width=0.33\linewidth]{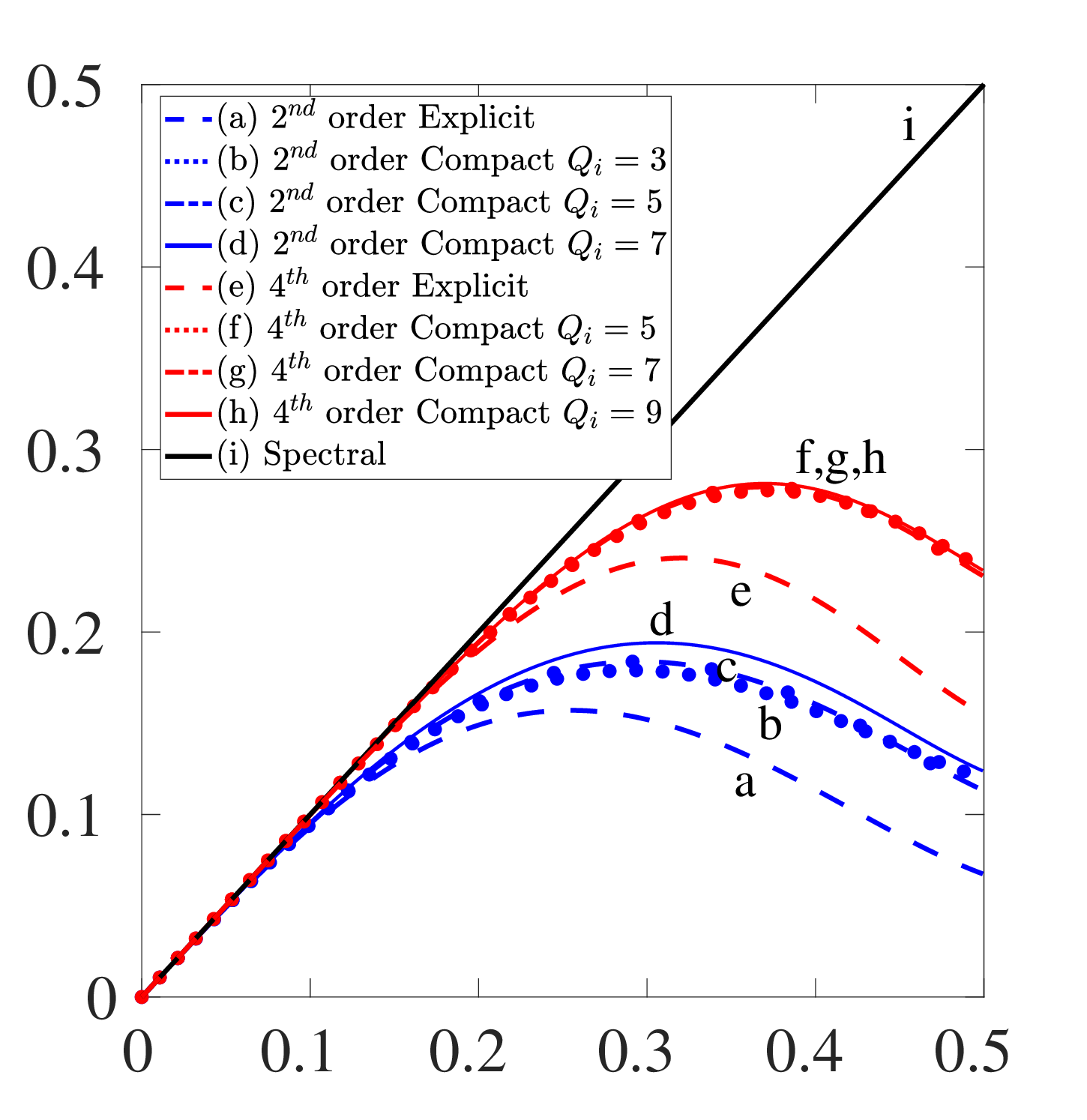}}
    \put(15.0,-0.25){{$k_x/k_{Ny}$}}
  \put(1.1, 4){\color{red}{\thickhyphen}} 
  \put(1.31, 4.04){\color{red}{{\scalebox{.6}{$\bullet$}}}}

  \put(1.09, 4.295){\color{red}{{\scalebox{.3}{$\bullet$}}}}
  \put(1.15, 4.295){\color{red}{{\scalebox{.3}{$\bullet$}}}}
  \put(1.2, 4.295){\color{red}{{\scalebox{.3}{$\bullet$}}}}
  \put(1.25, 4.295){\color{red}{{\scalebox{.3}{$\bullet$}}}}
  \put(1.31, 4.295){\color{red}{{\scalebox{.3}{$\bullet$}}}}

\put(1.1, 4.47){\color{red}{\thickhyphent}} 
\put(1.27, 4.47){\color{red}{\thickhyphent}} 

\put(1.1, 4.92){\color{blue}{\thickhyphen}} 
  \put(1.31, 4.96){\color{blue}{{\scalebox{.6}{$\bullet$}}}}

  \put(1.09, 5.215){\color{blue}{{\scalebox{.3}{$\bullet$}}}}
  \put(1.15, 5.215){\color{blue}{{\scalebox{.3}{$\bullet$}}}}
  \put(1.2, 5.215){\color{blue}{{\scalebox{.3}{$\bullet$}}}}
  \put(1.25, 5.215){\color{blue}{{\scalebox{.3}{$\bullet$}}}}
  \put(1.31, 5.215){\color{blue}{{\scalebox{.3}{$\bullet$}}}}
\put(1.1, 5.39){\color{blue}{\thickhyphent}} 
\put(1.27, 5.39){\color{blue}{\thickhyphent}} 

\put(7.23, 4.02){\color{red}{\thickhyphen}} 
  \put(7.42, 4.06){\color{red}{{\scalebox{.6}{$\bullet$}}}}

  \put(7.22, 4.315){\color{red}{{\scalebox{.3}{$\bullet$}}}}
  \put(7.27, 4.315){\color{red}{{\scalebox{.3}{$\bullet$}}}}
  \put(7.32, 4.315){\color{red}{{\scalebox{.3}{$\bullet$}}}}
  \put(7.38, 4.315){\color{red}{{\scalebox{.3}{$\bullet$}}}}
  \put(7.43, 4.315){\color{red}{{\scalebox{.3}{$\bullet$}}}}

\put(7.23, 4.48){\color{red}{\thickhyphent}} 
\put(7.4, 4.48){\color{red}{\thickhyphent}} 

\put(7.23, 4.94){\color{blue}{\thickhyphen}} 
  \put(7.42, 4.975){\color{blue}{{\scalebox{.6}{$\bullet$}}}}

  \put(7.22, 5.23){\color{blue}{{\scalebox{.3}{$\bullet$}}}}
  \put(7.27, 5.23){\color{blue}{{\scalebox{.3}{$\bullet$}}}}
  \put(7.32, 5.23){\color{blue}{{\scalebox{.3}{$\bullet$}}}}
  \put(7.37, 5.23){\color{blue}{{\scalebox{.3}{$\bullet$}}}}
  \put(7.43, 5.23){\color{blue}{{\scalebox{.3}{$\bullet$}}}}
\put(7.23, 5.4){\color{blue}{\thickhyphent}} 
\put(7.4, 5.4){\color{blue}{\thickhyphent}} 

\put(13.2, 3.985){\color{red}{\thickhyphen}} 
  \put(13.38, 4.015){\color{red}{{\scalebox{.6}{$\bullet$}}}}

  \put(13.19, 4.275){\color{red}{{\scalebox{.3}{$\bullet$}}}}
  \put(13.24, 4.275){\color{red}{{\scalebox{.3}{$\bullet$}}}}
  \put(13.29, 4.275){\color{red}{{\scalebox{.3}{$\bullet$}}}}
  \put(13.345, 4.275){\color{red}{{\scalebox{.3}{$\bullet$}}}}
  \put(13.4, 4.275){\color{red}{{\scalebox{.3}{$\bullet$}}}}

\put(13.19, 4.445){\color{red}{\thickhyphent}} 
\put(13.36, 4.445){\color{red}{\thickhyphent}} 

\put(13.2, 4.9){\color{blue}{\thickhyphen}} 
  \put(13.38, 4.945){\color{blue}{{\scalebox{.6}{$\bullet$}}}}

  \put(13.19, 5.19){\color{blue}{{\scalebox{.3}{$\bullet$}}}}
  \put(13.24, 5.19){\color{blue}{{\scalebox{.3}{$\bullet$}}}}
  \put(13.29, 5.19){\color{blue}{{\scalebox{.3}{$\bullet$}}}}
  \put(13.34, 5.19){\color{blue}{{\scalebox{.3}{$\bullet$}}}}
  \put(13.4, 5.19){\color{blue}{{\scalebox{.3}{$\bullet$}}}}
\put(13.19, 5.365){\color{blue}{\thickhyphent}} 
\put(13.36, 5.365){\color{blue}{\thickhyphent}} 

    \put(-0.1,2.25){\rotatebox{90}{$\Re\{k_{eff}\}/k_{Ny}$}}
    \put(6.0,2.25){\rotatebox{90}{$\Re\{k_{eff}\}/k_{Ny}$}}
     \put(12.0,2.25){\rotatebox{90}{$\Re\{k_{eff}\}/k_{Ny}$}}

    \end{picture}
    \end{center}
    \caption{Real part of resolving power of the gradient operator for second order accurate (blue lines) and fourth order (red lines) accurate operators, with schemes as labelled in Table \ref{tab:stencils} except for the spectral scheme (i) (exact differentiation). Left Panel $k_y=0$, Centre Panel $k_y=k_x$, Right Panel $k_y=2k_x$.}
    \label{fig:rp_grad_real}
\end{figure}
Fig.\ref{fig:rp_grad_real} compares the real part of the resolving power for gradient operators along various lines in $(k_x,k_y)$ space (note in all cases spectral accuracy is achieved if $k_{eff}=k_x$). For order $m=2$ operators we observe increasing improvement in resolution of moderate to small scales as the implicit stencil becomes larger. Resolving power appears most improved along the line $k_y=0$, along which two compact $2^{nd}$ order schemes outperform the explicit $4^{th}$ order scheme, with scheme (d) even outperforming one of the compact $4^{th}$ order schemes. For $4^{th}$ order schemes we again observe a consistent improvement in resolving power for compact approximations. Schemes (g) \& (h) are indistinguishable graphically for the three lines shown, with (f) also having very similar resolving power along the lines $k_y=k_x,2 k_x$. Given the reduced computational time associated with solving the global sparse linear system for smaller compact stencils, these graphs suggest that scheme (g) would likely be the optimum choice for a compact $4^{th}$ order approximation if we consider only the real part of the effective wavenumber. 

To further illustrate the resolving efficiency of a gradient approximation we introduce the measure
\begin{equation}
    \varepsilon_1(k_x,k_y)=\frac{|\Re\{k_{eff}(k_x,k_y)\}|}{k_x}
\end{equation}
which can be used to assess the range of ($k_x,k_y$) over which waves are well resolved. We wish for $\varepsilon_1$ to be close to zero over as wide a range of $(k_x,k_y)$ as possible. We can use the measure of resolving efficiency to quantify the improvements in resolving power shown in Fig.\ref{fig:rp_grad_real}. Table \ref{tab:grad_rp} gives values of $k_x$ at which the error of the real part of the effective wavenumber becomes greater than $0.1\%,1\%,10\%$ (corresponding to resolving efficiencies of $\varepsilon_1=0.001,0.01,0.1$ respectively) of $k_x$ for various lines in wavenumber space for each scheme. Note that higher order convergence results in $\varepsilon_1\rightarrow0$ faster in the limit $k_x/k_{Ny}\rightarrow 0$, which is why the $4^{th}$ order explicit scheme remains below $0.1\%$ error for longer than $2^{nd}$ order compact schemes. Within a family of schemes of a given order, as the size of the implicit stencil increases the resolving power metric improves across all of wavenumber space, though the greatest improvement is typically seen along $k_y=0$. In particular, along this line in wavenumber space the $2^{nd}$ order scheme (d) has an error of under $10\%$ for over half of wavenumbers, compared to under a quarter for the comparable explicit scheme. Schemes (c) \& (d) outperform the $4^{th}$ order explicit schemes here, with $\varepsilon_1=0.1$ being crossed at higher values of $k_x$. Similar improvements are seen for $4^{th}$ order schemes, although again marginal difference between (g) \& (h) is observed.
\begin{table}[]
    \centering
    \begin{tabular}{ *{10}{c|} }
  & \multicolumn{3}{c|}{$\varepsilon_1=0.001$}
            & \multicolumn{3}{c|}{$\varepsilon_1=0.01$}
             & \multicolumn{3}{c|}{$\varepsilon_1=0.1$}
            \\

   \quad Numerical Scheme \quad & \   $k_y=0$ \quad   & \ $k_y=k_x$  \quad & \ $k_y=2k_x$ \quad & \   $k_y=0$ \quad   & \ $k_y=k_x$  \quad & \  $k_y=2k_x$ \quad &\  $k_y=0$ \quad   & \ $k_y=k_x$  \quad & \  $k_y=2k_x$ \quad \\
    \hline
 (a) &0.03  &0.02  &0.01   &0.07  &0.05  &0.04   &0.24 &0.17 &0.12\\
 (b) &0.05&0.02  &0.01  & 0.16  &0.07  & 0.04 & 0.43&0.22 & 0.14 \\
 (c) &0.06 &0.02  &0.01  &0.19  &0.07  &0.04  &0.50 &0.23  & 0.14\\
 (d) & 0.08 &0.02  &0.01  &0.24  &0.08  &0.04  &0.57 &0.25 & 0.15 \\
 (e) &0.14&0.1  &0.07  & 0.26   &0.18&0.13  &0.49 &0.34  &0.24 \\
 (f) &0.19&0.13  &0.08  &0.33  &0.23  &0.15  & 0.58&0.41  &0.28 \\
 (g)& 0.21&0.13  &0.08  &0.36  &0.23  &0.15  & 0.62&0.42  &0.28 \\
 (h) &0.21&0.13  &0.09  & 0.36 &0.24  & 0.16 &0.62 &0.43  &0.29 \\
    \end{tabular}
    \caption{Table showing values of $k_x/k_{Ny}$ at which error in resolving power of gradient operator as measured by the metric $e_1$ becomes larger than certain values. Schemes correspond to those given in Table \ref{tab:stencils} ((a)-(d) are $2^{nd}$ order, $(e)-(h)$ $4^{th}$ order).}
    \label{tab:grad_rp}
\end{table}

\begin{figure}[t]
    \begin{center}
    \setlength{\unitlength}{1cm}
    \begin{picture}(18,5)(0,0)
    \put(0.3,-0.5){\includegraphics[width=0.33\linewidth]{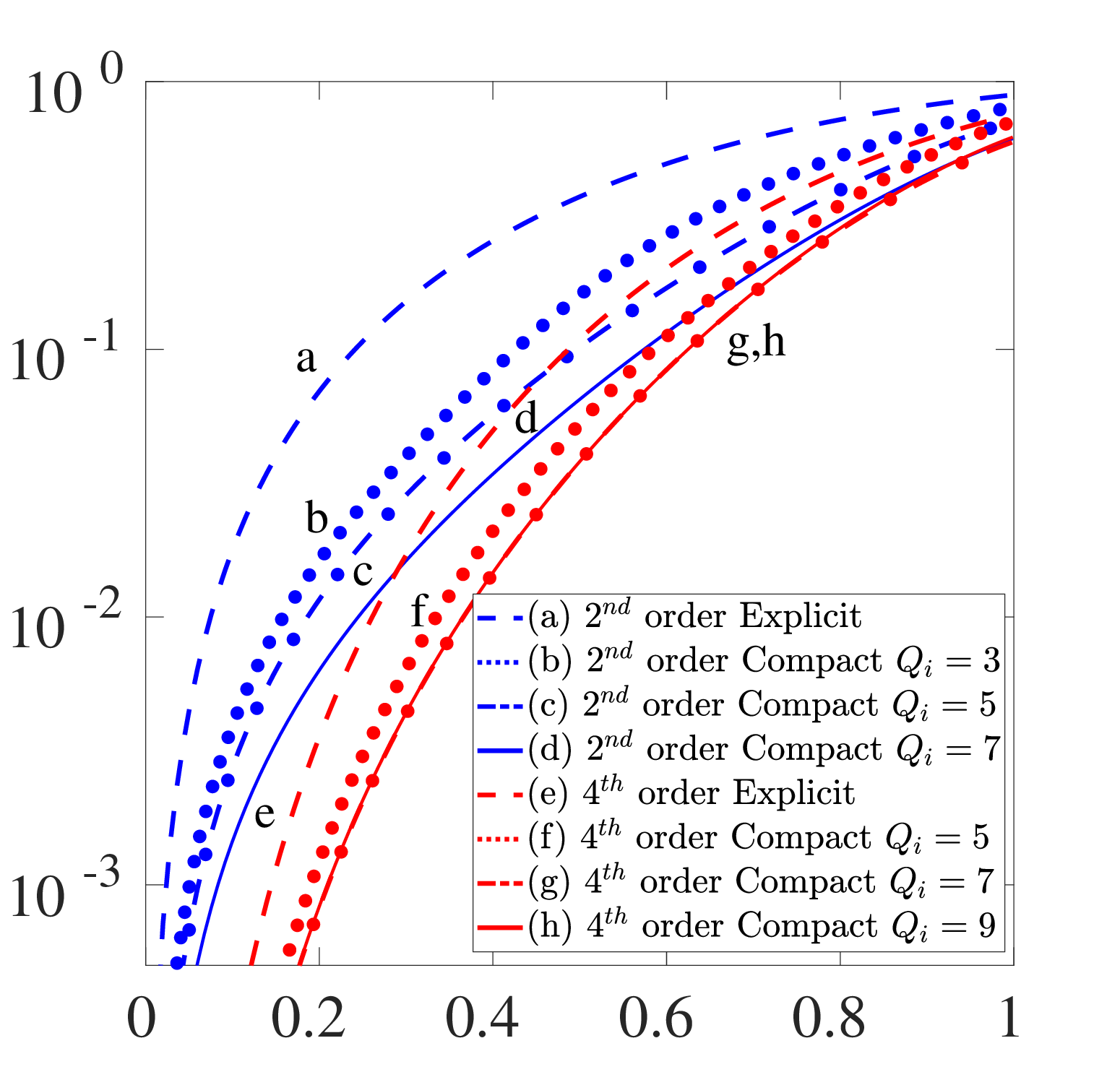}}
    \put(2.9,-0.65){$k_x/k_{Ny}$}
    \put(6.3,-0.5){\includegraphics[width=0.33\linewidth]{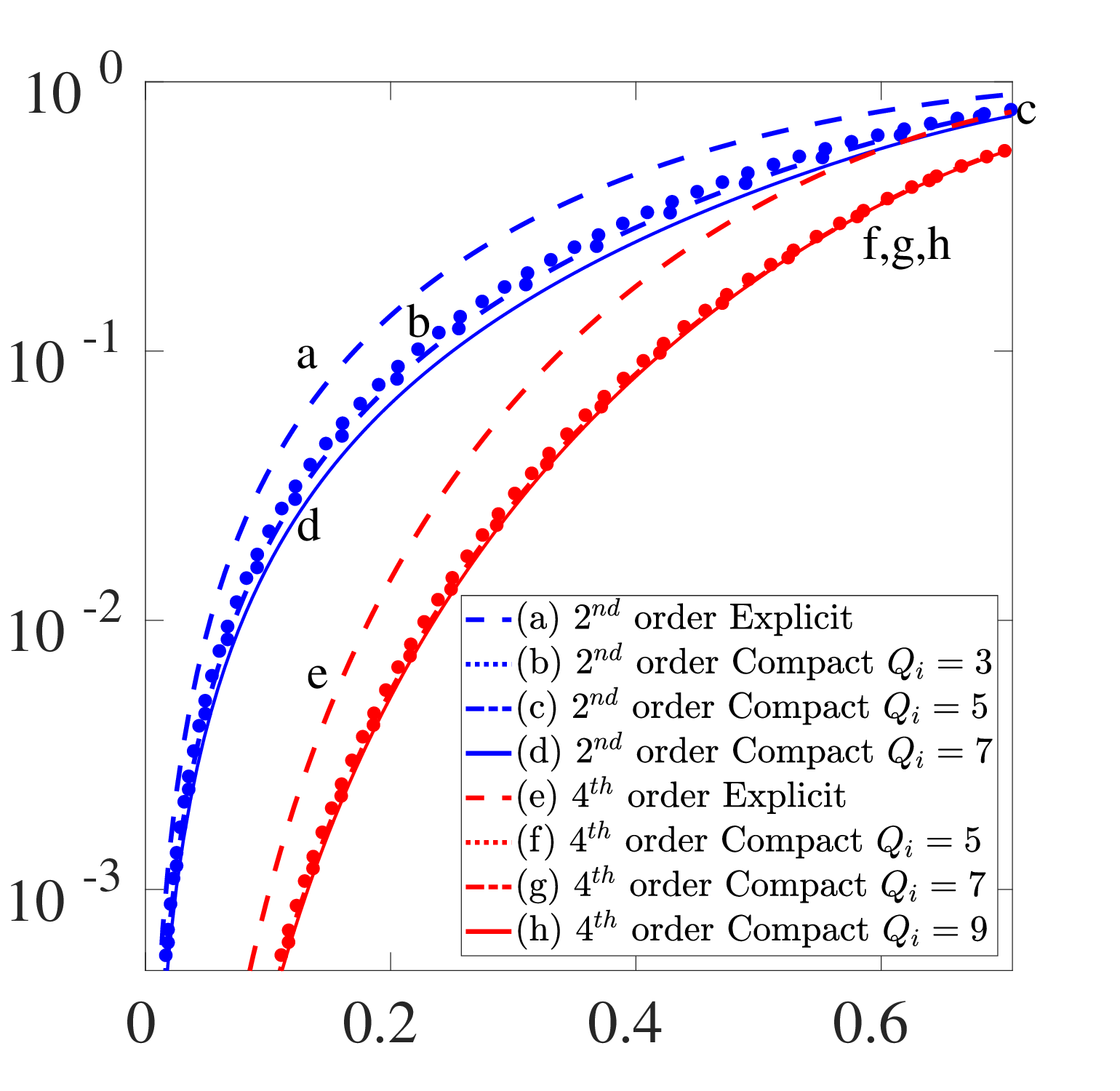}}
    \put(8.95,-0.65){$k_x/k_{Ny}$}
    \put(12.3,-0.5){\includegraphics[width=0.33\linewidth]{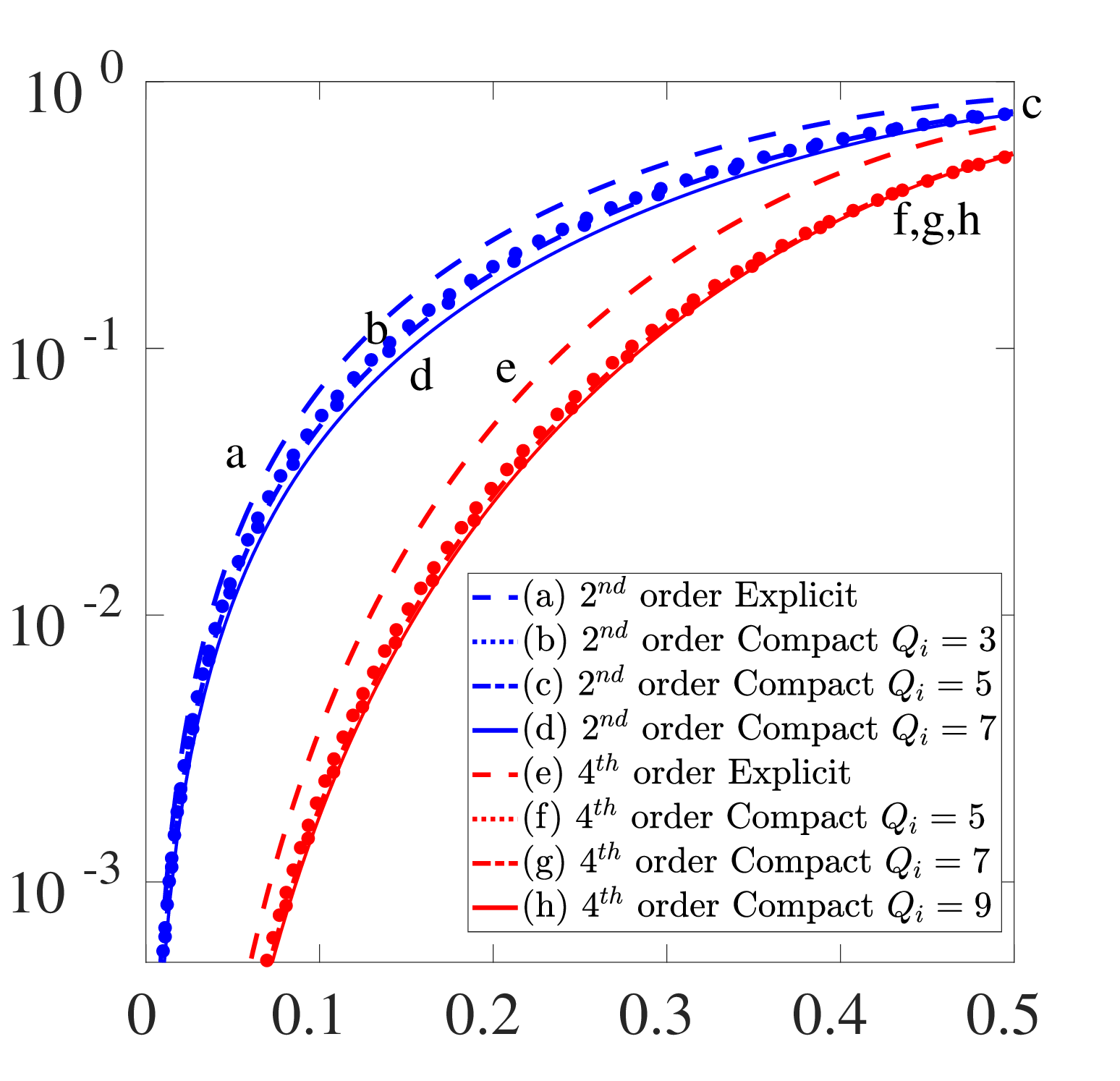}}
    \put(15.0,-0.65){$k_x/k_{Ny}$}
  
\put(2.82, 0.48){\color{red}{\thickhyphen}} 
  \put(3.0, 0.52){\color{red}{{\scalebox{.6}{$\bullet$}}}}

\put(2.82, 0.95){\color{red}{\thickhyphent}} 
\put(2.97, 0.95){\color{red}{\thickhyphent}} 

\put(2.82, 1.41){\color{blue}{\thickhyphen}} 
  \put(3.0, 1.45){\color{blue}{{\scalebox{.6}{$\bullet$}}}}

\put(2.82, 1.88){\color{blue}{\thickhyphent}} 
\put(2.97, 1.88){\color{blue}{\thickhyphent}} 

\put(8.76, 0.49){\color{red}{\thickhyphen}} 
  \put(8.94, 0.53){\color{red}{{\scalebox{.6}{$\bullet$}}}}

\put(8.76, 0.96){\color{red}{\thickhyphent}} 
\put(8.91, 0.96){\color{red}{\thickhyphent}} 

\put(8.76, 1.43){\color{blue}{\thickhyphen}} 
  \put(8.94, 1.47){\color{blue}{{\scalebox{.6}{$\bullet$}}}}

\put(8.76, 1.9){\color{blue}{\thickhyphent}} 
\put(8.91, 1.9){\color{blue}{\thickhyphent}} 

\put(14.795, 0.565){\color{red}{\thickhyphen}} 
  \put(14.975, 0.6){\color{red}{{\scalebox{.6}{$\bullet$}}}}

\put(14.79, 1.025){\color{red}{\thickhyphent}} 
\put(14.94, 1.025){\color{red}{\thickhyphent}} 

\put(14.795, 1.5){\color{blue}{\thickhyphen}} 
  \put(14.975, 1.54){\color{blue}{{\scalebox{.6}{$\bullet$}}}}

\put(14.79, 1.97){\color{blue}{\thickhyphent}} 
\put(14.94, 1.97){\color{blue}{\thickhyphent}}

    \put(-0.1,2.25){\rotatebox{90}{$\varepsilon_1$}}
    \put(6.0,2.25){\rotatebox{90}{$\varepsilon_1$}}
     \put(12.0,2.25){\rotatebox{90}{$\varepsilon_1$}}

    \end{picture}
    \end{center}
    \caption{Error $\varepsilon_1$ of real part of resolving power of the gradient operator for second order accurate (blue lines) and fourth order (red lines) accurate operators, with schemes as labelled in Table \ref{tab:stencils}. Left Panel $k_y=0$, Centre Panel $k_y=k_x$, Right Panel $k_y=2k_x$.}
    \label{fig:grad_phase_error}
\end{figure}
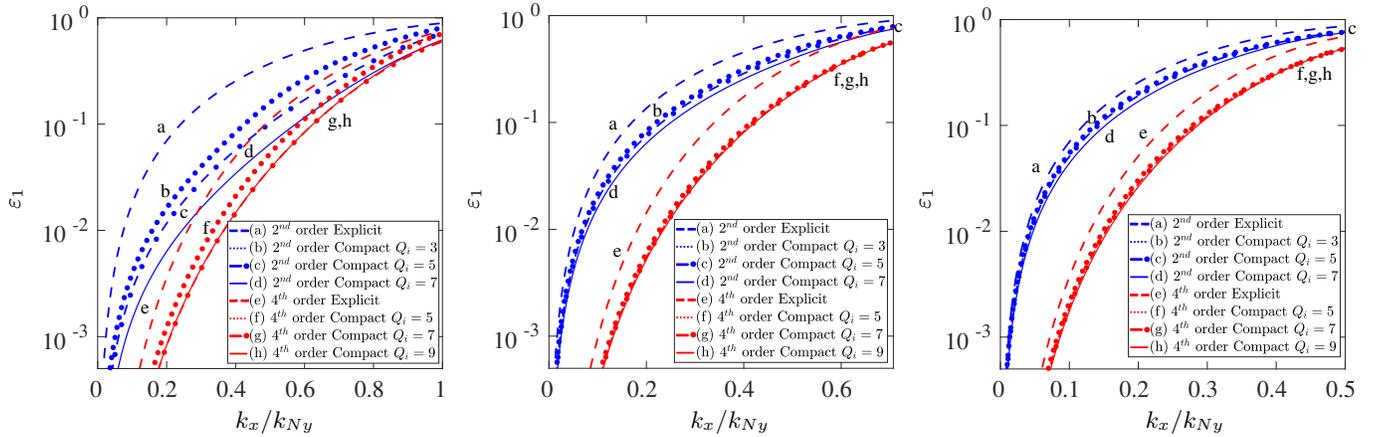

Fig.\ref{fig:grad_phase_error} provides a graphical illustration of the data presented in Table \ref{tab:grad_rp} by showing how $\varepsilon_1$ varies for the different schemes with $k_x$ for the same lines in wavenumber space considered in Fig.\ref{fig:rp_grad_real}. For functions aligned with the $x-$axis we observe a much smaller error, and much larger improvement compared to the explicit cases, for compact schemes.

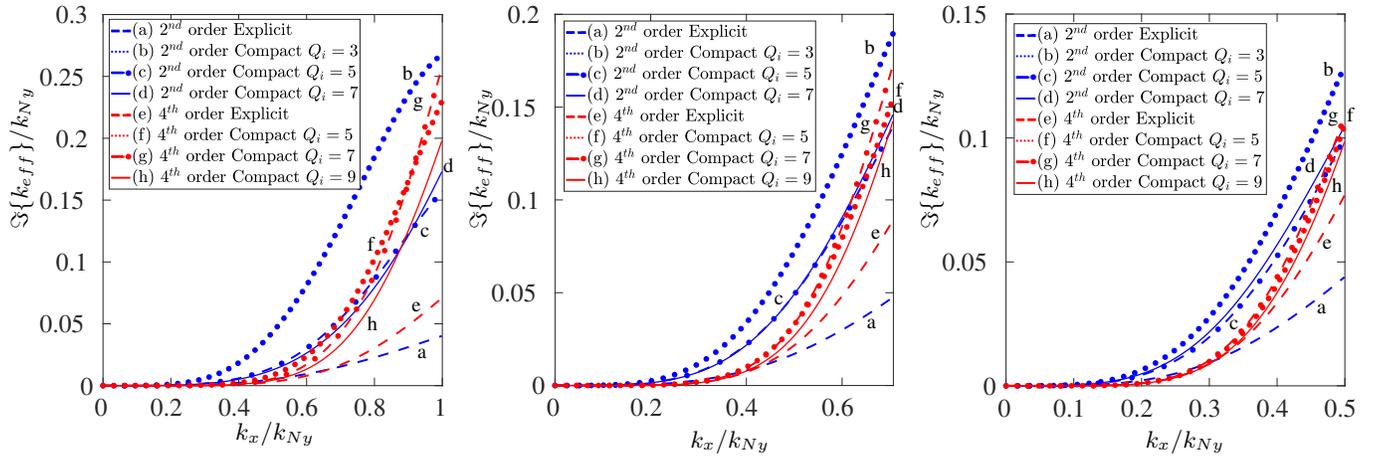
\begin{figure}[t]
    \begin{center}
    \setlength{\unitlength}{1cm}
    \begin{picture}(18,5)(0,0)
    \put(0.3,-0.5){\includegraphics[width=0.33\linewidth]{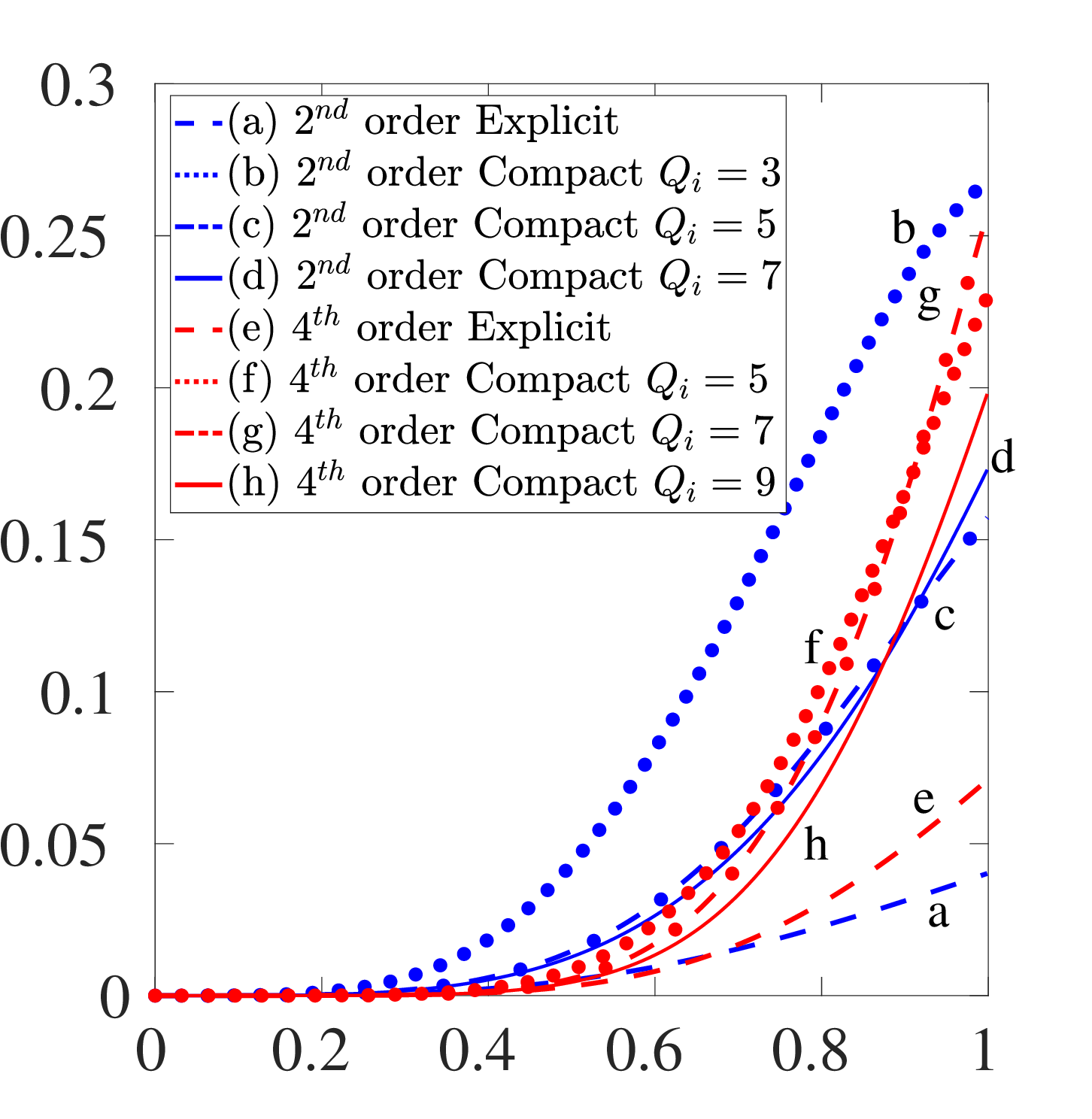}}
    \put(2.9,-0.55){$k_x/k_{Ny}$}
    \put(6.3,-0.5){\includegraphics[width=0.33\linewidth]{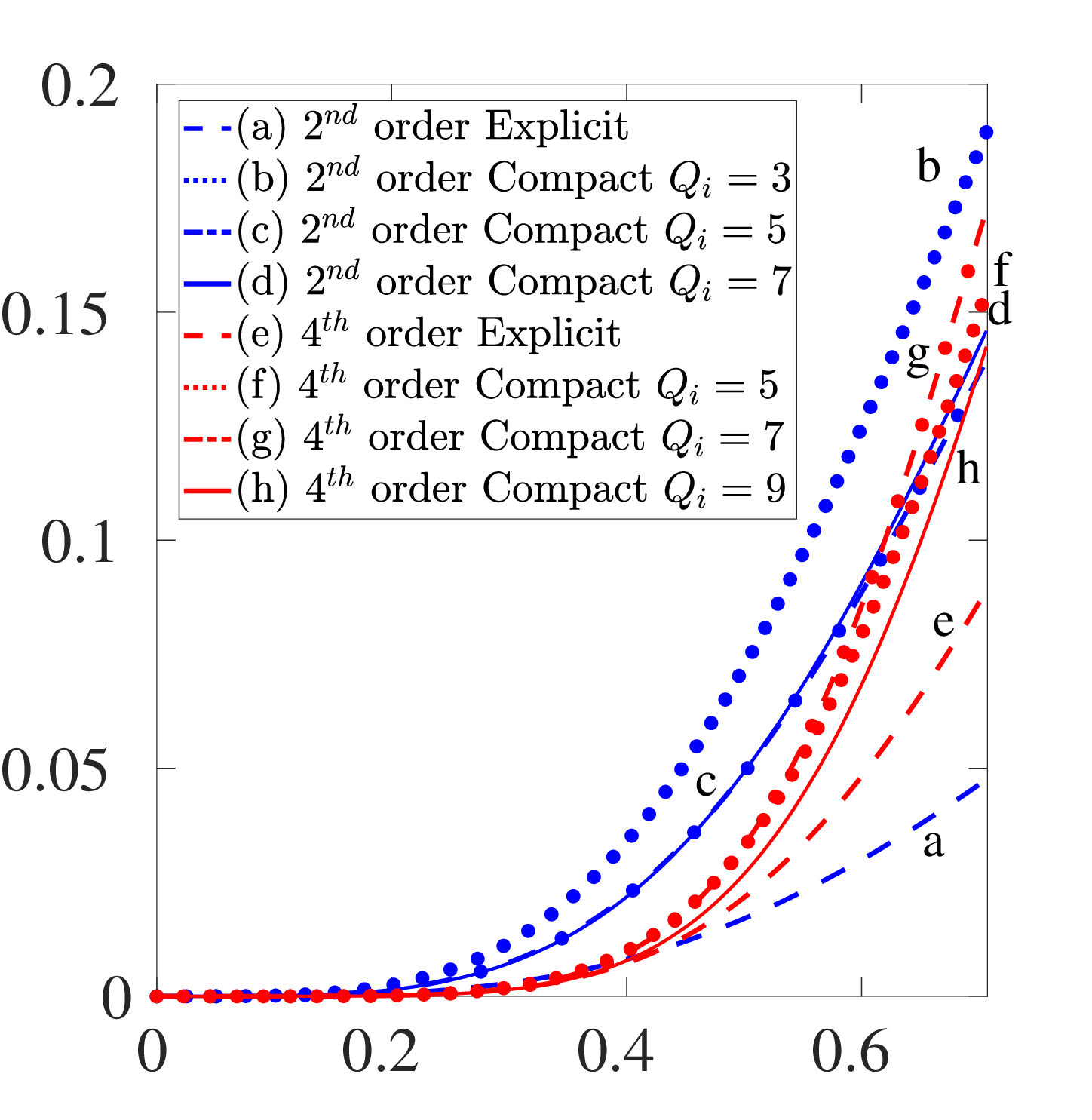}}
    \put(8.95,-0.65){$k_x/k_{Ny}$}
    \put(12.3,-0.5){\includegraphics[width=0.33\linewidth]{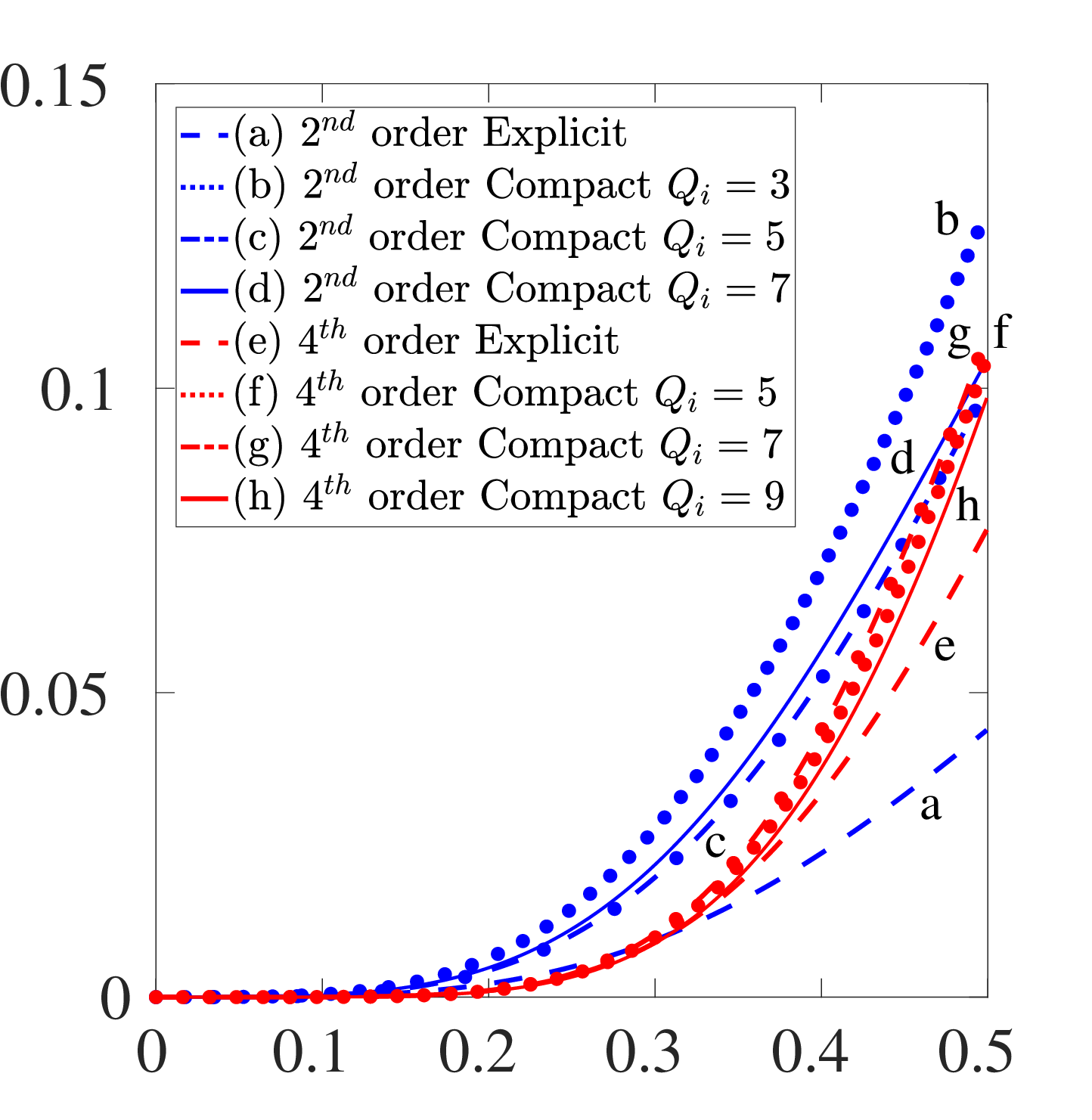}}
    \put(15.0,-0.65){$k_x/k_{Ny}$}
  
\put(1.25, 3.13){\color{red}{\thickhyphen}} 
  \put(1.43, 3.17){\color{red}{{\scalebox{.6}{$\bullet$}}}}

\put(1.25, 3.68){\color{red}{\thickhyphent}} 
\put(1.4, 3.68){\color{red}{\thickhyphent}} 

\put(1.25, 4.245){\color{blue}{\thickhyphen}} 
  \put(1.43, 4.285){\color{blue}{{\scalebox{.6}{$\bullet$}}}}

\put(1.25, 4.8){\color{blue}{\thickhyphent}} 
\put(1.4, 4.8){\color{blue}{\thickhyphent}} 

\put(7.3, 3.10){\color{red}{\thickhyphen}} 
  \put(7.48, 3.14){\color{red}{{\scalebox{.6}{$\bullet$}}}}

\put(7.3, 3.66){\color{red}{\thickhyphent}} 
\put(7.45, 3.66){\color{red}{\thickhyphent}} 

\put(7.3, 4.225){\color{blue}{\thickhyphen}} 
  \put(7.48, 4.255){\color{blue}{{\scalebox{.6}{$\bullet$}}}}

\put(7.3, 4.78){\color{blue}{\thickhyphent}} 
\put(7.45, 4.78){\color{blue}{\thickhyphent}} 

\put(13.28, 3.06){\color{red}{\thickhyphen}} 
  \put(13.46, 3.1){\color{red}{{\scalebox{.6}{$\bullet$}}}}

\put(13.28, 3.62){\color{red}{\thickhyphent}} 
\put(13.43, 3.62){\color{red}{\thickhyphent}} 

\put(13.28, 4.185){\color{blue}{\thickhyphen}} 
  \put(13.46, 4.22){\color{blue}{{\scalebox{.6}{$\bullet$}}}}

\put(13.28, 4.74){\color{blue}{\thickhyphent}} 
\put(13.43, 4.74){\color{blue}{\thickhyphent}}

    \put(-0.1,2.25){\rotatebox{90}{$\Im\{k_{eff}\}/k_{Ny}$}}
    \put(6.0,2.25){\rotatebox{90}{$\Im\{k_{eff}\}/k_{Ny}$}}
     \put(12.0,2.25){\rotatebox{90}{$\Im\{k_{eff}\}/k_{Ny}$}}

    \end{picture}
    \end{center}
    \caption{Imaginary part of the resolving power of the gradient operator for second order accurate (blue lines) and fourth order (red lines) accurate operators, with schemes as labelled in Table \ref{tab:stencils}. Left Panel $k_y=0$, Centre Panel $k_y=k_x$, Right Panel $k_y=2k_x$.}
    \label{fig:rp_grad_im}
\end{figure}

Thus far, we have considered only the real part of the effective wavenumber. Fig.\ref{fig:rp_grad_im} shows the imaginary part of $k_{eff}$ (note that for exact differentiation we desire $\Im\{k_{eff}\}\equiv 0$). Compact schemes typically perform worse than their explicit counterparts, though the increase in error remains comparatively small up to moderate wavenumbers for all schemes except (b). There is little difference between the other two $2^{nd}$ order compact schemes (c) \& (d), however (h) has a noticeably smaller error than both (f) \& (g), indicating there are advantages to use this scheme despite $\Re\{k_{eff}\}$ being graphically identical for (g) \& (h). The effect of this on gradient approximations of functions will be seen in \S\ref{subsec:convergence}.

\begin{figure}[t]
    \begin{center}
    \setlength{\unitlength}{1cm}
    \begin{picture}(18,5)(0,0)
    \put(0.3,-0.5){\includegraphics[width=0.33\linewidth]{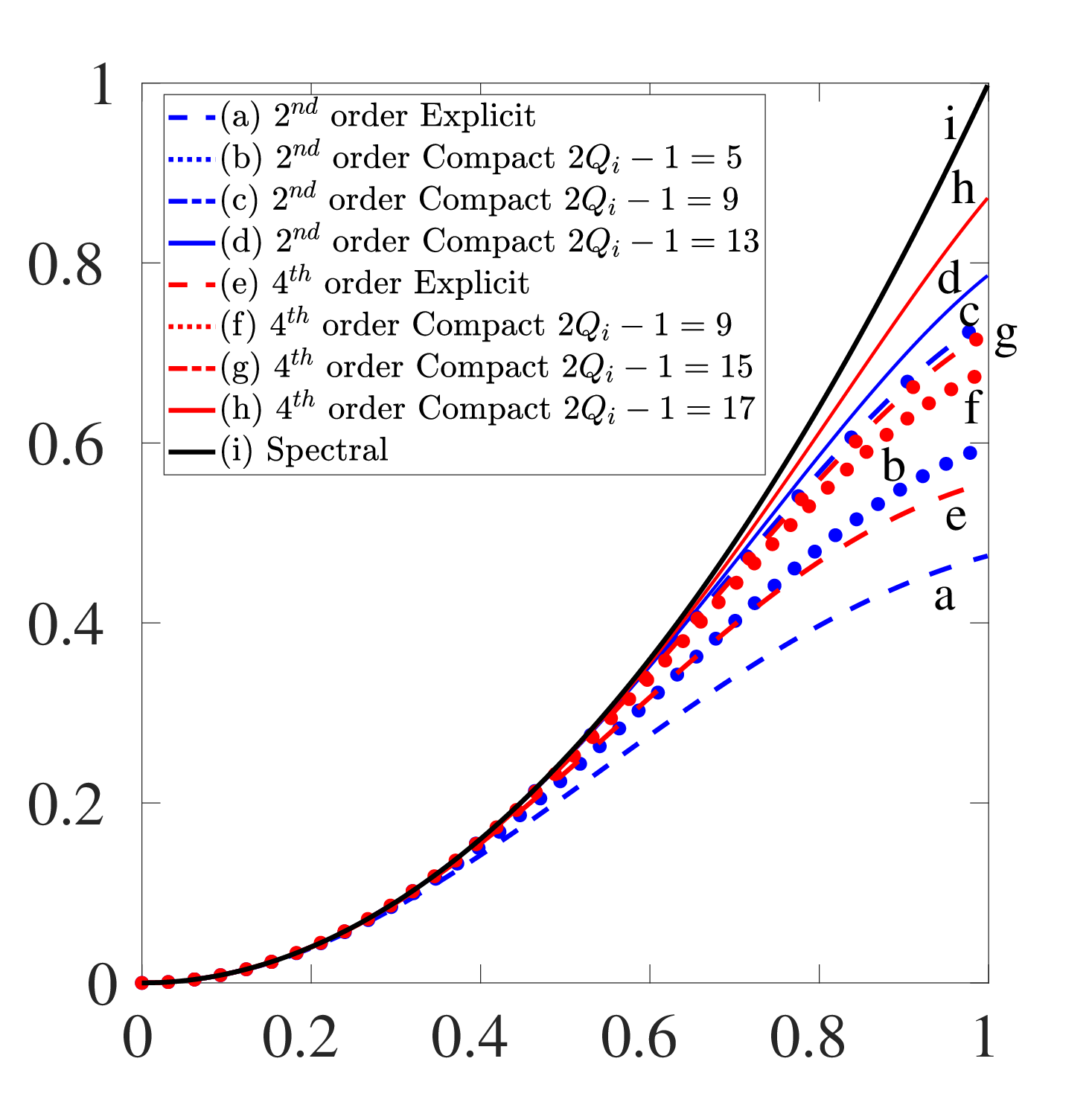}}
    \put(2.9,-0.65){$q/k_{Ny}$}
    \put(6.3,-0.5){\includegraphics[width=0.33\linewidth]{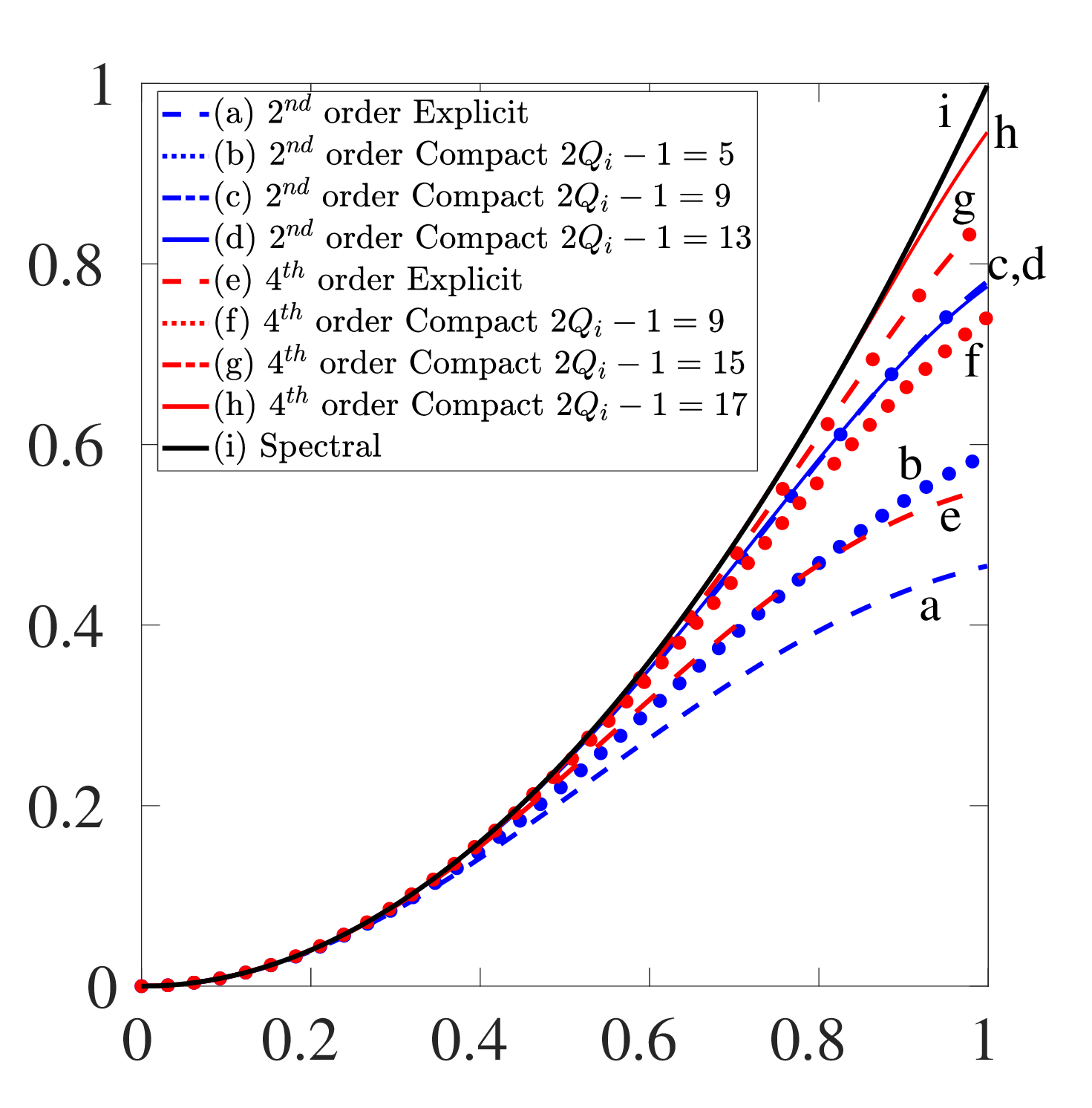}}
    \put(8.95,-0.65){$q/k_{Ny}$}
    \put(12.3,-0.5){\includegraphics[width=0.33\linewidth]{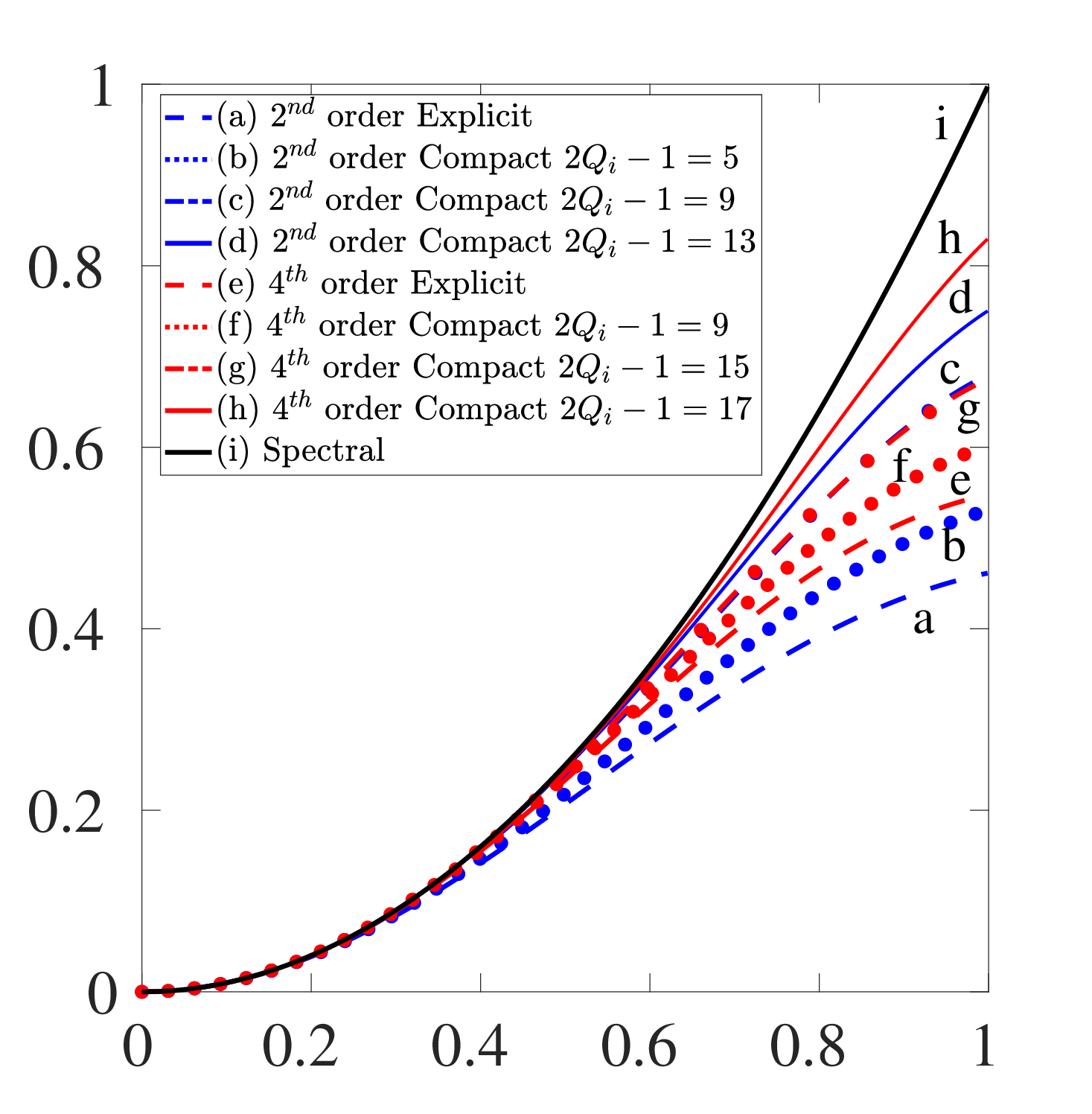}}
    \put(15.0,-0.65){$q/k_{Ny}$}

\put(1.22, 3.407){\color{red}{\thickhyphen}} 
  \put(1.4, 3.45){\color{red}{{\scalebox{.6}{$\bullet$}}}}

\put(1.21, 3.86){\color{red}{\thickhyphent}} 
\put(1.36, 3.86){\color{red}{\thickhyphent}} 

\put(1.22, 4.307){\color{blue}{\thickhyphen}} 
  \put(1.4, 4.35){\color{blue}{{\scalebox{.6}{$\bullet$}}}}

\put(1.21, 4.77){\color{blue}{\thickhyphent}} 
\put(1.36, 4.77){\color{blue}{\thickhyphent}} 

\put(7.18, 3.44){\color{red}{\thickhyphen}} 
  \put(7.36, 3.475){\color{red}{{\scalebox{.6}{$\bullet$}}}}

\put(7.18, 3.89){\color{red}{\thickhyphent}} 
\put(7.33, 3.89){\color{red}{\thickhyphent}} 

\put(7.18, 4.34){\color{blue}{\thickhyphen}} 
  \put(7.36, 4.375){\color{blue}{{\scalebox{.6}{$\bullet$}}}}

\put(7.18, 4.8){\color{blue}{\thickhyphent}} 
\put(7.33, 4.8){\color{blue}{\thickhyphent}} 

\put(13.19, 3.45){\color{red}{\thickhyphen}} 
  \put(13.37, 3.49){\color{red}{{\scalebox{.6}{$\bullet$}}}}

\put(13.19, 3.91){\color{red}{\thickhyphent}} 
\put(13.34, 3.91){\color{red}{\thickhyphent}} 

\put(13.19, 4.355){\color{blue}{\thickhyphen}} 
  \put(13.37, 4.40){\color{blue}{{\scalebox{.6}{$\bullet$}}}}

\put(13.19, 4.81){\color{blue}{\thickhyphent}} 
\put(13.34, 4.81){\color{blue}{\thickhyphent}} 
  
    \put(-0.1,2.25){\rotatebox{90}{$\Re\{{q}_{eff}^2\}/k_{Ny}^2$}}
    \put(6.0,2.25){\rotatebox{90}{$\Re\{{q}_{eff}^2\}/k_{Ny}^2$}}
     \put(12.0,2.25){\rotatebox{90}{$\Re\{{q}_{eff}^2\}/k_{Ny}^2$}}

    \end{picture}
    \end{center}
    \caption{Resolving power of the Laplacian operator for second order accurate (blue lines) and fourth order (red lines) accurate operators, with schemes as labelled in Table \ref{tab:stencils} except for the spectral scheme (i) (exact differentiation). Left Panel $k_y=0$, Centre Panel $k_y=k_x$, Right Panel $k_x=0$.}
    \label{fig:rp_lap_real}
\end{figure}

We now turn to Laplace operators. Fig.\ref{fig:rp_lap_real} shows how the real part of the resolving power $\Re\{q_{eff}^2\}$ varies with $q/k_{Ny}\coloneqq |\boldsymbol{k}|/k_{Ny}$ for various lines in wavenumber space. All compact operators show significant improvement across all three lines in wavenumber space compared to explicit schemes. The implicit operators remain close to spectral accuracy for a much larger proportion of wavenumbers, with the $4^{th}$ order scheme with $2Q_i-1=17$ performing particularly well. It is interesting to note that the $2^{nd}$ order schemes (c) \& (d) actually outperform all $4^{th}$ order compact schemes save (h) at high wavenumber along the lines $k_x=0, k_y=0$. This is likely due to the prescribed form of the implicit coefficients defined in \eqref{eq:implicit_coeffs} and the optimisation procedure outlined in \S\ref{subsec:Stencil_choice} which may be open to improvement for $4^{th}$ order schemes. In addition, scheme (c) is closer to exact differentiation than (d) along $k_y=k_x$ at high wavenumbers, thus larger compact stencils do not always result in improved resolution characteristics.

To quantify the improvement in the real part of the effective wavenumber by compact schemes we introduce the metric
\begin{equation}
    \varepsilon_2(k_x,k_y)=\frac{\left|\Re\{q_{eff}^2(k_x,k_y)\}\right|}{q^2}.
\end{equation}
As was the case for $\varepsilon_1$ when considering the gradient operator, we wish for $\varepsilon_2$ to remain as close to zero as possible as $|\boldsymbol{k}|$ increases. Table \ref{tab:lap_rp} shows values of $q/k_{Ny}$ at which the error in the real part of the effective wavenumber becomes greater than $0.1\%, 1\%$, and $10\%$ respectively along the three lines in wavenumber space show in Fig. \ref{fig:rp_lap_real}. Viewed in this manner, it clear that there is a vast improvement provided by compact LABFM even for low wavenumbers. For $2^{nd}$ order schemes, the percentage of wavenumbers with an error under $0.1\%$ is increased by a factor of at least six when using scheme (d) compared to explicit LABFM. The error only reaches $10\%$ when the magnitude of the wavenumber is approximately $\frac{4}{5}k_{Ny}$. For $4^{th}$ order schemes, although there is improvement across all three lines considered the resolving efficiency along $k_y=k_x$ is most striking. Here, for scheme (h), there is less than $0.1\%$ error up $q=\frac{5}{6}k_{Ny}$, and there is a smaller than $10\%$ error for all wavenumbers to up to $k_{Ny}$. The resolution characteristics of this scheme appear spectral-like. 

\begin{table}[t]
    \centering
    \begin{tabular}{ *{10}{c|} }
  & \multicolumn{3}{c|}{$\varepsilon_2=0.001$}
            & \multicolumn{3}{c|}{$\varepsilon_2=0.01$}
             & \multicolumn{3}{c|}{$\varepsilon_2=0.1$}
            \\

   \quad Numerical Scheme \quad & \quad   $k_y=0$ \quad   & \quad $k_y=k_x$  \quad & \quad  $k_x=0$ \quad & \quad   $k_y=0$ \quad   & \quad $k_y=k_x$  \quad & \quad  $k_x=0$ \quad &\quad   $k_y=0$ \quad   & \quad $k_y=k_x$  \quad & \quad  $k_x=0$ \quad \\
    \hline
 (a) &0.04  &0.04  &0.04  &0.12  &0.12  &0.12  &0.38 &0.38 &0.37\\
 (b) &0.06  &0.05  &0.04  &0.18  &0.16  &0.14 &0.54 &0.49 &0.45 \\
 (c) &0.31  & 0.17 &0.10  &0.46  & 0.42 &0.29 &0.77 &0.82 & 0.68 \\
 (d) &0.37  &0.33  &0.24  & 0.51 &0.50  &0.44 &0.83 &0.82 &0.79 \\
 (e) & 0.16 &0.16  &0.16  &0.30  &0.30  &0.30 &0.57 &0.57 &0.57 \\
 (f) & 0.21 &0.22  &0.18  &0.38  &0.40  &0.32 &0.71 &0.75 &0.61 \\
 (g) &0.23  &0.28  &0.21  &0.42  & 0.52 &0.37 &0.75 &0.93 &0.69 \\
 (h) & 0.32 &0.83  &0.29  &0.57  &0.89  &0.51 &0.95 & --- &0.83 \\
    \end{tabular}
    \caption{Table showing values at which error in resolving power of the Laplacian operator as measured by the metric $\varepsilon_2$ becomes larger than certain values. Schemes correspond to those given in Table \ref{tab:stencils} ((a)-(d) are $2^{nd}$ order, $(e)-(h)$ $4^{th}$ order).}
    \label{tab:lap_rp}
\end{table}

\begin{figure}[t]
    \begin{center}
    \setlength{\unitlength}{1cm}
    \begin{picture}(18,5)(0,0)
    \put(0.3,-0.5){\includegraphics[width=0.33\linewidth]{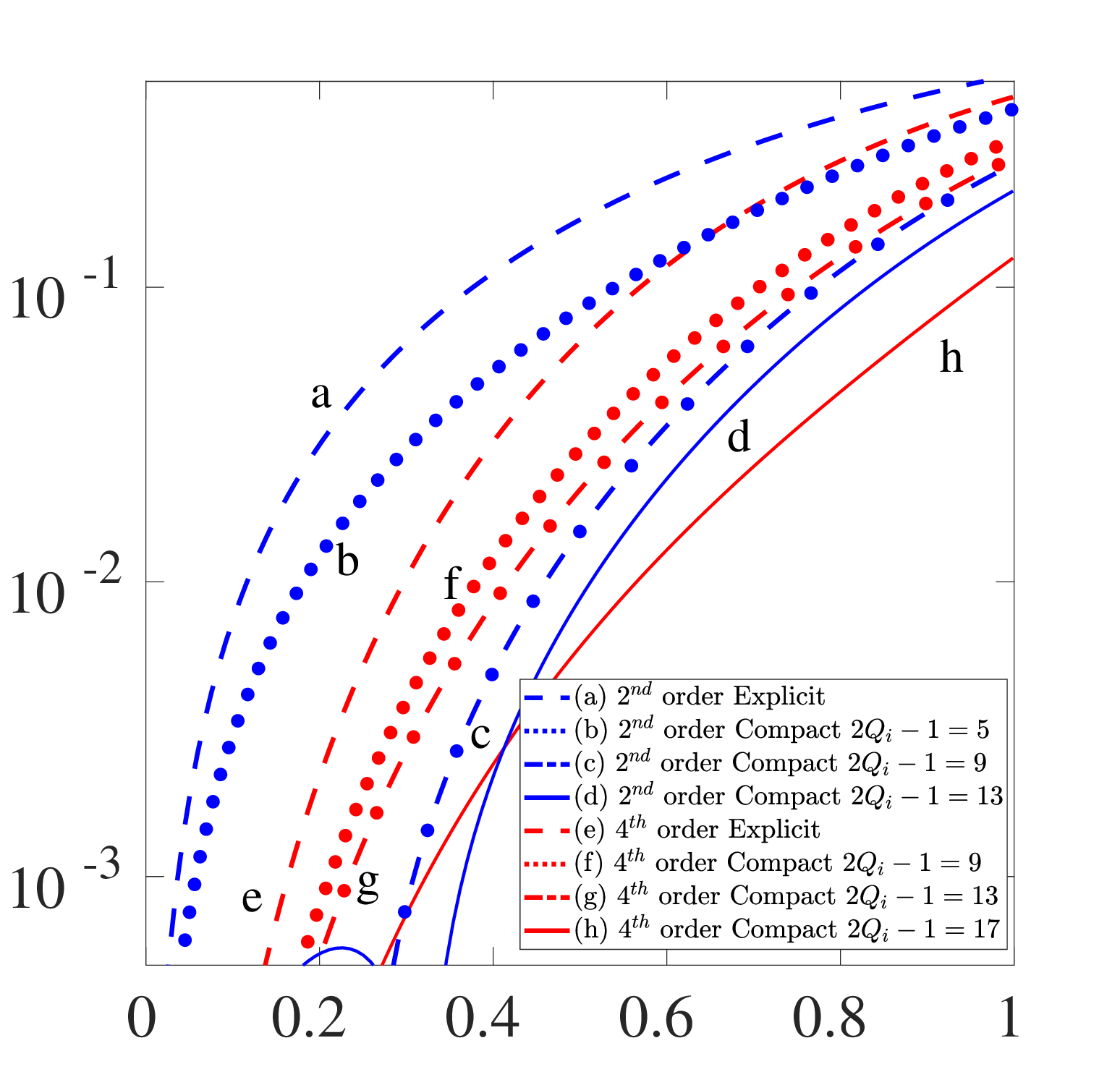}}
    \put(2.9,-0.65){$q/k_{Ny}$}
    \put(6.3,-0.5){\includegraphics[width=0.33\linewidth]{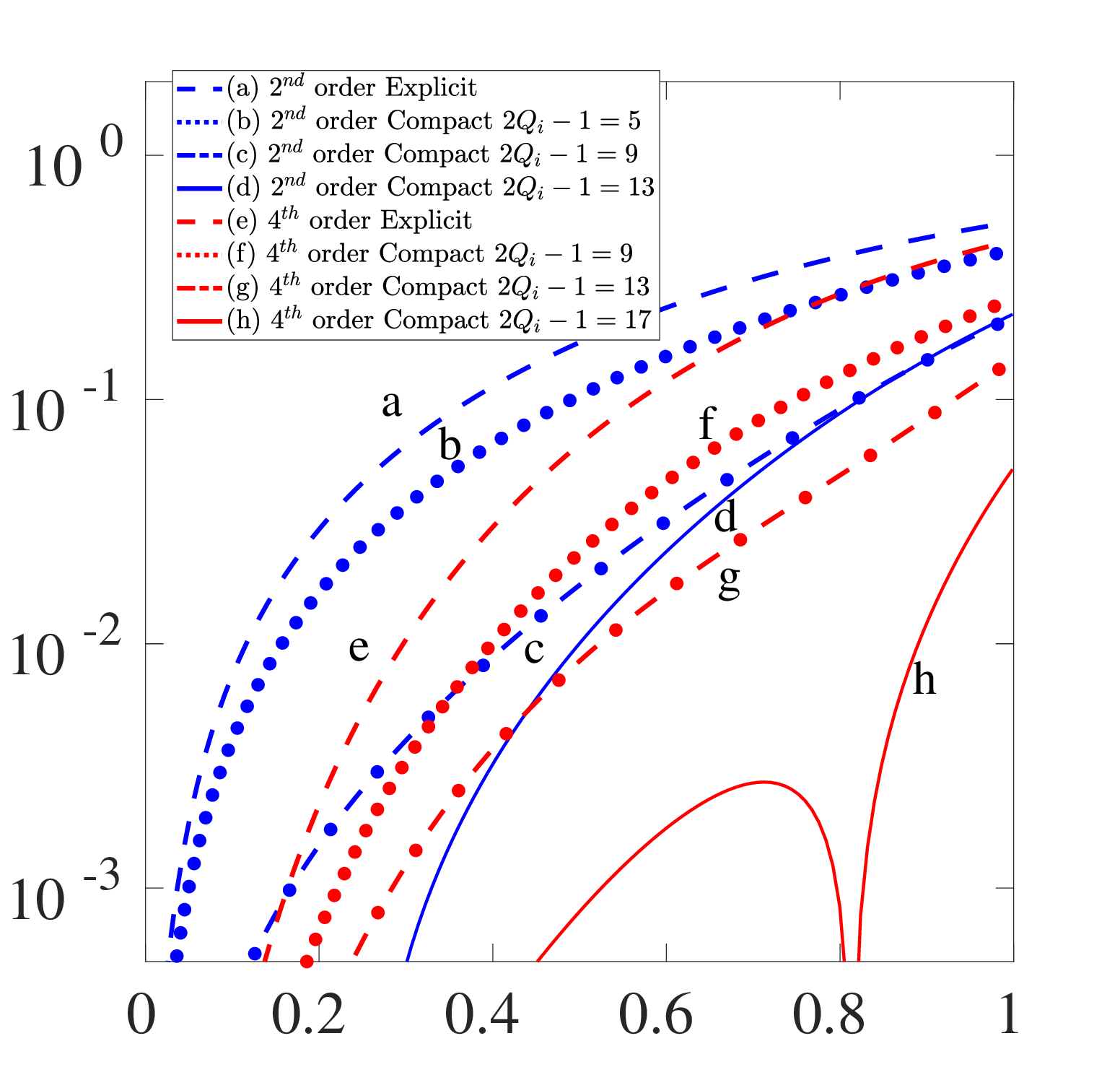}}
    \put(8.95,-0.65){$q/k_{Ny}$}
    \put(12.3,-0.5){\includegraphics[width=0.33\linewidth]{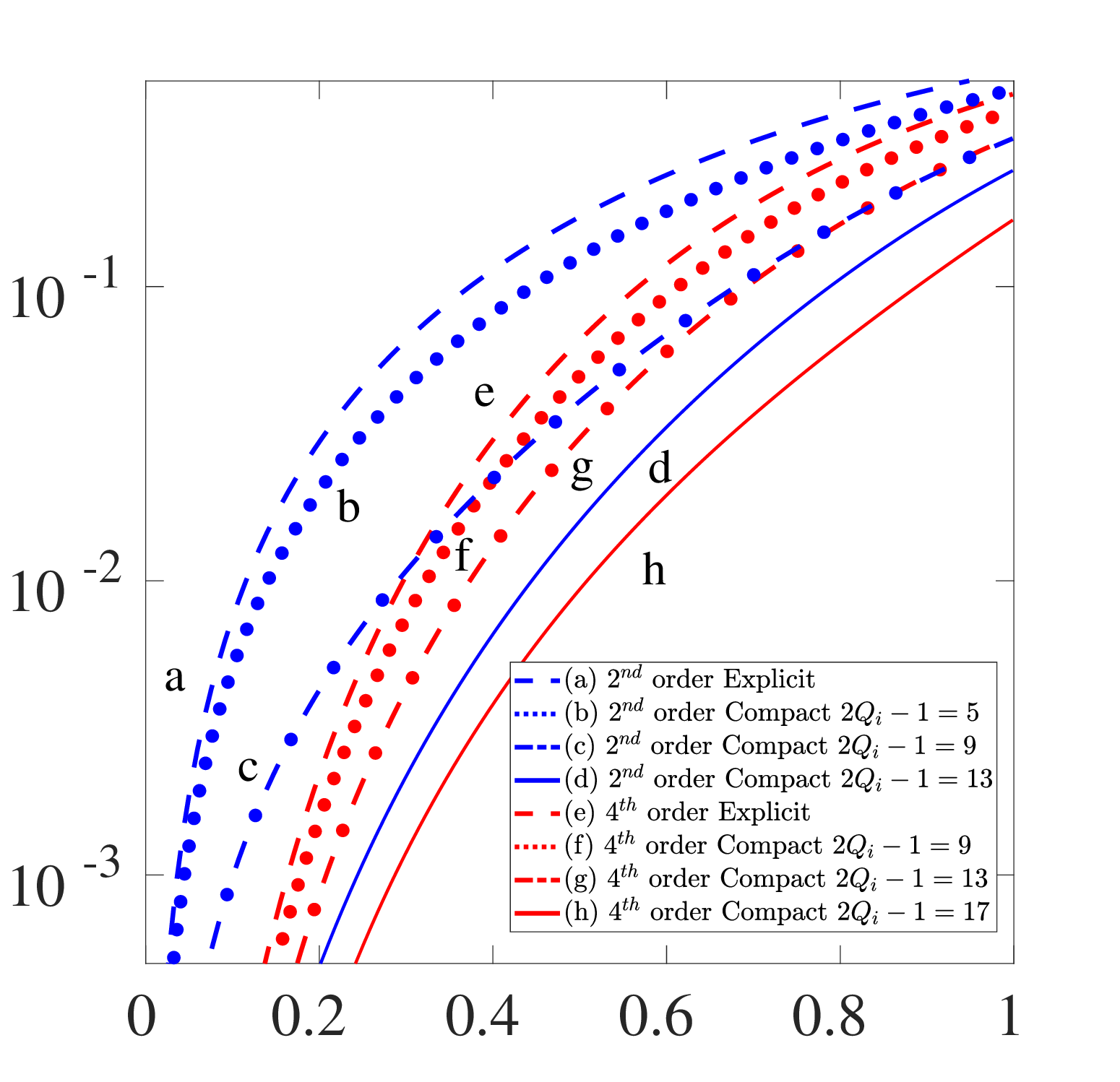}}
    \put(15.0,-0.65){$q/k_{Ny}$}
  
    \put(-0.1,2.25){\rotatebox{90}{$ \varepsilon_2$}}
    \put(6.0,2.25){\rotatebox{90}{$ \varepsilon_2$}}
     \put(12.0,2.25){\rotatebox{90}{$ \varepsilon_2$}}

\put(3.06, 0.405){\color{red}{\thickhyphen}} 
  \put(3.24, 0.44){\color{red}{{\scalebox{.6}{$\bullet$}}}}

\put(3.06, 0.76){\color{red}{\thickhyphent}} 
\put(3.21, 0.76){\color{red}{\thickhyphent}} 

\put(3.06, 1.11){\color{blue}{\thickhyphen}} 
  \put(3.24, 1.15){\color{blue}{{\scalebox{.6}{$\bullet$}}}}

\put(3.06, 1.46){\color{blue}{\thickhyphent}} 
\put(3.21, 1.46){\color{blue}{\thickhyphent}} 

\put(7.23, 3.595){\color{red}{\thickhyphen}} 
  \put(7.41, 3.63){\color{red}{{\scalebox{.6}{$\bullet$}}}}

\put(7.23, 3.95){\color{red}{\thickhyphent}} 
\put(7.38, 3.95){\color{red}{\thickhyphent}} 

\put(7.23, 4.292){\color{blue}{\thickhyphen}} 
  \put(7.41, 4.33){\color{blue}{{\scalebox{.6}{$\bullet$}}}}

\put(7.23, 4.65){\color{blue}{\thickhyphent}} 
\put(7.38, 4.65){\color{blue}{\thickhyphent}} 

\put(15.02, 0.475){\color{red}{\thickhyphen}} 
  \put(15.19, 0.52){\color{red}{{\scalebox{.6}{$\bullet$}}}}

\put(15.01, 0.825){\color{red}{\thickhyphent}} 
\put(15.16, 0.825){\color{red}{\thickhyphent}} 

\put(15.02, 1.185){\color{blue}{\thickhyphen}} 
  \put(15.19, 1.22){\color{blue}{{\scalebox{.6}{$\bullet$}}}}

\put(15.01, 1.53){\color{blue}{\thickhyphent}} 
\put(15.16, 1.53){\color{blue}{\thickhyphent}} 

    \end{picture}
    \end{center}
    \caption{Error $\varepsilon_2$ of resolving power of the Laplacian operator for second order accurate (blue lines) and fourth order (red lines) accurate operators, with schemes as labelled in Table \ref{tab:stencils}. Left Panel $k_y=0$, Centre Panel $k_y=k_x$, Right Panel $k_x=0$.}
    \label{fig:lap_phase_error}
\end{figure}

Fig.\ref{fig:lap_phase_error} provides a graphical representation of the error $\varepsilon_2$ of resolving power of compact Laplace operators over a wide range of wavenumbers. Here the improvements provided by the compact formulation are striking, particularly in schemes (c), (d), and (h). The $4^{th}$ order scheme (h) always has an error lower than $20\%$, whereas the explicit $4^{th}$ order scheme passes this threshold at approximately $q/k_{Ny}\approx 0.7$. Similar improvements can be observed between the $2^{nd}$ order explicit scheme (a) and the $2^{nd}$ order compact scheme (d) across all wavenumbers shown. Note that scheme (d) in the left panel (up to $q\approx 0.3 k_{Ny}$) and scheme (h) in the centre panel (up to $q\approx0.8 k_{Ny}$) have regions where $\Re\{q_{eff}^2\}>q^2$, however this remains below the threshold of $0.5\%$ excitation set in the algorithm of \S\ref{subsec:Optimise_coefficients}.

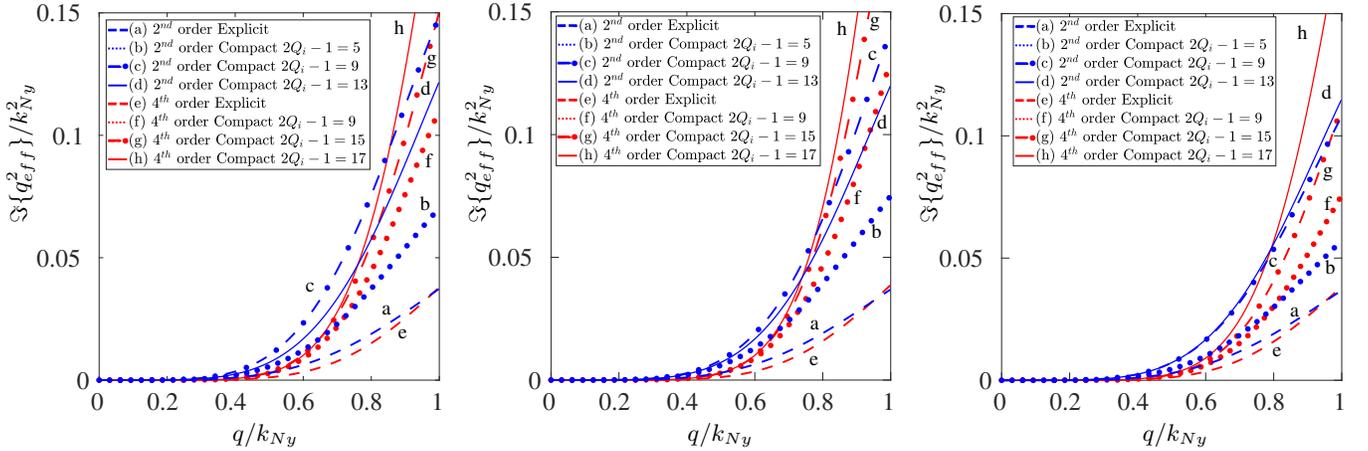
\begin{figure}[t]
    \begin{center}
    \setlength{\unitlength}{1cm}
    \begin{picture}(18,5)(0,0)
    \put(0.3,-0.5){\includegraphics[width=0.33\linewidth]{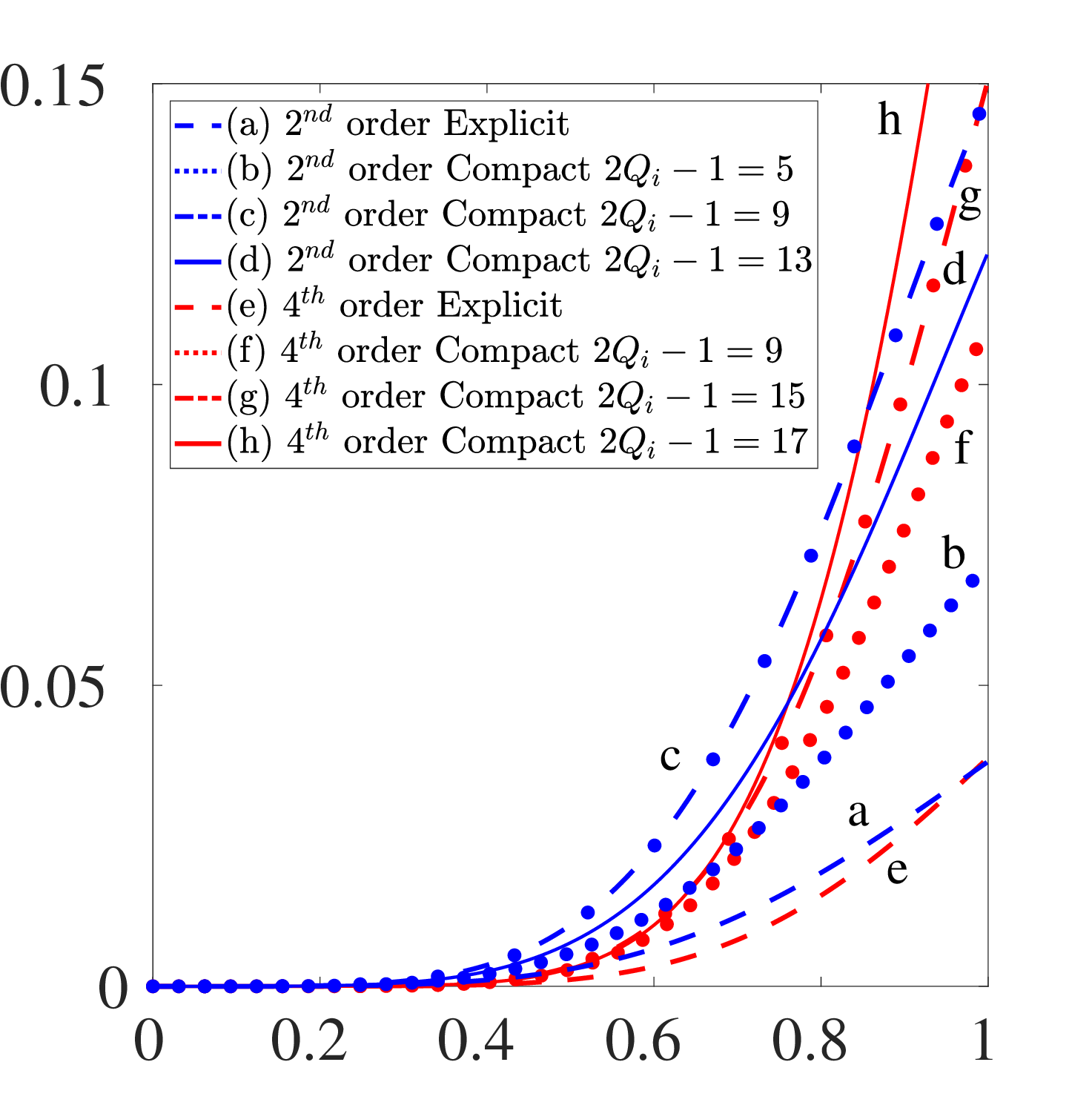}}
    \put(2.9,-0.65){$q/k_{Ny}$}
    \put(6.3,-0.5){\includegraphics[width=0.33\linewidth]{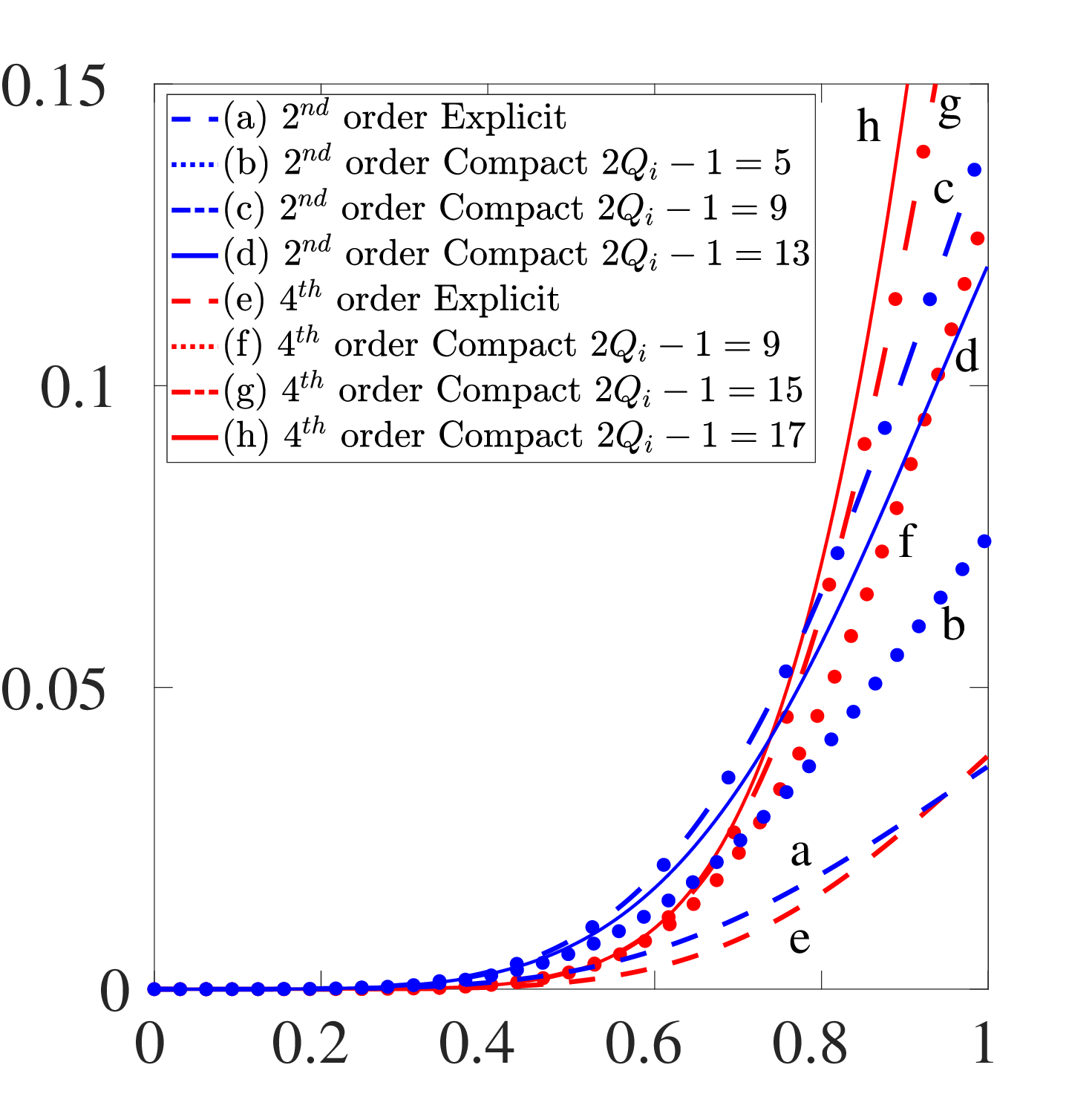}}
    \put(8.95,-0.65){$q/k_{Ny}$}
    \put(12.3,-0.5){\includegraphics[width=0.33\linewidth]{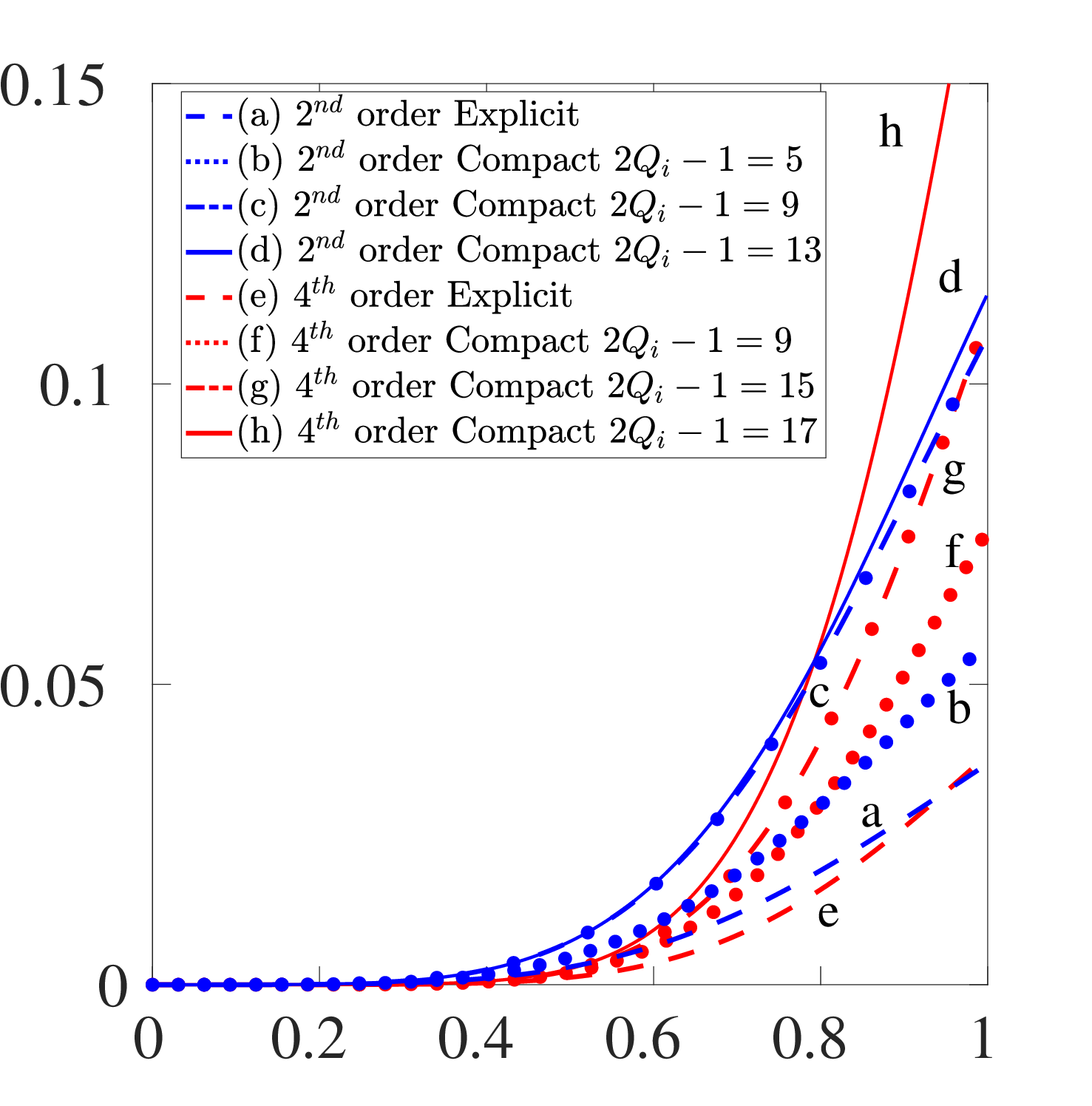}}
    \put(15.0,-0.65){$q/k_{Ny}$}

\put(1.25, 3.26){\color{red}{\thickhyphen}} 
  \put(1.43, 3.3){\color{red}{{\scalebox{.6}{$\bullet$}}}}

\put(1.25, 3.75){\color{red}{\thickhyphent}} 
\put(1.4, 3.75){\color{red}{\thickhyphent}} 

\put(1.25, 4.245){\color{blue}{\thickhyphen}} 
  \put(1.43, 4.285){\color{blue}{{\scalebox{.6}{$\bullet$}}}}

\put(1.25, 4.73){\color{blue}{\thickhyphent}} 
\put(1.4, 4.73){\color{blue}{\thickhyphent}} 

\put(7.22, 3.31){\color{red}{\thickhyphen}} 
  \put(7.4, 3.345){\color{red}{{\scalebox{.6}{$\bullet$}}}}

\put(7.22, 3.8){\color{red}{\thickhyphent}} 
\put(7.37, 3.8){\color{red}{\thickhyphent}} 

\put(7.22, 4.292){\color{blue}{\thickhyphen}} 
  \put(7.4, 4.33){\color{blue}{{\scalebox{.6}{$\bullet$}}}}

\put(7.22, 4.79){\color{blue}{\thickhyphent}} 
\put(7.37, 4.79){\color{blue}{\thickhyphent}} 

\put(13.31, 3.3){\color{red}{\thickhyphen}} 
  \put(13.49, 3.34){\color{red}{{\scalebox{.6}{$\bullet$}}}}

\put(13.31, 3.795){\color{red}{\thickhyphent}} 
\put(13.46, 3.795){\color{red}{\thickhyphent}} 

\put(13.31, 4.285){\color{blue}{\thickhyphen}} 
  \put(13.49, 4.325){\color{blue}{{\scalebox{.6}{$\bullet$}}}}

\put(13.31, 4.78){\color{blue}{\thickhyphent}} 
\put(13.46, 4.78){\color{blue}{\thickhyphent}}

    \put(-0.1,2.25){\rotatebox{90}{$\Im\{q_{eff}^2\}/k_{Ny}^2$}}
    \put(6.0,2.25){\rotatebox{90}{$\Im\{q_{eff}^2\}/k_{Ny}^2$}}
     \put(12.0,2.25){\rotatebox{90}{$\Im\{q_{eff}^2\}/k_{Ny}^2$}}

    \end{picture}
    \end{center}
    \caption{Imaginary part of the resolving power of the Laplacian operator for second order accurate (blue lines) and fourth order (red lines) accurate operators, with schemes as labelled in Table \ref{tab:stencils}. Left Panel $k_y=0$, Centre Panel $k_y=k_x$, Right Panel $k_x=0$.}
    \label{fig:rp_lap_im}
\end{figure}

To complete the resolving power analysis of the Laplacian operator we now turn to the imaginary part of the effective wavenumber. Fig.\ref{fig:rp_lap_im} shows $\Im\{q^2_{eff}\}$ (normalised by the square of the Nyquist wavenumber) for schemes (a)--(h). For moderate and large wavenumbers the compact schemes have a larger error than the explicit schemes and in general the larger implicit stencil the larger the error. However, these errors are smaller than the improvements in the real part of the effective wavenumber, especially at relatively low wavenumbers (i.e. up to $|\boldsymbol{k}|=0.5k_{Ny}$).

In this section we have demonstrated that compact LABFM operators typically remain close to exact differentiation for a much larger range of wavenumbers than explicit LABFM. Compact Laplacian operators in particular remain close to spectral accuracy across a wide range of wavenumbers. However, the imaginary part of the effective wavenumber can often have greater error for compact operators. We now perform further numerical tests on the compact LABFM, using explicit LABFM as a reference, to determine if this is important in operator approximations.

\section{Convergence and Stability Tests}
\label{sec:Convergence_Studies}
\subsection{Convergence}
\label{subsec:convergence}
We first compare the convergence of the implicit operators to explicit operators. We evaluate the accuracy of a given approximation to an operator by using the $L_2$-norm defined through
\begin{equation}
    \label{eq::l2_norm_defn}
    \rnumber(\cdot )= \frac{\left(\sum_{i=1}^N \left[L_{\mathrm{a}}(\cdot)_i-L_{n}(\cdot)_i \right]^2\right)^{1/2}}{\left(\sum_{i=1}^N \left[L_{\mathrm{a}}(\cdot)_i\right]^2\right)^{1/2}}
\end{equation}
where $L_a(\cdot)_i$ is the analytic operator in question, and $L_n(\cdot)_i$ is the numerical approximation. To measure the improvement provided by using compact LABFM we follow \cite{MKLABFM_Arxiv} in defining the ratio $\mathcal{R}$ between the $\rnumber$ of the approximation produced by explicit LABFM (E) to that of compact LABFM (C)
\begin{equation}
    \label{eq::l2_norm_ratio}
    \mathcal{R}\coloneqq \frac{\rnumber(\mathrm{C})}{\rnumber(\mathrm{E})}.
\end{equation}

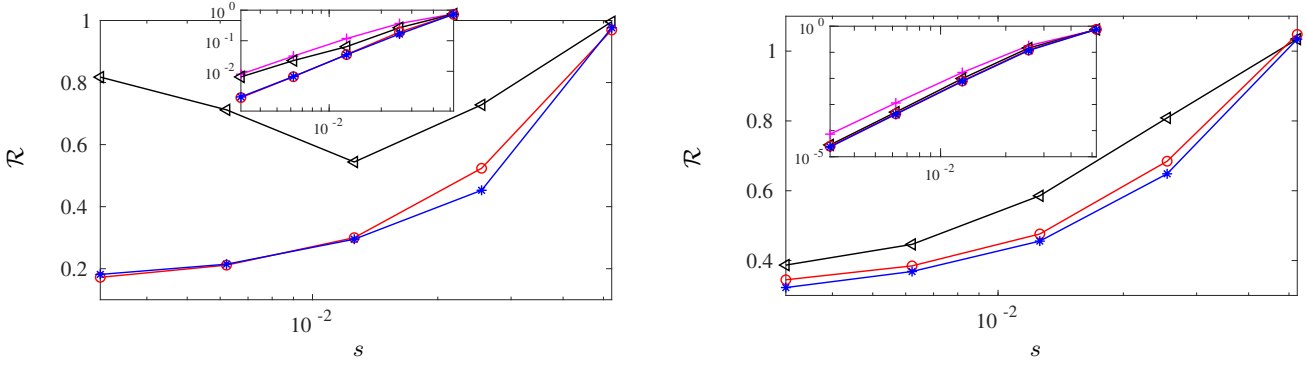
\begin{figure}
    \begin{center}
    \setlength{\unitlength}{1cm}
    \begin{picture}(18,5)(0,0)
    \put(0,0){\includegraphics[width=0.49\linewidth]{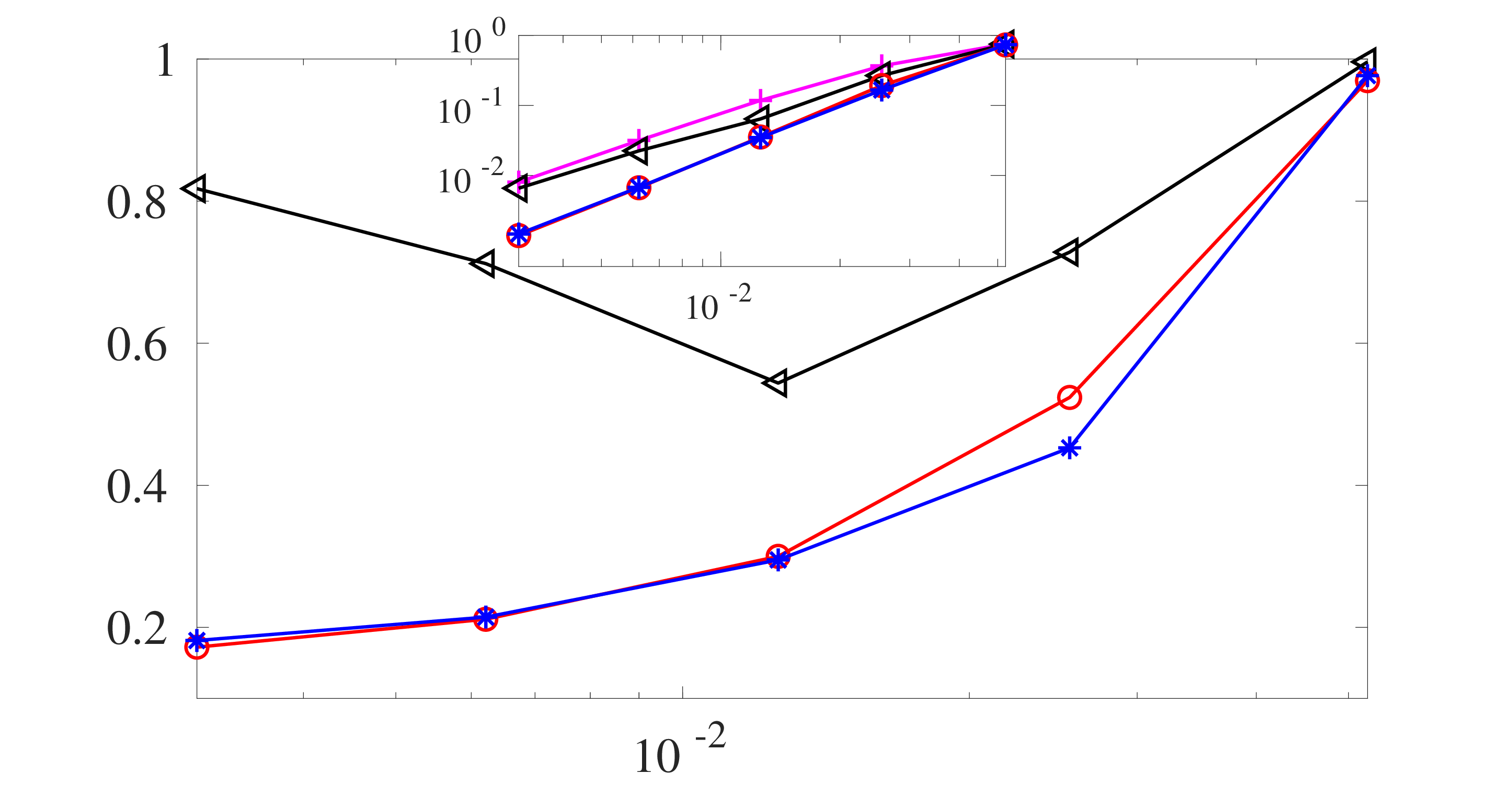}}
    \put(9,0){\includegraphics[width=0.49\linewidth]{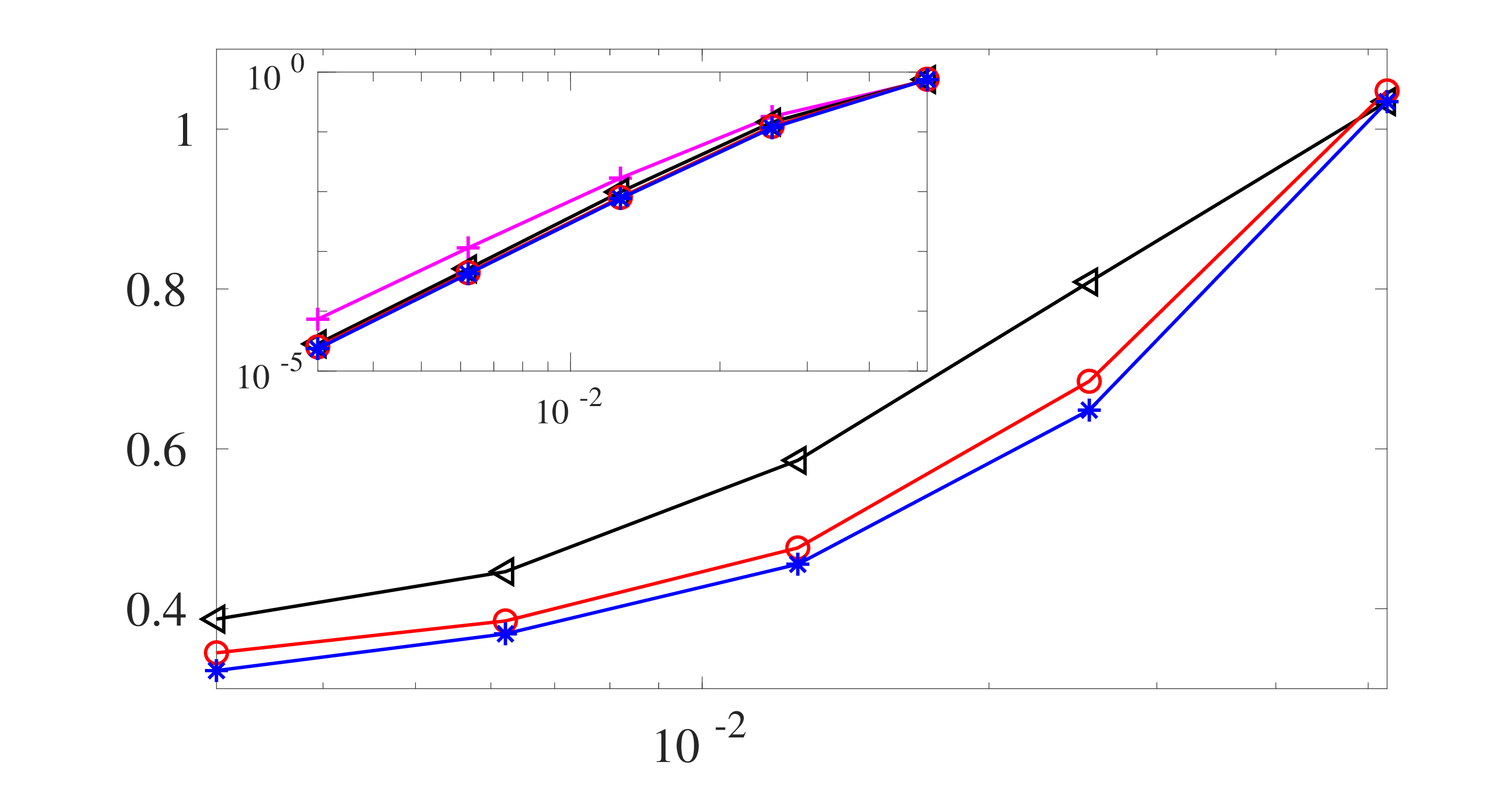}}
    \put(-0.1,2.25){\rotatebox{90}{$\mathcal{R}$}}
   \put(8.9,2.25){\rotatebox{90}{$\mathcal{R}$}}
    \put(4.5,-0.25){$s$}
    \put(13.5,-0.25){$s$}
\end{picture}
    \end{center}
    \caption{Ratio $\mathcal{R}$ of errors arising from explicit and compact approximations to the gradient operator acting on function \eqref{eq::test_function}. Left Panel, $2^{nd}$ order approximations: black $\triangleleft$, $Q_i=3$ (scheme b); red $\circ$ $Q_i=5$ (scheme c); blue $*$, $Q_i=7$ (scheme d). Right Panel, $4^{th}$ order approximations; black $\triangleleft$, $Q_i=5$ (scheme f); red $\circ$ $Q_i=7$ (scheme g); blue $*$, $Q_i=9$ (scheme h). Insets show convergence of $L_2-$ norm with magenta line corresponding to explicit LABFM.}
    \label{fig:grad_conv}
    \end{figure}

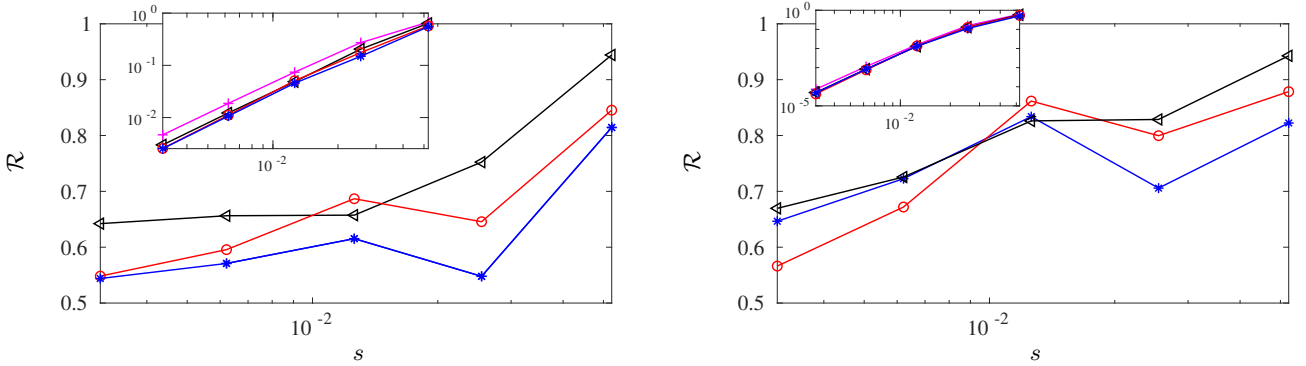
\begin{figure}
    \begin{center}
    \setlength{\unitlength}{1cm}
    \begin{picture}(18,5)(0,0)
    \put(0,0){\includegraphics[width=0.49\linewidth]{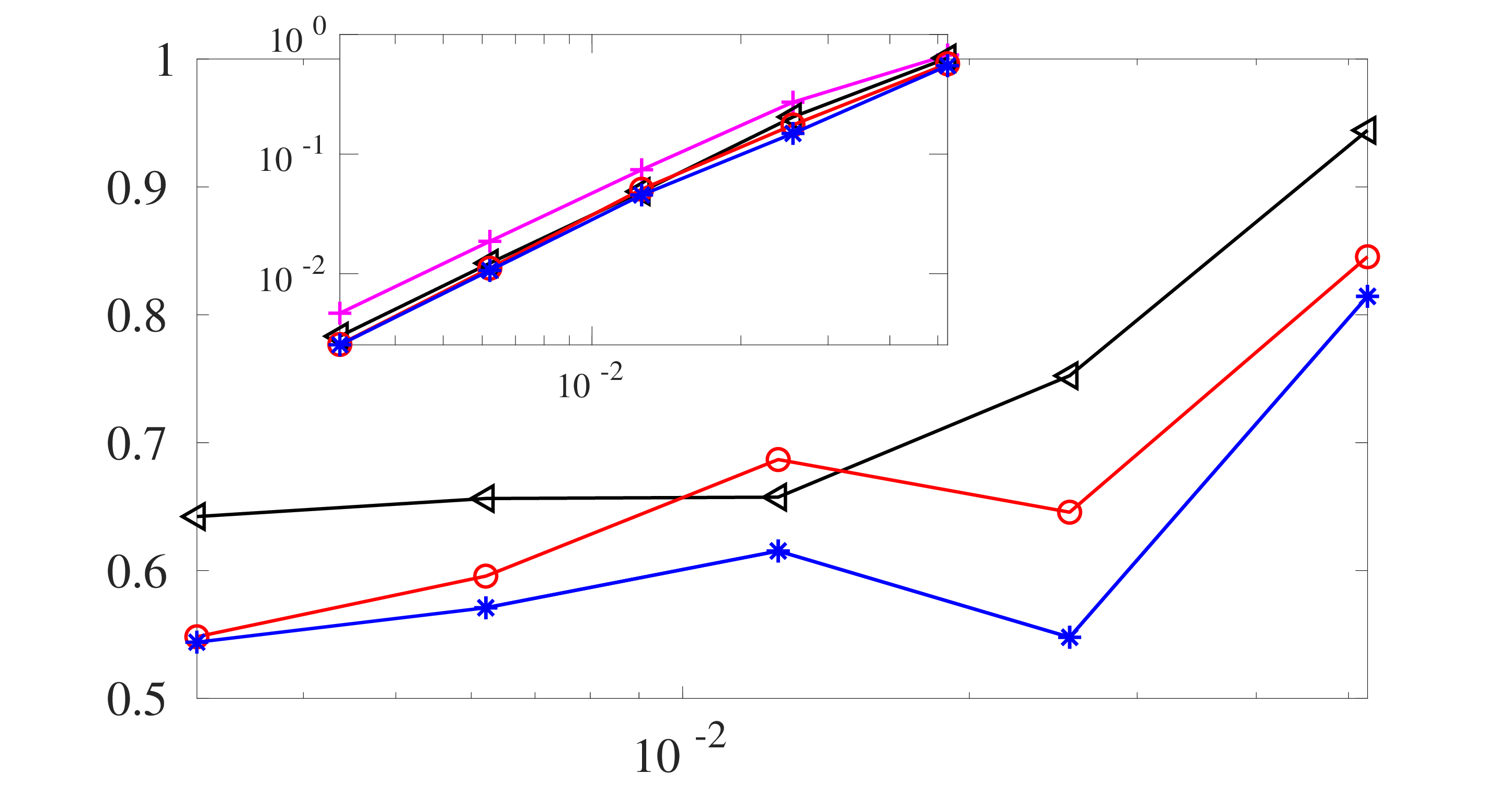}}
    \put(9,0){\includegraphics[width=0.49\linewidth]{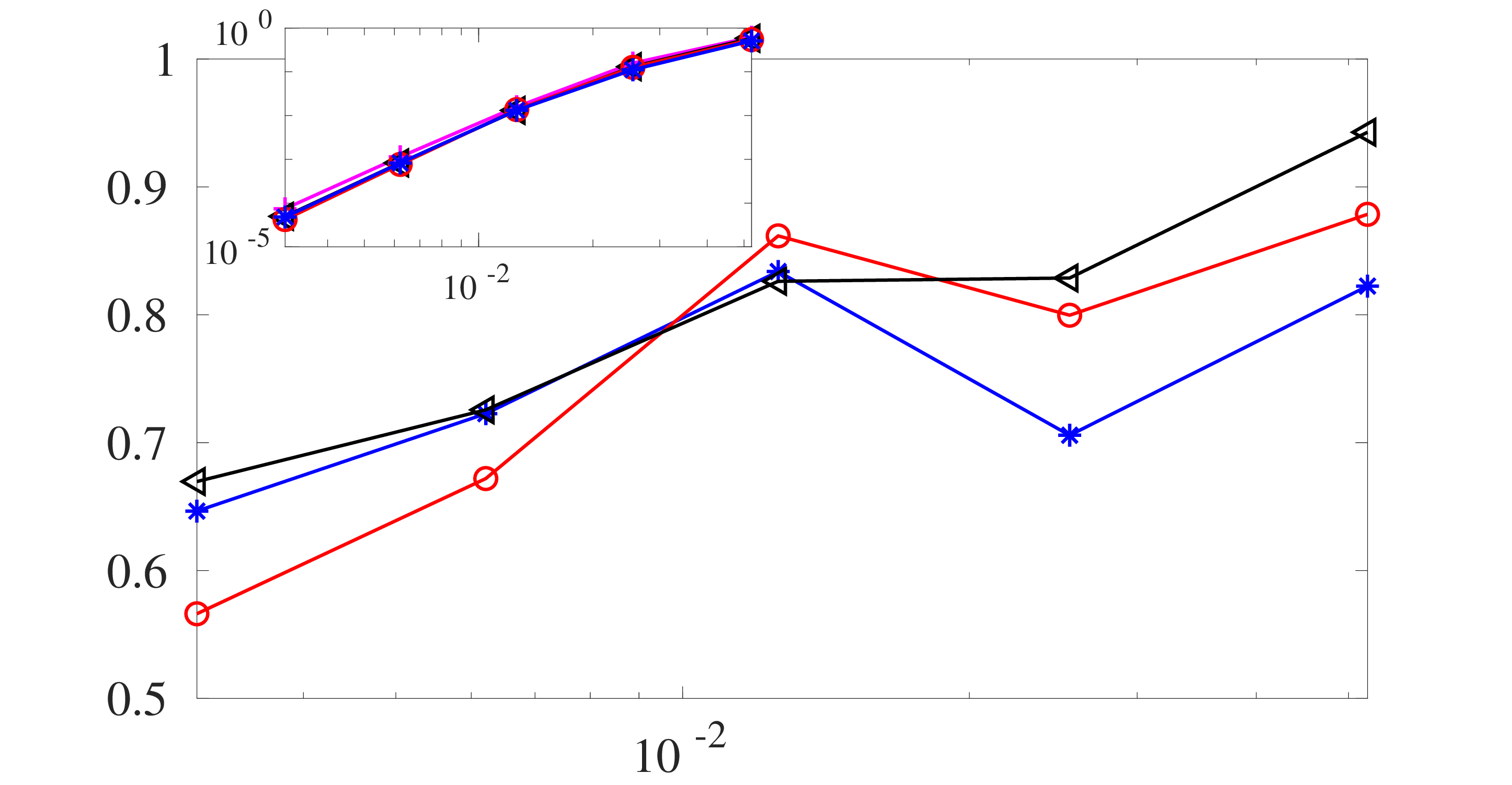}}
    \put(-0.1,2.25){\rotatebox{90}{$\mathcal{R}$}}
   \put(8.9,2.25){\rotatebox{90}{$\mathcal{R}$}}
    \put(4.5,-0.25){$s$}
    \put(13.5,-0.25){$s$}
\end{picture}
    \end{center}
    \caption{Ratio $\mathcal{R}$ of errors arising from explicit and compact approximations to the Laplacian operator acting on function \eqref{eq::test_function}. Left Panel, $2^{nd}$ order approximations: black $\triangleleft$, $2Q_i-1=5$ (scheme b); red $\circ$ $2Q_i-1=9$ (scheme c); blue $*$, $2Q_i-1=13$ (scheme d). Right Panel, $4^{th}$ order approximations; black $\triangleleft$, $2Q_i-1=9$ (scheme f); red $\circ$ $2Q_i-1=13$ (scheme g); blue $*$, $2Q_i-1=17$ (scheme h). Insets show convergence of $L_2-$ norm with magenta line corresponding to explicit LABFM.}
    \label{fig:lap_conv}
    \end{figure}

To compare the convergence of explicit LABFM to compact LABFM we discretise a doubly-periodic square domain $(x,y) \in [0,1]\times [0,1]$ and follow \cite{king_2022} in defining a test function
\begin{equation}
    \label{eq::test_function}
    \phi(x,y)=\sin(2\pi y)\frac{4}{\pi}\sum_{k=1}^8\frac{\sin\left(2(2k-1)\pi(x-1/4)\right)}{2k-1},
\end{equation}
which has sharp gradients at $x=0.25, 0.75$ and is the first eight terms of the Fourier series of a top-hat function. Approximating this function, which has a maximum wavenumber of $|\boldsymbol{k}|=\sqrt{904}\pi$, across a range of resolutions allows us to see improvements in accuracy provided by the compact formulation across a range of wavenumbers. Note that at the coarsest resolution presented in Figs.\ref{fig:grad_conv}\ and \ref{fig:lap_conv} the shortest wavelength in \eqref{eq::test_function} cannot be resolved.

Fig.\ref{fig:grad_conv} shows how $\mathcal{R}$ varies with the nodal spacing $s$ for gradient approximations. At all but the coarsest resolution the compact operators provide a clear improvement. For schemes (c) \& (d) presented in \S\ref{sec:Resolving Power} there is an $80\%$ reduction in error compared to the $2^{nd}$ order explicit scheme. Even at the coarsest resolution which can resolve all wavenumbers of \eqref{eq::test_function} these schemes have approximately a $50\%$ reduction in error. Fourth order schemes provide similar improvements, with over a $60\%$ error reduction on the finest resolution for all schemes.

We now consider Laplacian approximations to \eqref{eq::test_function}. Fig. \ref{fig:lap_conv} compares the compact and explicit schemes with second and fourth order convergence properties. Consistent improvement across all resolutions can be seen for both orders of approximations. The second order compact schemes (c) \& (d) have a $30-50\%$ error reduction for all resolutions. The $4^{th}$ order compact schemes (f)--(h) achieve similar improvements at finer resolutions, though at coarser resolutions the improvements are less clear. Note that the greatest improvement in the resolving efficiency of Laplacian operators was observed along the line $k_y=k_x$, along which lies the smallest non-zero wavenumber in function \eqref{eq::test_function}. It is therefore reasonable to expect the improvement in approximations to Laplacian operators to be better for functions which are dominated by wavenumbers close to $k_y=k_x$ unlike \eqref{eq::test_function}.

\subsection{Stability Tests}
\label{subsec:Stability_Tests}
In time-dependent PDE systems, stability of numerical approximations of operators is also important. For explicit approximations, the stability of the operator can be determined by solving the eigenproblem
\begin{equation}
    \label{eq:Explict_eval_problem}
    \boldsymbol{A}\boldsymbol{v}=\lambda \boldsymbol{v}
\end{equation}
where $\boldsymbol{A}$ is a matrix containing the weights for a given operator, $\lambda$ is the eigenvalue, and  $\boldsymbol{v}$ is the corresponding eigenvector. Explicit LABFM approximations to gradient operators are known to be unstable, similar to other meshfree methods. To determine the stability of implicit approximations, instead of \eqref{eq:Explict_eval_problem} we must instead solve the generalised eigenvalue problem 
\begin{equation}
    \label{sec:Implicit_eval_problem}
    \boldsymbol{A}\boldsymbol{v}=\lambda \boldsymbol{B}\boldsymbol{v}
\end{equation}
where $\boldsymbol{B}$ is a matrix containing the implicit coefficients (i.e. the $\alpha_{q,i}^d$ terms). We solve this problem using the OpenBLAS library.

\begin{figure}[t]
    \begin{center}
    \setlength{\unitlength}{1cm}
    \begin{picture}(18,7)(0,0)
    \put(0,3.5){\includegraphics[width=0.32\linewidth]{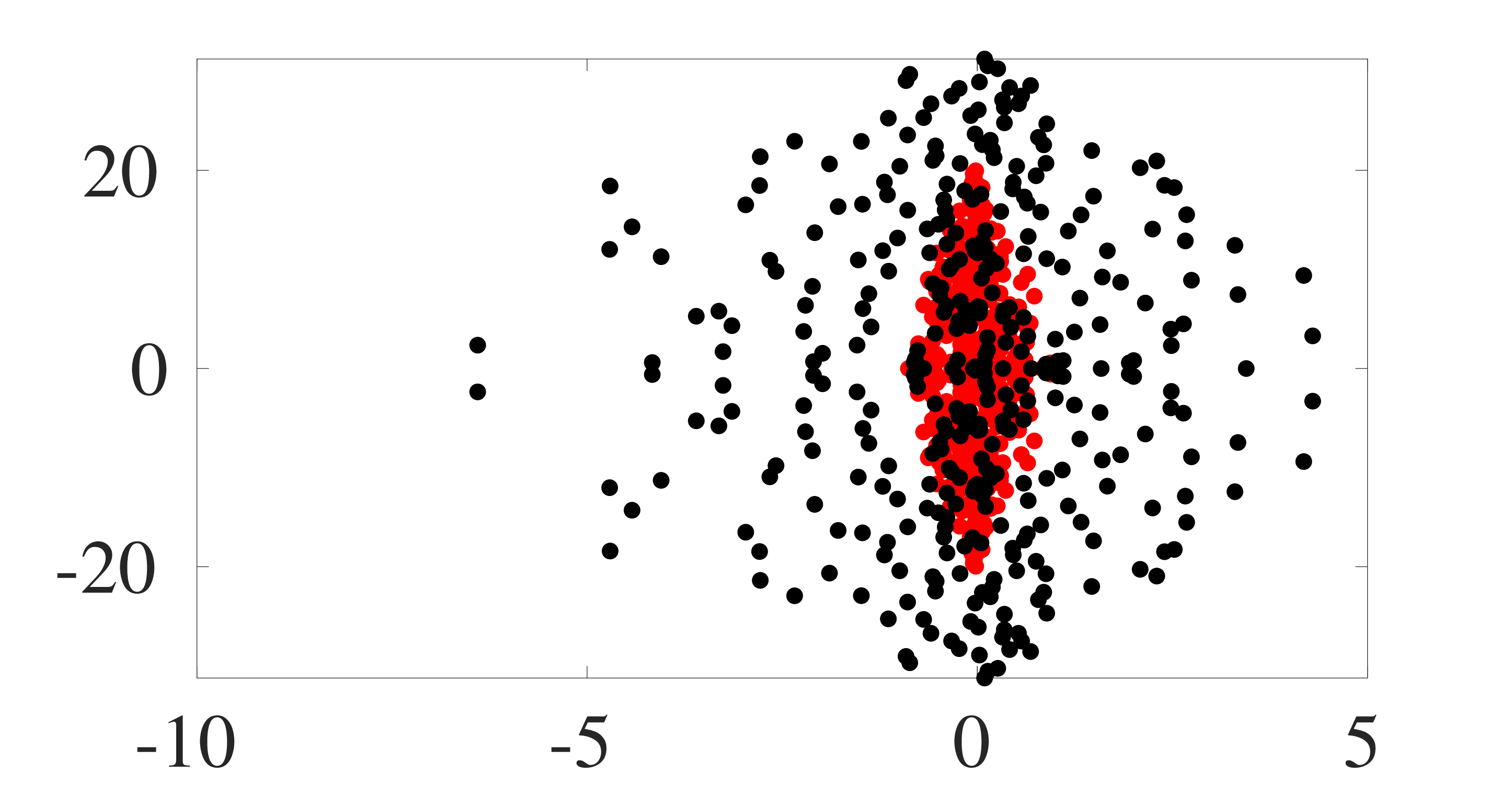}}
    \put(6,3.5){\includegraphics[width=0.32\linewidth]{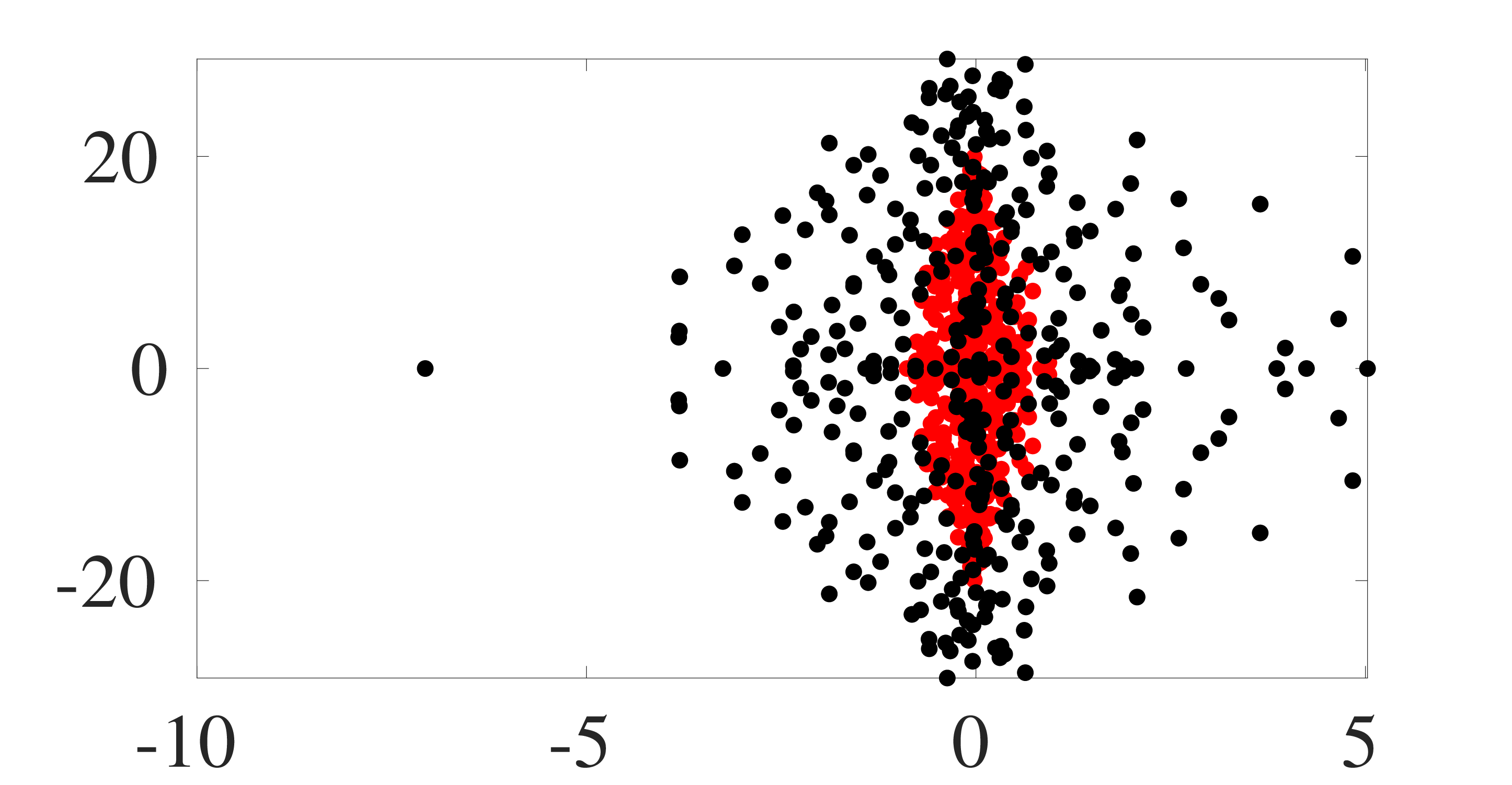}}
    \put(12,3.5){\includegraphics[width=0.32\linewidth]{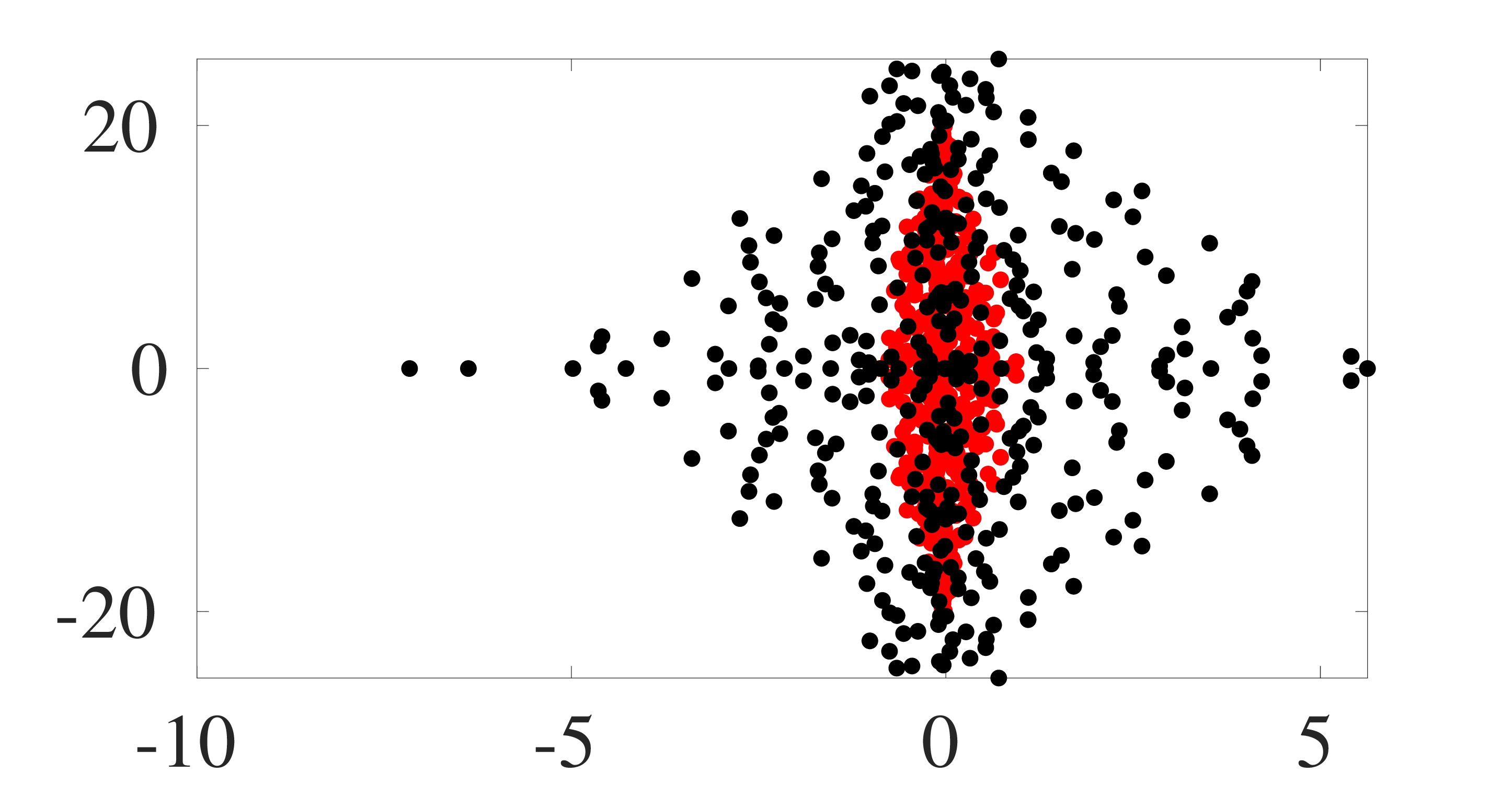}}
    \put(2.5,3.1){$\Re\{\lambda\}$}
    \put(8.5,3.1){$\Re\{\lambda\}$}
    \put(14.5,3.1){$\Re\{\lambda\}$}
     \put(-0.25,4.65){\rotatebox{90}{$\Im\{\lambda\}$}}
     \put(5.75,4.65){\rotatebox{90}{$\Im\{\lambda\}$}}
     \put(11.75,4.65){\rotatebox{90}{$\Im\{\lambda\}$}}
     \put(0,0){\includegraphics[width=0.32\linewidth]{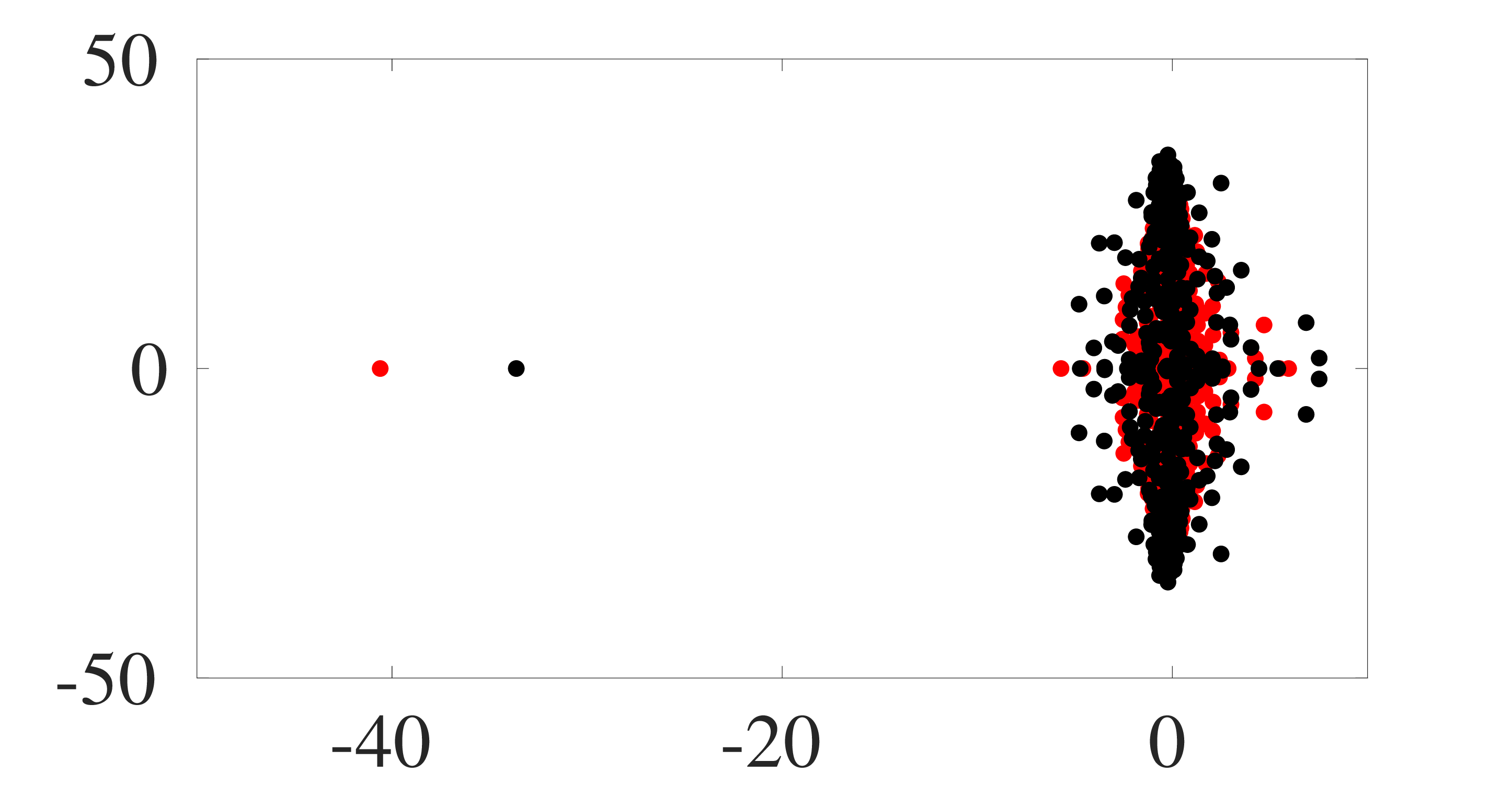}}
    \put(6,0){\includegraphics[width=0.32\linewidth]{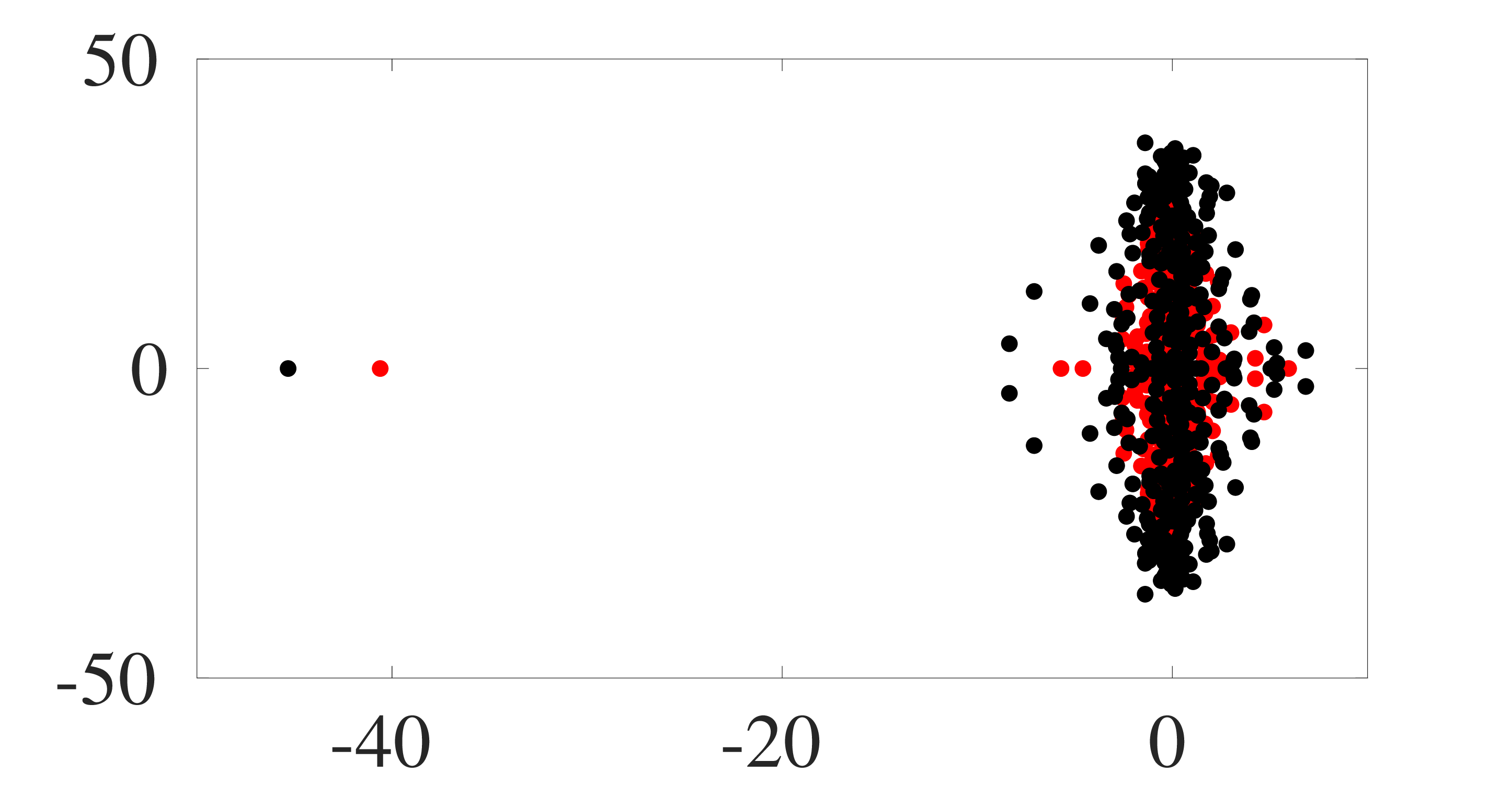}}
    \put(12,0){\includegraphics[width=0.32\linewidth]{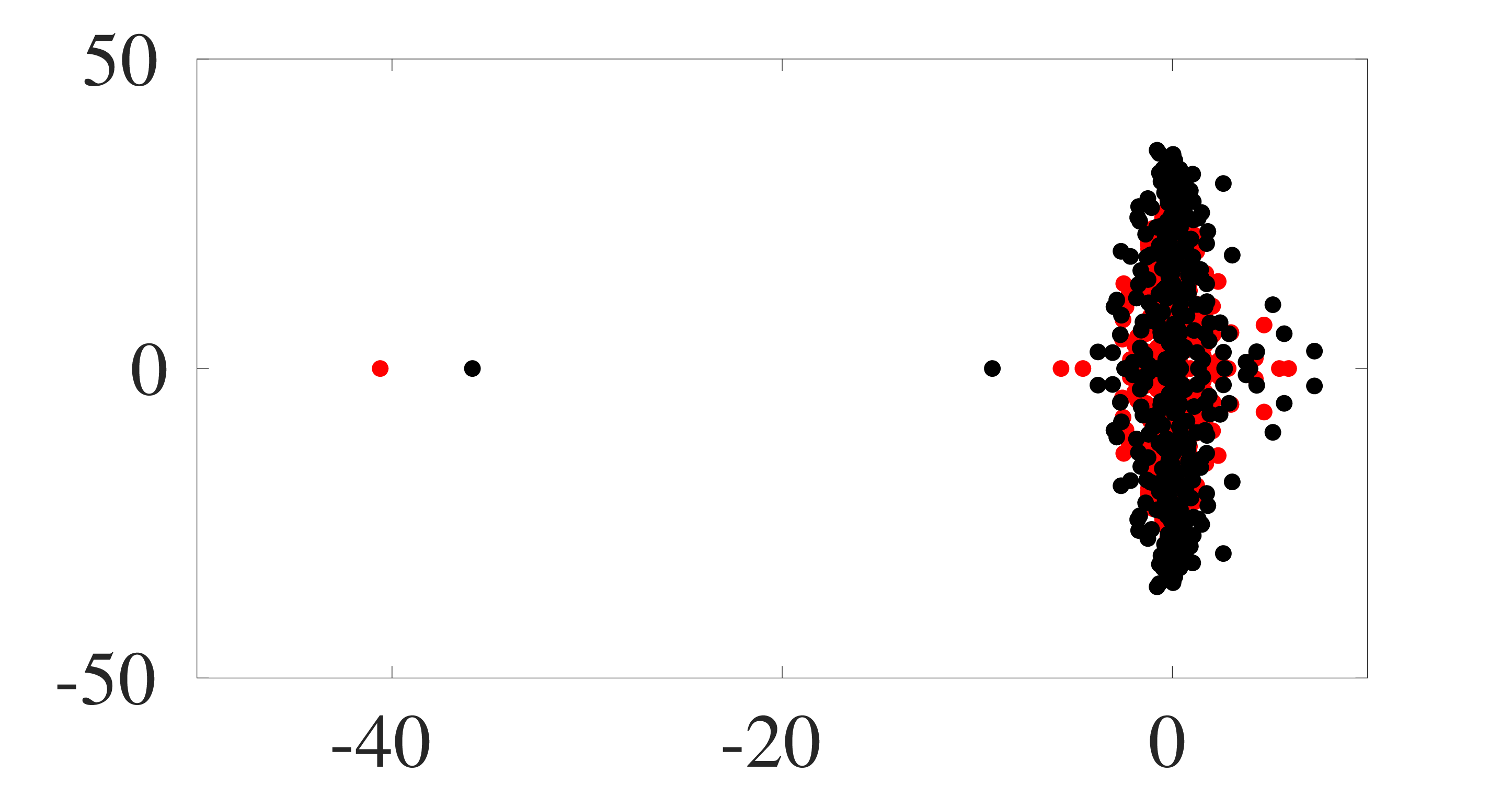}}
    \put(2.5,-0.5){$\Re\{\lambda\}$}
    \put(8.5,-0.5){$\Re\{\lambda\}$}
    \put(14.5,-0.5){$\Re\{\lambda\}$}
     \put(-0.25,1.15){\rotatebox{90}{$\Im\{\lambda\}$}}
     \put(5.75,1.15){\rotatebox{90}{$\Im\{\lambda\}$}}
     \put(11.75,1.15){\rotatebox{90}{$\Im\{\lambda\}$}}
    \end{picture}
    \end{center}
    \caption{Stability of a gradient operator for $2^{nd}$ order (top row) and $4^{th}$ order (bottom row) approximations on a typical discretisation with $441$ points. Red dots correspond to explicit LABFM, black dots to those produced using compact LABFM. Top row: left panel $Q_i=3$ (scheme b); centre panel $Q_i=5$ (scheme c); right panel: $Q_i=7$ (scheme d). Bottom row: left panel $Q_i=5$ (scheme f); centre panel $Q_i=7$ (scheme g); right panel: $Q_i=9$ (scheme h).}
    \label{fig:grad_stab}
\end{figure}
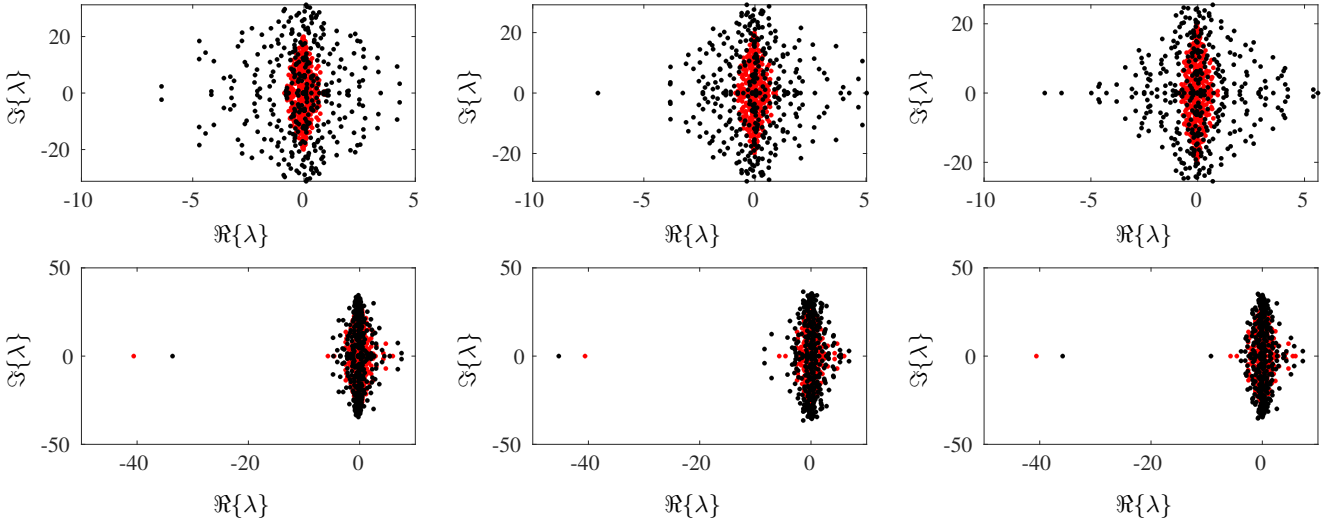

Fig.\ref{fig:grad_stab} compares the stability of gradient operators for $2^{nd}$ and $4^{th}$ order compact schemes to their explicit counterparts. The $2^{nd}$ order compact schemes have a relatively consistent maximal growth rate of $\Re\{\lambda\}\approx 5$. These schemes are more unstable than the explicit scheme, which has $\Re\{\lambda\}\approx 1$. The difference in stability for $4^{th}$ order schemes is smaller, however the explicit scheme remains slightly less unstable. In any case, purely advective problems are unstable for all approximations. If a more complex PDE system is advection dominated with insufficient diffusion then it will require numerical filtering to maintain stability for all schemes presented herein. 

\begin{figure}[t]
    \begin{center}
    \setlength{\unitlength}{1cm}
    \begin{picture}(18,7)(0,0)
    \put(0,3.5){\includegraphics[width=0.32\linewidth]{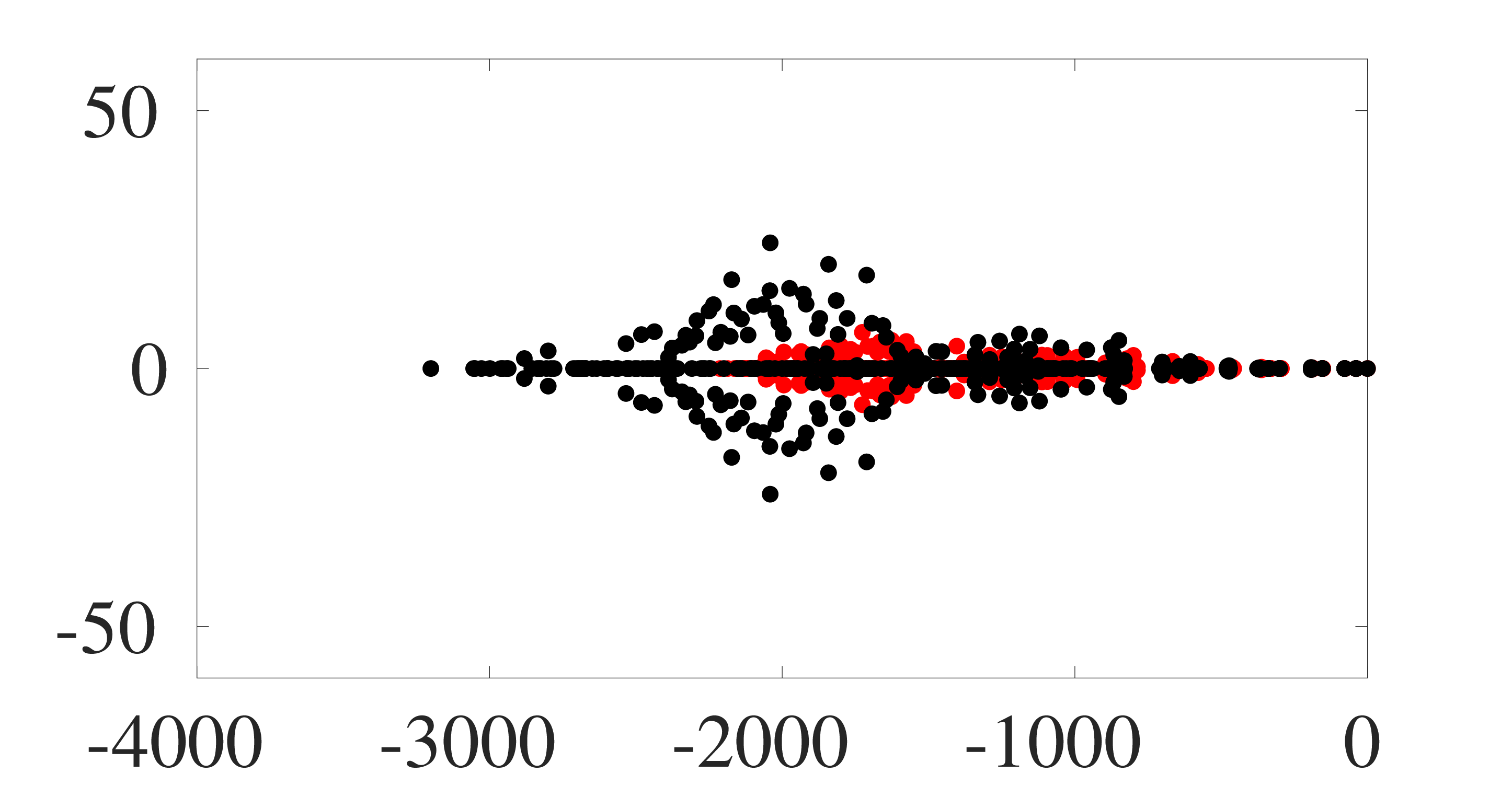}}
    \put(6,3.5){\includegraphics[width=0.32\linewidth]{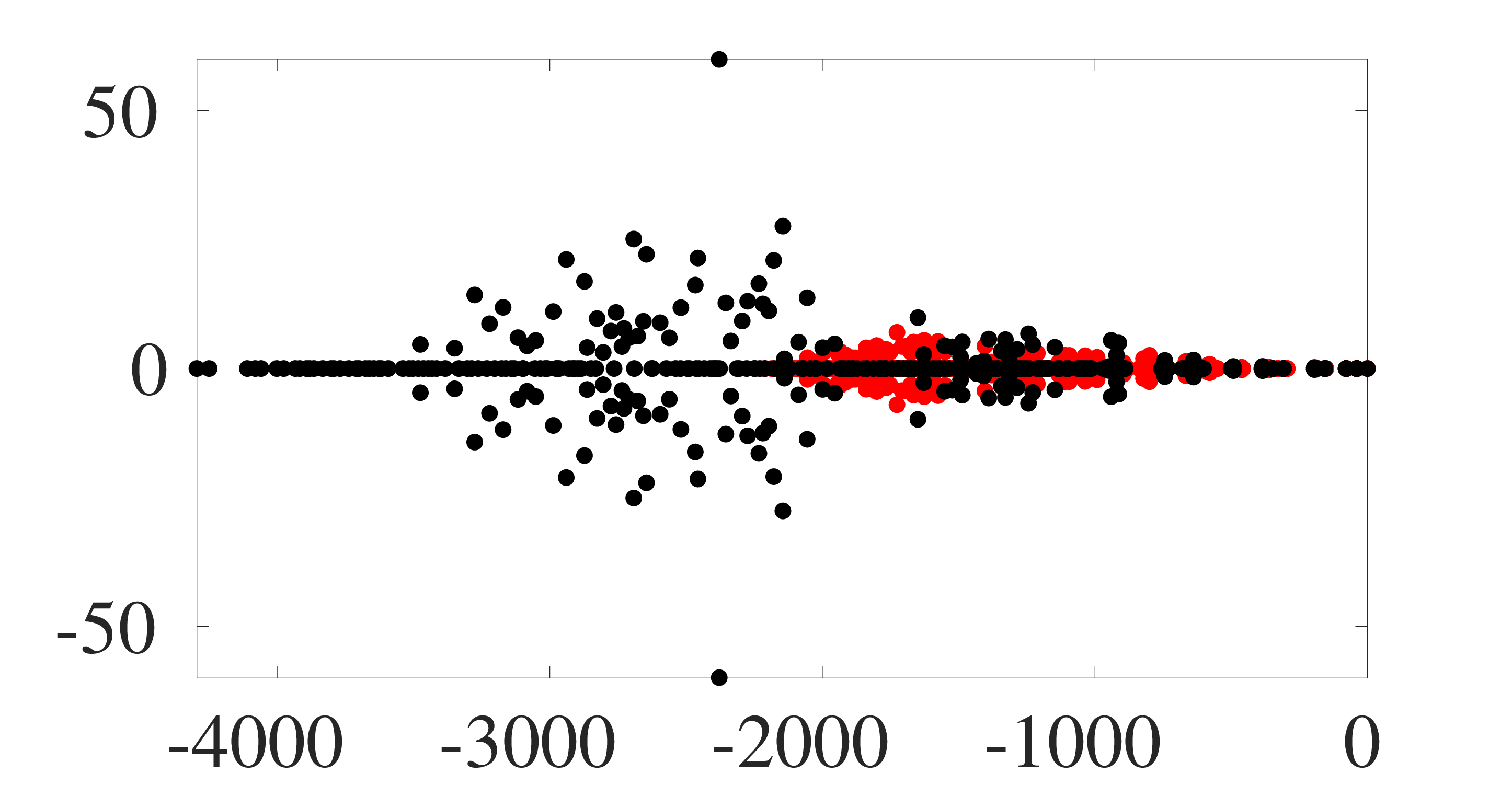}}
    \put(12,3.5){\includegraphics[width=0.32\linewidth]{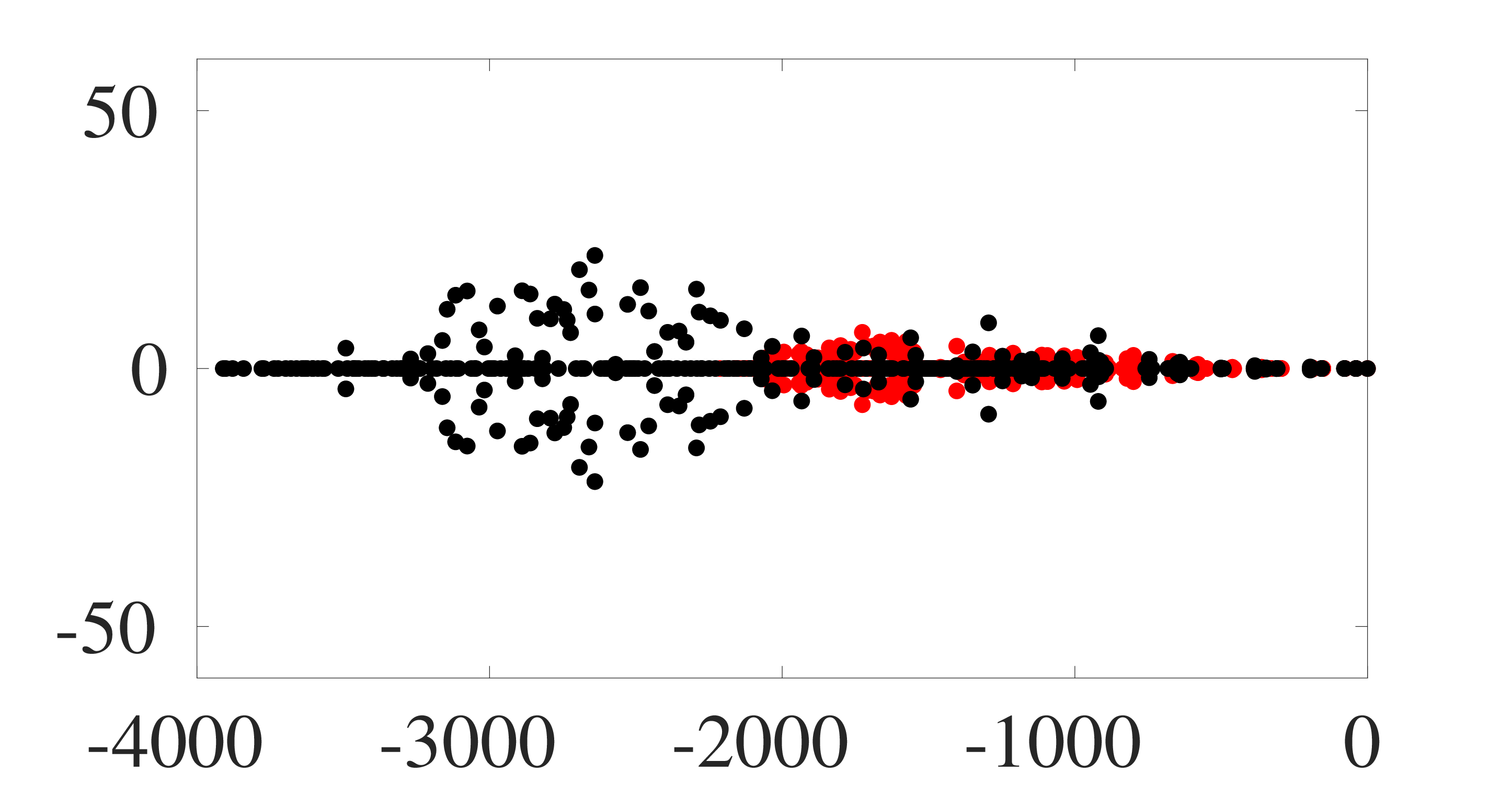}}
    \put(2.5,3.1){$\Re\{\lambda\}$}
    \put(8.5,3.1){$\Re\{\lambda\}$}
    \put(14.5,3.1){$\Re\{\lambda\}$}
     \put(-0.25,4.65){\rotatebox{90}{$\Im\{\lambda\}$}}
     \put(5.75,4.65){\rotatebox{90}{$\Im\{\lambda\}$}}
     \put(11.75,4.65){\rotatebox{90}{$\Im\{\lambda\}$}}
     \put(0,0){\includegraphics[width=0.32\linewidth]{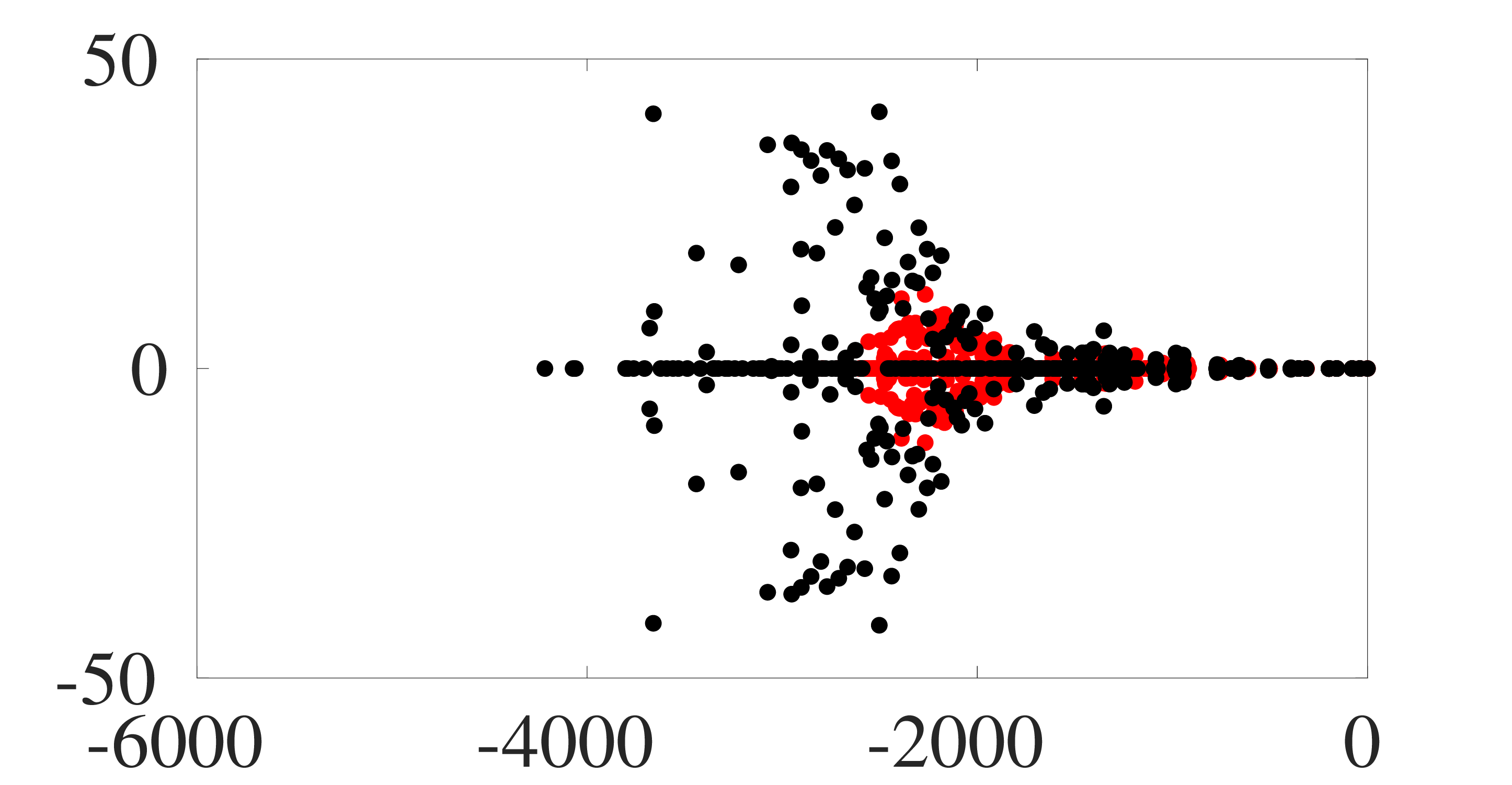}}
    \put(6,0){\includegraphics[width=0.32\linewidth]{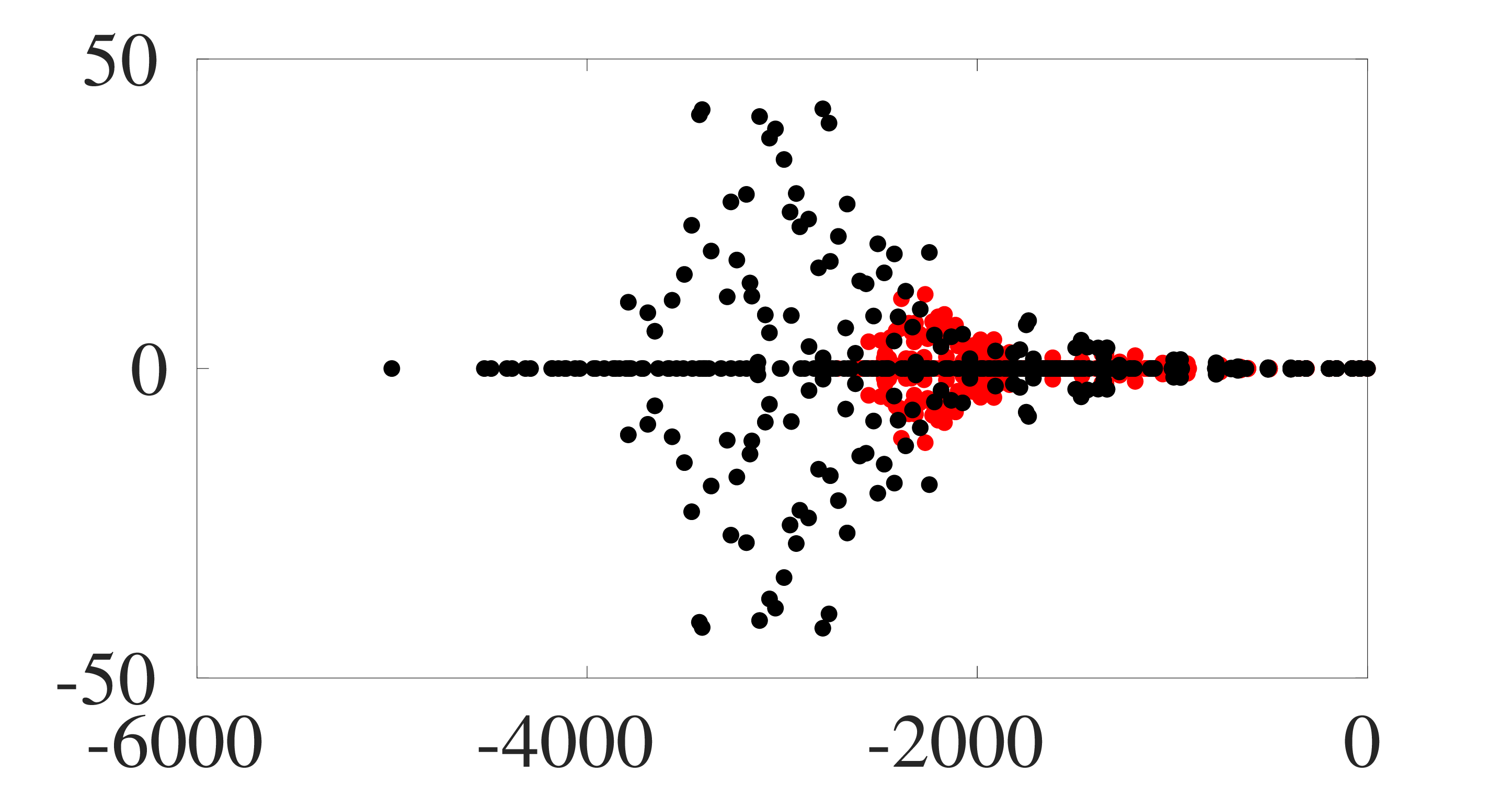}}
    \put(12,0){\includegraphics[width=0.32\linewidth]{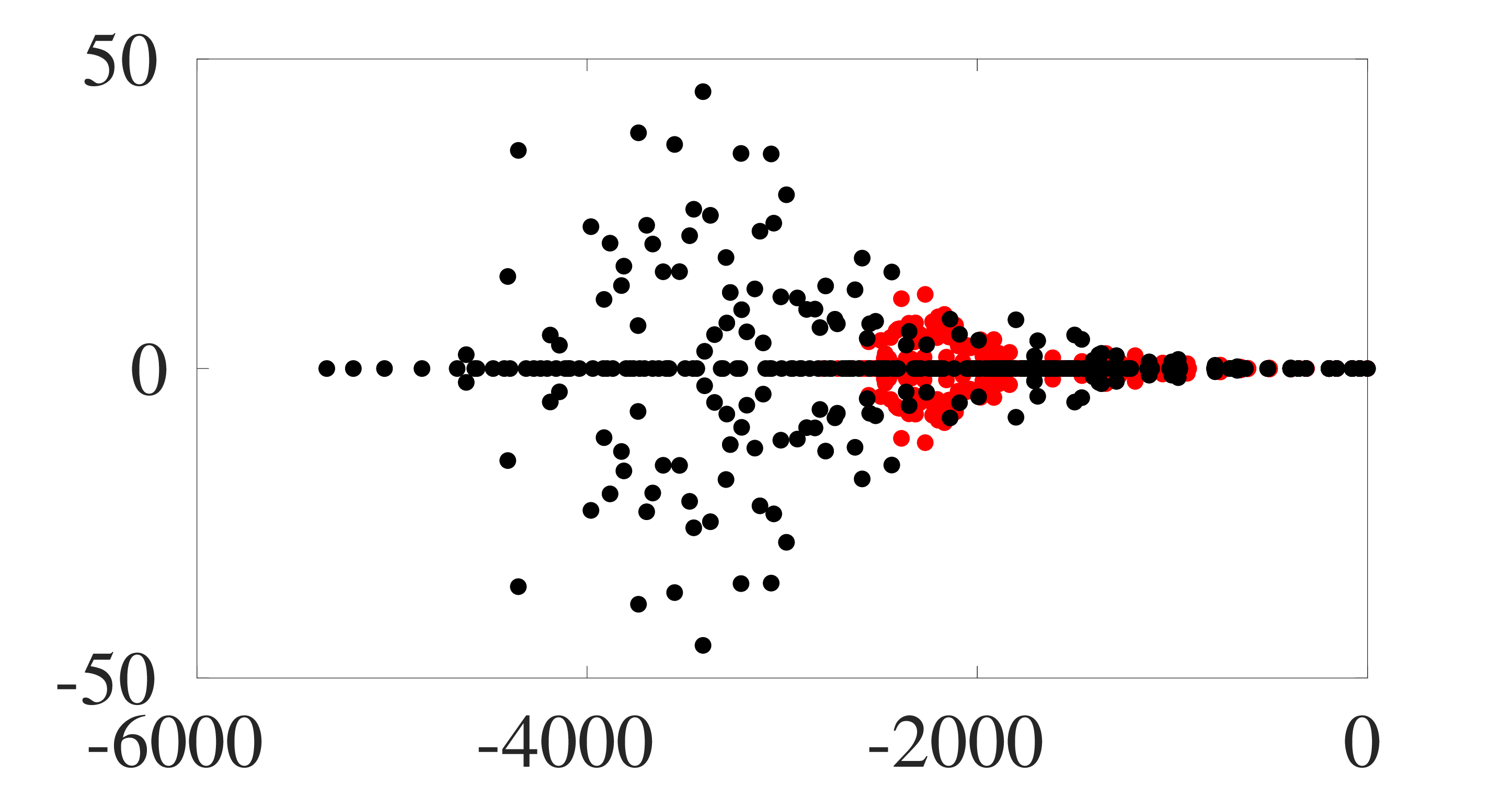}}
    \put(2.5,-0.5){$\Re\{\lambda\}$}
    \put(8.5,-0.5){$\Re\{\lambda\}$}
    \put(14.5,-0.5){$\Re\{\lambda\}$}
     \put(-0.25,1.15){\rotatebox{90}{$\Im\{\lambda\}$}}
     \put(5.75,1.15){\rotatebox{90}{$\Im\{\lambda\}$}}
     \put(11.75,1.15){\rotatebox{90}{$\Im\{\lambda\}$}}
    \end{picture}
    \end{center}
    \caption{Stability of the Laplacian operator for $2^{nd}$ order (top row) and $4^{th}$ order (bottom row) approximations on a typical discretisation with $441$ points. Red dots correspond to explicit LABFM, black dots to those produced using compact LABFM.  Top row: left panel $2Q_i-1=5$ (scheme b); centre panel $2Q_i-1=9$ (scheme c); right panel: $2Q_i-1=13$ (scheme d). Bottom row: left panel $2Q_i-1=9$ (scheme f); centre panel $2Q_i-13=13$ (scheme g); right panel: $2Q_i-1=17$ (scheme h).}
    \label{fig:lap_stab}
\end{figure}
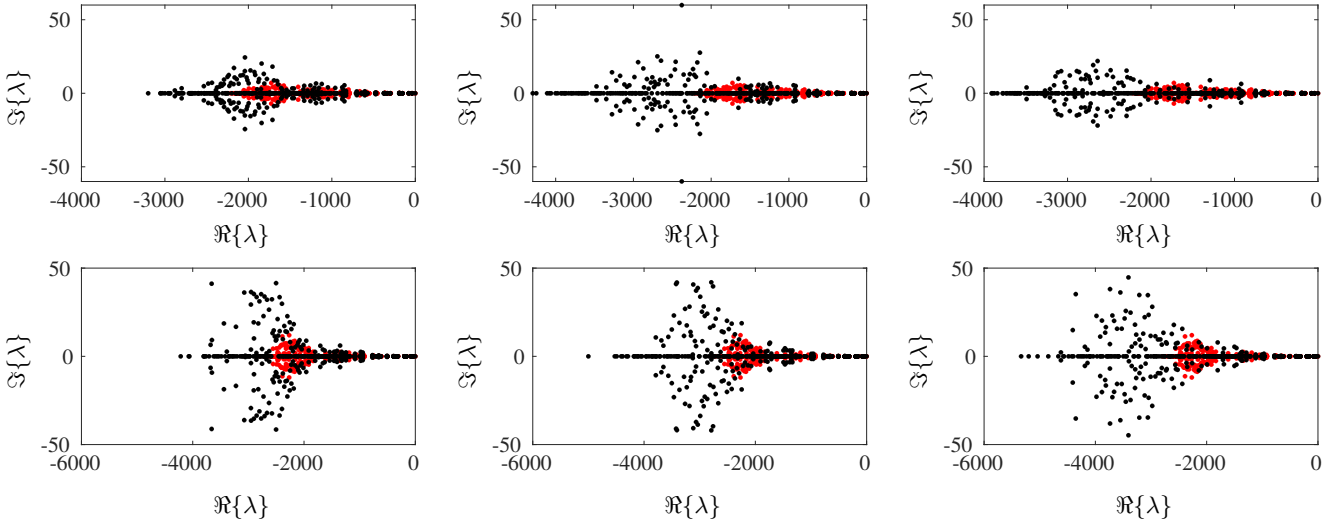

We now consider the stability of Laplacian operators. Fig.\ref{fig:lap_stab} compares the stability of compact and explicit operators on a typical discretisation. All operators are temporally stable, though compact schemes typically have a larger max$\{\Im\{\lambda\}\}$. The imaginary parts of these eigenvalues remain the same magnitude as explicit LABFM and given their large, negative real parts (thus rapid temporal decay), they are unlikely to affect the accuracy of approximations.

\section{Applications to PDEs}
\label{sec:PDE_solutions}
Now that we have considered the properties of compact LABFM, we turn to using the new implicit formulation to solve PDE systems. 

\subsection{Viscous Burgers Equation}
\label{subsec:Burgers_Equation}
To demonstrate the benefit of using compact approximations to resolve fine spatial structures in a time-dependent PDE system we consider viscous Burgers equation, a prototypical equation for the development, and subsequent dissipation, of shock waves. We solve the system
\begin{gather}
    \frac{\partial \boldsymbol{u}}{\partial t}+\boldsymbol{u}\cdot (\nabla \boldsymbol{u})=\frac{1}{Re} \nabla^2 \boldsymbol{u}
\end{gather}
in a periodic domain $(x,y)\in [0,1]\times [0,1]$, with the initial condition
\begin{gather}
    u(x,y,0)=\sin (2\pi x),\qquad v(x,y,0)\equiv 0
\end{gather}
which has solution \cite{cole}
\begin{equation}
    \label{BE_soln}
    u(x,t)=\frac{4\pi}{Re}\frac{\sum_{n=1}^\infty nA_n\sin\left(2n\pi x\right)e^{-\frac{n^2\pi^2t}{Re}}}{A_0+\sum_{n=1}^\infty A_n\cos\left(2n\pi x\right)e^{-\frac{n^2\pi^2t}{Re}}}
\end{equation}
where
\begin{equation}
    A_0=e^{-\frac{Re}{4\pi}}I_0\left(\frac{Re}{4\pi}\right),\qquad A_n=2e^{-\frac{Re}{4\pi}}I_n\left(\frac{Re}{4\pi}\right)
\end{equation}
and $I_n$ are the modified Bessel functions of the first kind. The $A_n$ terms decay relatively quickly, and only the first $30$ terms are required to obtain an analytic solution converged to machine precision. The initial sinusoidal velocity profile develops towards a saw-tooth profile due to non-linear advection, before decaying due to viscous dissipation. To ensure temporal stability as we evolve the system we set the time-step condition 
\begin{equation}
\label{eq::time_step}
    \delta t=\min \left( \frac{0.1\ \mathrm{min}(s_i)}{\mathrm{max}(|u|)} ,\frac{0.05\ \mathrm{min}(s^2_i)}{Re}\right)
\end{equation}
and use the classical fourth order Runge-Kutta scheme (RK4). This relatively small $\delta t$ is chosen to ensure that errors due to time-stepping remain small compared to errors arising from the approximations of spatial derivatives, and we note our simulations remained stable for significantly larger choices of $\delta t$ for both explicit and compact simulations. 

\begin{figure}
    \begin{center}
    \setlength{\unitlength}{1cm}
    \begin{picture}(18,5)(0,0)
    \put(0,0){\includegraphics[width=0.49\linewidth]{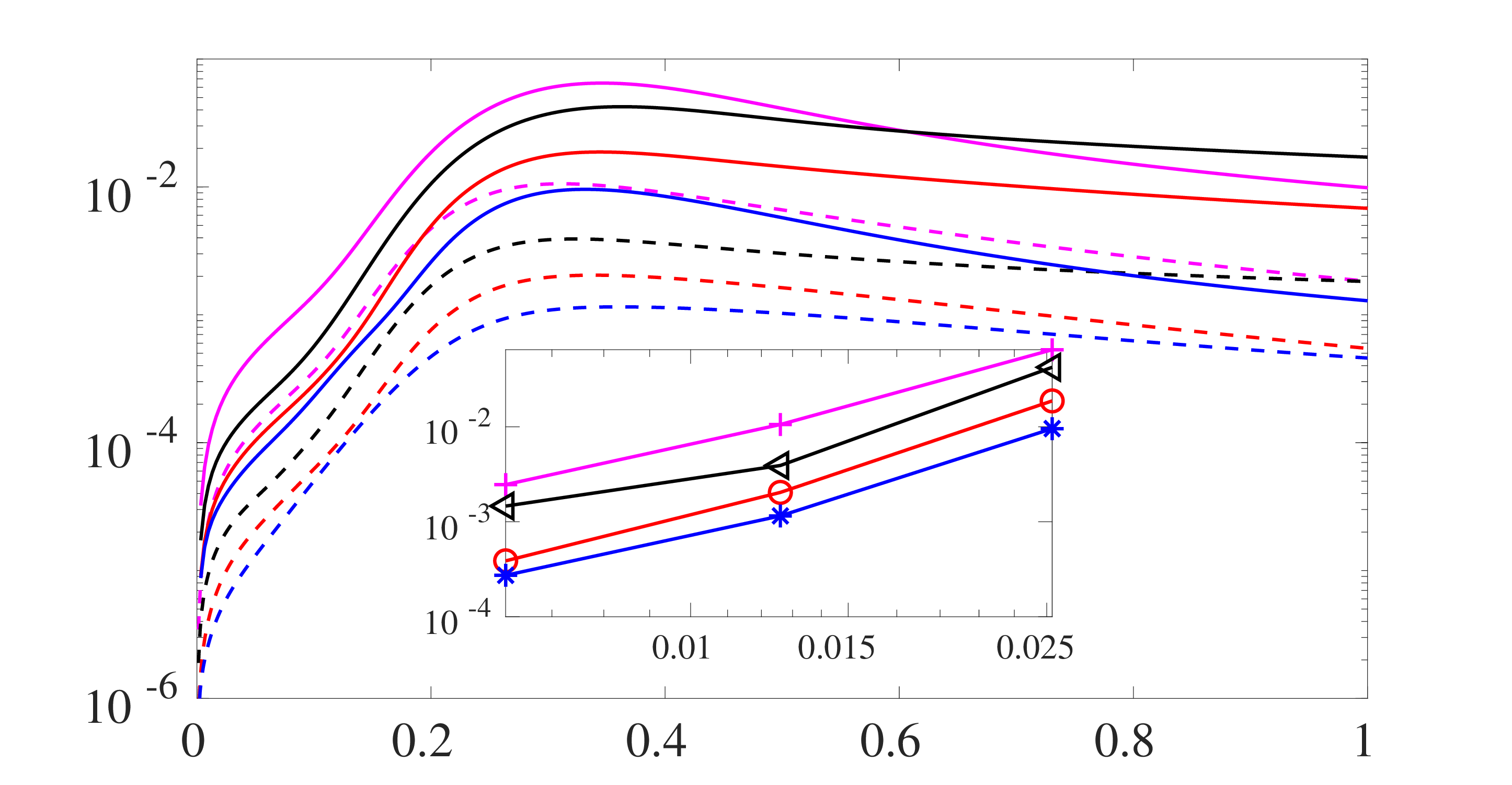}}
    \put(9,0){\includegraphics[width=0.49\linewidth]{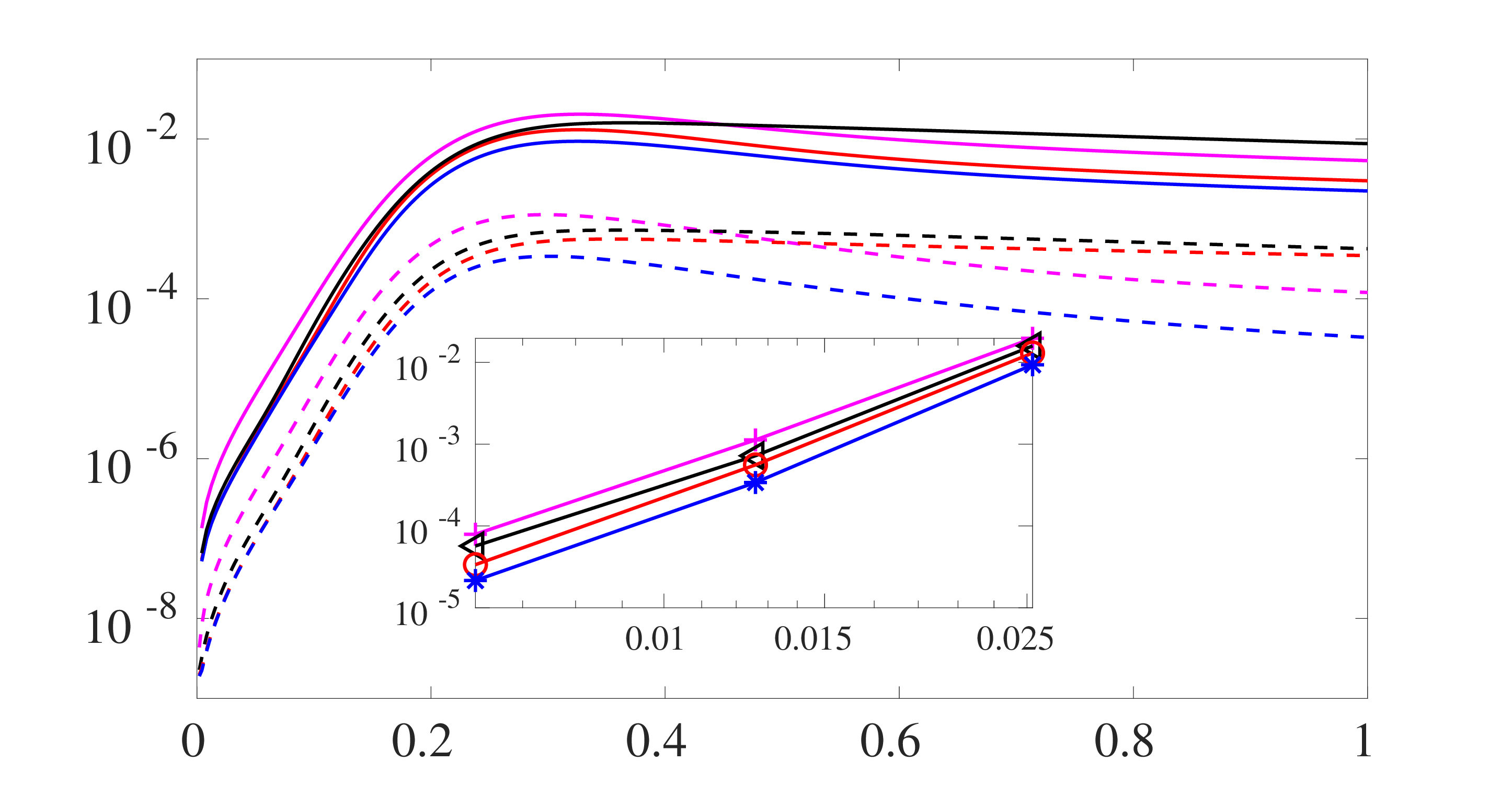}}
    \put(13.5,-0.25){$t$}
    \put(-0.1,1.75){\rotatebox{90}{$L_2-$norm}}
    \put(8.9,1.75){\rotatebox{90}{$L_2-$norm}}
   
    \put(4.5,-0.25){$t$}
    \put(4.25,0.65){{{$s$}}}
    \put(2.25,0.6){\rotatebox{90}{{\tiny{max($L_2-$norm)}}}}

    \put(13.25,0.65){{{$s$}}}
    \put(11.05,0.6){\rotatebox{90}{{\tiny{max($L_2-$norm)}}}}
    
    \end{picture}
    \end{center}
    \caption{Time evolution of $L_2-$norm of solution to Burgers equation with $Re=100$. Left Panel, $2^{nd}$ order schemes: magenta lines explicit scheme (a); black lines scheme (b); red lines scheme (c); blue lines scheme (d). Right Panel $4^{th}$ order schemes: magenta lines scheme (e); black lines scheme (f); red lines scheme (g); blue lines scheme (h). Schemes as denoted in Table \ref{tab:stencils}. Solid lines in main figure denote resolution $s=1/40$, dashed lines $s=1/80$. Insets show convergence of maximum value of $L_2-$norm for increasing resolution. Note that finest resolution of inset is not shown in main figure.}
    \label{fig:Re100}
\end{figure}
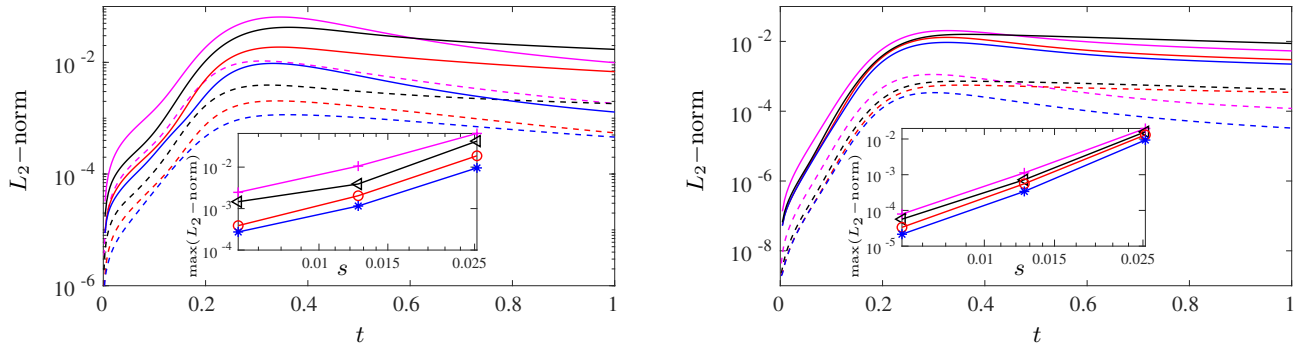

We simulate the case $Re=100$. Fig. \ref{fig:Re100} shows the evolutions of the $\rnumber$ with time for all compact and explicit schemes. In this case, the early time dynamics (up to $t\approx 0.4$) are dominated by large gradients which occur due to the formation of the shock wave. Beyond this time, dissipation is the most important physical feature of the system, however there are other important numerical factors such as global conservation of (for example) linear momentum. In all simulations, compact schemes outperform the explicit schemes in the advection dominated early time regime, with more points in the implicit stencil leading to more accurate simulations. The $2^{nd}$ order compact schemes offer up to an order of magnitude improvement, and the error of $4^{th}$ order schemes are reduced by up to $80\%$. In fact, on the coarsest resolution the $2^{nd}$ order compact scheme with $Q_i=7$ has a smaller error than the explicit $4^{th}$ order scheme. However, for later times compact schemes do not necessarily produce more accurate results. This is because of other factors such as conservation properties of the numerical scheme. We note that the schemes with the largest implicit stencils retain the smallest $\rnumber$ for all times, with the magnitude of improvement in approximation remaining similar. 

The improvements in resolving power demonstrated \S\ref{sec:Resolving Power} are consistently observed in the advection dominated dynamics. To illustrate this improvement, the insets in Fig.\ref{fig:Re100} show the convergence of the temporal maximum of the $\rnumber$ (which occurs during this phase) for each scheme. For both $2^{nd}$ and $4^{th}$ order approximations the larger the implicit stencil the smaller the maximum error at all resolutions.


\subsection{Poisson's Equation}
\label{subsec:Poisson_Equation}
\begin{figure}
\begin{center}
    \setlength{\unitlength}{1cm}
    \begin{picture}(18,5)(0,0)
     \put(2,0){\includegraphics[width=0.4\linewidth]{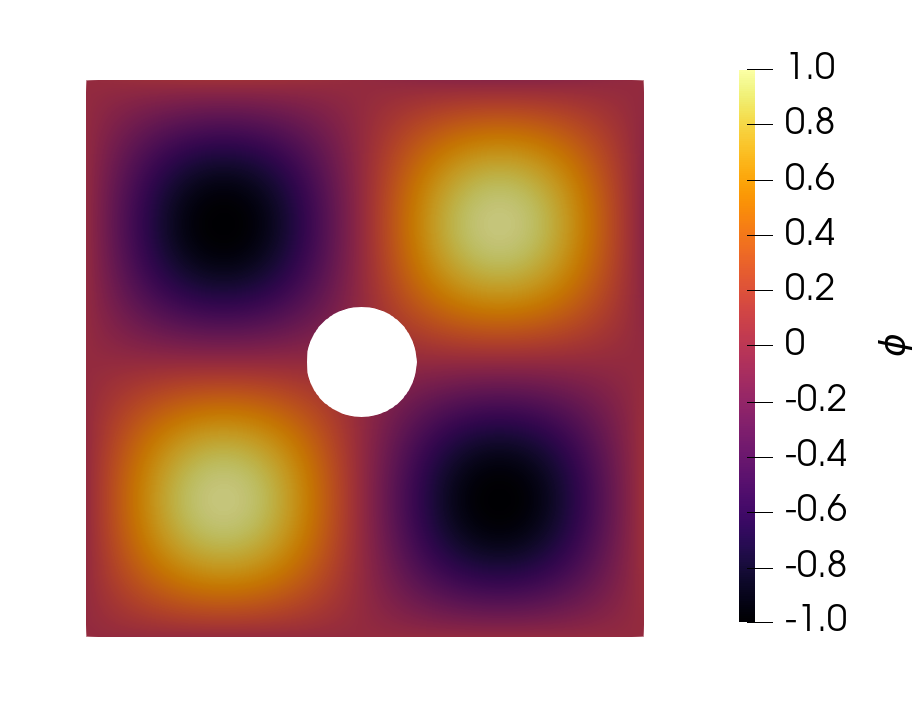}}
     \put(10,0.3){\includegraphics[width=0.3\linewidth]{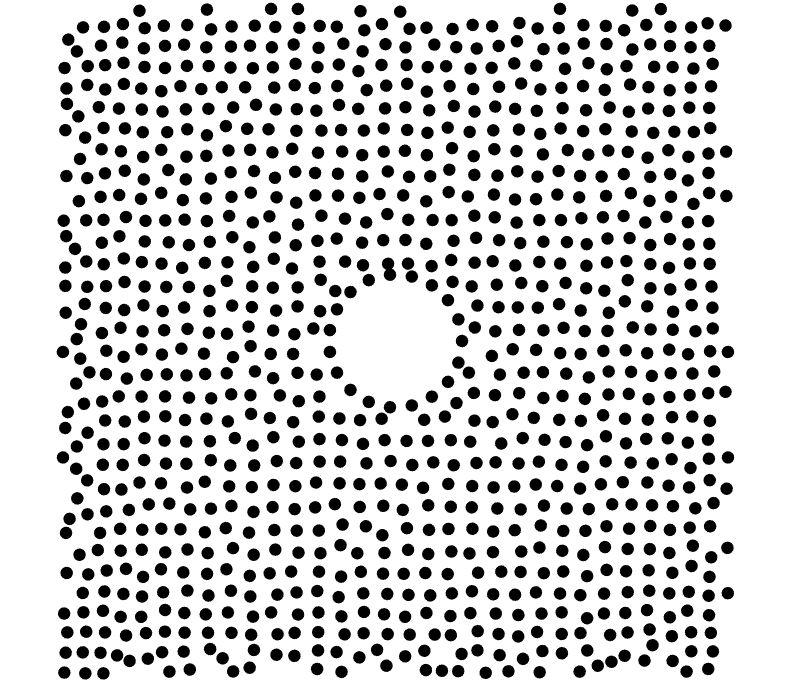}}
    \end{picture}
    \end{center}
    \caption{Left: the domain \eqref{eq:Poisson_domain} and solution \eqref{eq:poisson_solution} for the Poisson problem given by \eqref{eq:Poisson_function}. Right: example discretisation.}
    \label{fig:Poisson_Complex_geometry}
\end{figure}

Another prototypical PDE is Poisson's equation. We wish to find solutions to the equation 
\begin{equation}
    \label{eq::Poisson_equation}
    \nabla^2 \phi=f(x,y)
\end{equation}
in a domain $\Omega\subseteq \mathbb{R}^2$. Analytic solutions to this equation can be constructed trivially even in non-trivial geometries, enabling the benefits of meshless methods to be quantified. To this end, we consider the domain 
\begin{equation}
\label{eq:Poisson_domain}
    \Omega\coloneqq \{(x,y)\in[0,1]\times [0,1]| \sqrt{(x-0.5)^2+(y-0.5)^2}\geq 0.1\},
\end{equation}
an example discretisation of which can be seen in the right panel of Fig.\ref{fig:Poisson_Complex_geometry}.
We set 
\begin{equation}
\label{eq:Poisson_function}
f(x,y)=-8\pi^2 \sin(2\pi x)\sin(2\pi y)
\end{equation}
for which the analytic solution is
\begin{equation}
\label{eq:poisson_solution}
\phi=\sin(2\pi x)\sin(2\pi y)    
\end{equation}
which is shown in the left panel of Fig.\ref{fig:Poisson_Complex_geometry}. We solve \eqref{eq::Poisson_equation} subject to Dirichlet boundary conditions on the inner boundary $\Gamma_d$ where $\sqrt{x^2+y^2}=0.1$, and periodic boundary conditions on the outer boundary. For both explicit and implicit LABFM a global linear system must be solved to obtain the numerical approximation to $\phi$. Once we have discretised our domain with $N$ particles, the compact LABFM approximation of \eqref{eq::Poisson_equation} at each non-boundary node $i$ can be written as
\begin{equation}
    \sum_{j\in \mathcal{N}_i}\phi_{ji}w_{ji}^L=\sum_{q\in \mathcal{M}_i}\alpha_{q,i}^Lf(x_q,y_q).
\end{equation}
For the whole domain this can be expressed in terms of vectors and matrices
\begin{equation}
\label{eq:poisson_mat_eqn}
    \boldsymbol{A}^L \boldsymbol{\phi}=\boldsymbol{\alpha}^L \boldsymbol{f}
\end{equation}
where $\boldsymbol{\phi}$ is the solution vector, $\boldsymbol{f}$ is a vector of values $f$ at each node, and the matrices $\boldsymbol{A}^L, \boldsymbol{\alpha}^L$ take the form
\begin{subequations}
\begin{equation}
    A_{ij}^L=\begin{cases}\qquad w_{ji}^L \ &\mathrm{if} \quad {i\notin \Gamma_d}, \ j\neq i \ \mathrm{and} \ j\in\mathcal{N}_i\\
    -\sum_{j\in\mathcal{N}_i}w_{ji}^L \ &\mathrm{if} \quad i=i, \ {i\notin \Gamma_d}  , \\
   \qquad \delta_{ij} &\mathrm{if} \ i\in\Gamma_d,,\\
   \qquad 0 & \mathrm{otherwise},
    \end{cases}
\end{equation}
\begin{equation}
    \alpha_{iq}=\begin{cases}
        \alpha^L_{q,i} &\mathrm{if} \quad i\notin\Gamma_d, \ q\in{\mathcal{M}_i},\\
        \delta_{iq} &\mathrm{if} \ i\in\Gamma_d,\\
        0 &\mathrm{otherwise}.
        
    \end{cases}
\end{equation}
\end{subequations}
Here $\delta_{ij}$ is the Kronecker delta. For explicit methods the matrix $\boldsymbol{A}^L$ is of the same form as the compact case (though the weights obviously take different values), however, $\alpha_{iq}=\delta_{iq}$ for all $i,j$. The only difference in computational efficiency between compact and explicit formulations (once weights have been calculated) is therefore in the calculation of the right-hand side of \eqref{eq:poisson_mat_eqn} which is small compared to the cost of solving the global linear system. Any improvement from compact operators therefore comes with close to zero increase in computational cost when solving Poisson's equation and other elliptic problems.

\begin{figure}
    \begin{center}
    \setlength{\unitlength}{1cm}
    \begin{picture}(18,5)(0,0)
    \put(0,0){\includegraphics[width=0.49\linewidth]{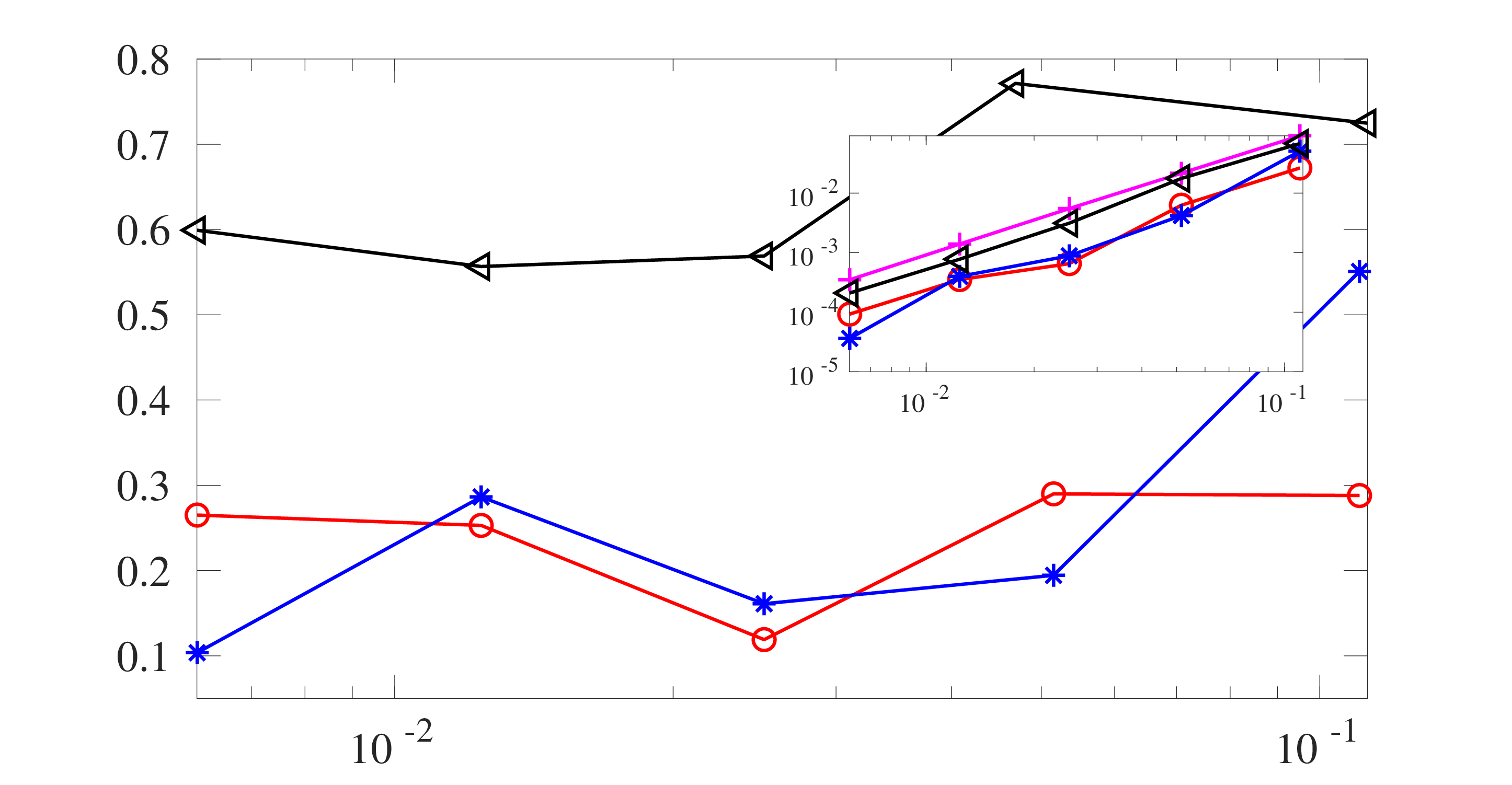}}
    \put(9,0){\includegraphics[width=0.49\linewidth]{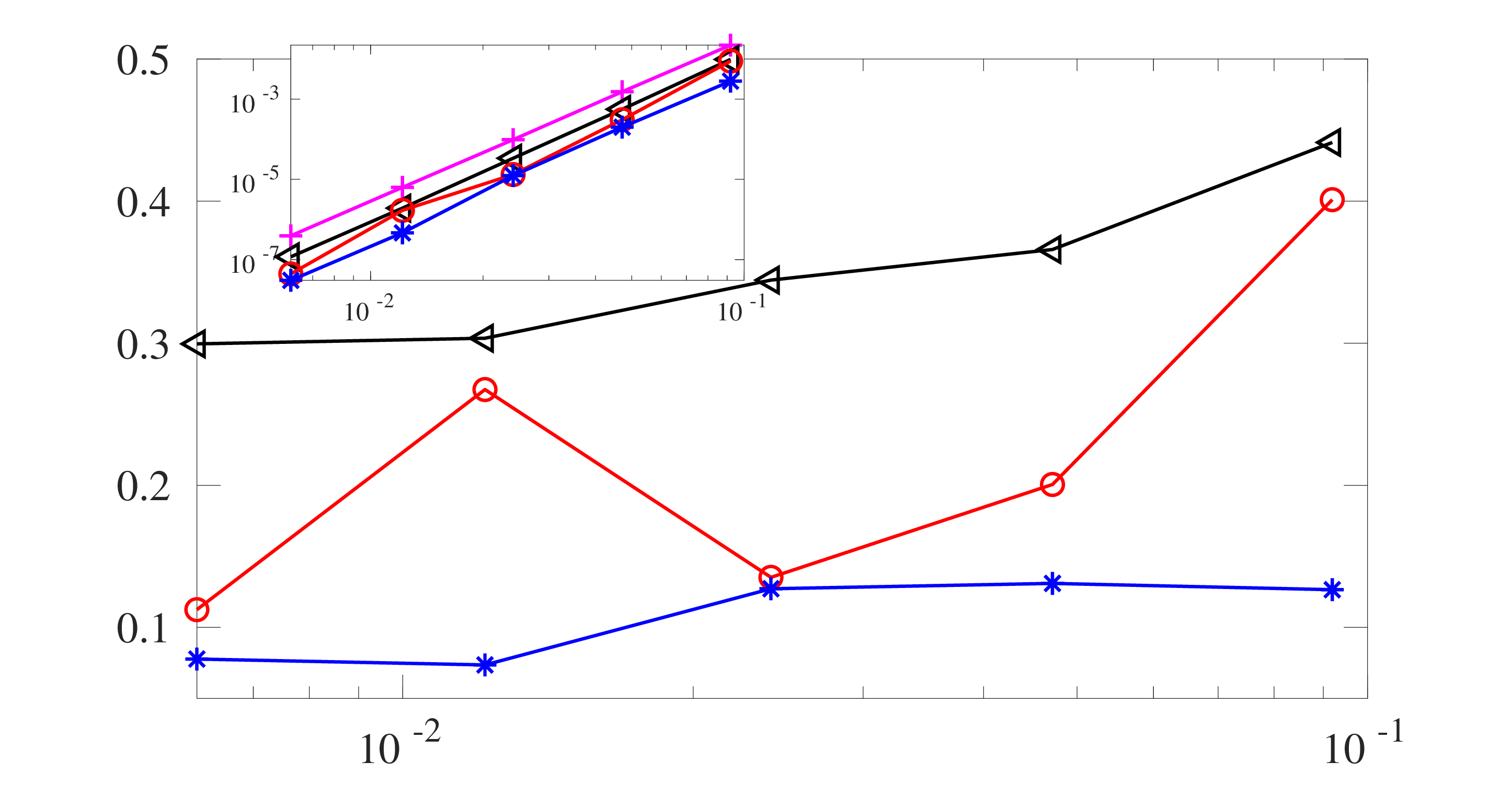}}
    \put(-0.1,2.25){\rotatebox{90}{$\mathcal{R}$}}
   \put(8.9,2.25){\rotatebox{90}{$\mathcal{R}$}}
    \put(4.5,-0.25){$s$}
    \put(13.5,-0.25){$s$}
\end{picture}
    \end{center}
    \caption{Ratio $\mathcal{R}$ of errors arising from explicit and compact approximations to Poisson's equation \eqref{eq::Poisson_equation} with right-hand side given by \eqref{eq:Poisson_function}. Left Panel, $2^{nd}$ order approximations: black $\triangleleft$, $2Q_i-1=5$ (scheme b); red $\circ$ $2Q_i-1=9$ (scheme c); blue $*$, $2Q_i-1=13$ (scheme d). Right Panel, $4^{th}$ order approximations; black $\triangleleft$, $2Q_i-1=9$ (scheme f); red $\circ$ $2Q_i-1=13$ (scheme g); blue $*$, $2Q_i-1=17$ (scheme h). Insets show convergence of $L_2-$ norm with magenta line corresponding to explicit LABFM.}
    \label{fig:poisson_eq}
    \end{figure}
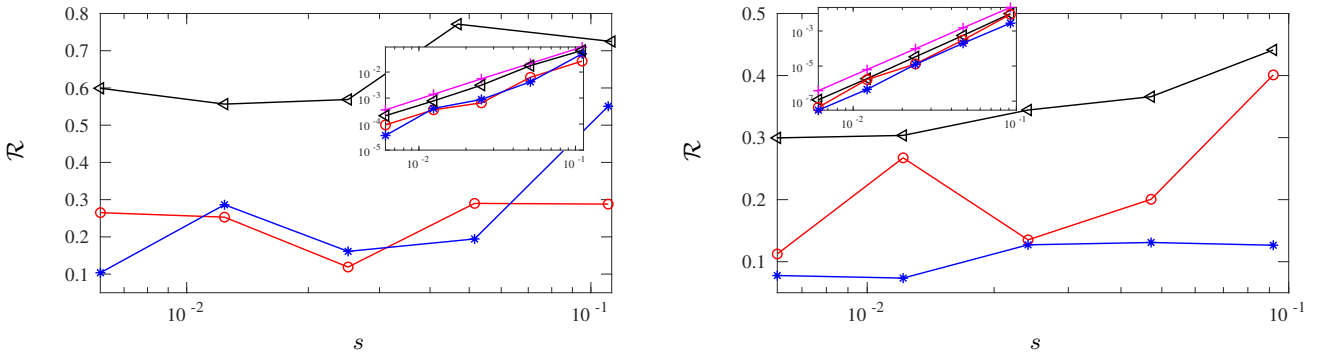

Fig.\ref{fig:poisson_eq} compares the solutions produced by explicit operators to those from various compact stencils. For $2^{nd}$ order operators scheme (b) with only five nodes in the compact stencil has a small but consistent error reduction of $20-40\%$. Larger stencils on the other hand have a consistently much larger reduction in error up to an order of magnitude. Similar behaviour is observed for $4^{th}$ order operators, though here there is both greater improvement for the smallest implicit stencil (always a reduction in error of more than a half) and the largest implicit stencil, which performs best for all nodal spacings. When using scheme (h), an error reduction of an order of magnitude is consistently observed, even for the relatively simple test function given by \eqref{eq:Poisson_function}. This demonstrates the potential for compact LABFM to provide significant improvements in the accuracy of solutions to PDEs in complex geometries, for negligible increase in computational cost.

\section{Conclusions}
\label{sec:Conclusions}
We have presented a new compact meshless method which approximates operators using Lele-type schemes. Our approach optimises the resolving power at each node, enabling significant improvements in this metric across all of wavenumber space. The resolving efficiency of the compact operators has been shown to compare especially favourably to their explicit counterparts, with in general larger implicit stencils resulting in better resolution characteristics. Our framework allows for the construction of $4^{th}$ order consistent operators with spectral-like accuracy. The new method shows improved approximation of operators in convergence tests, and offers up to an order of magnitude reduction in error when simulating canonical PDE systems.

Our method of choosing the implicit stencil and coefficients encourages diagonal dominance and ensures the global linear system remains well-conditioned whilst maximising resolving power, however it may not yield the optimal resolution characteristics. No investigations were performed into varying the prescribed form of the coefficients. It is possible that a different form, for example similar to a Wendland C2 kernel or a different exponent will lead to improved accuracy. It is also possible that choosing the nodes in the implicit stencil in a different manner may improve the resolving power of the scheme. Such an investigation is out of the scope of this proof of concept study, however our group remains interested in exploring this topic further.

Although we have focused on applying our methodology to generate a compact formulation of LABFM, in theory our approach is transferable to other meshless methods. For example, one could generate a Lele-type scheme for compact RBF-FDs with resolving power optimisation. Our framework provides a rigorous footing for such schemes, enabling the next generation of meshless methods to aim for high-order, spectral-like resolving power of (for example) turbulent flow simulations.

\begin{acknowledgements}
This work was partially funded by the Engineering and Physical Sciences Research Council (EPSRC) grant EP/W005247/2. JK is supported by the Royal Society via a University Research Fellowship (URF\textbackslash R1\textbackslash 221290). We are grateful for assistance given by Research IT and the use of the Computational Shared Facility at the University of Manchester. 
\end{acknowledgements}

\bibliographystyle{elsarticle-num} 
  \bibliography{HB_bib}
\end{document}